\documentclass[11pt,a4paper,reqno]{amsart}
\usepackage[a4paper,margin=29mm]{geometry}
\usepackage{amsmath,amssymb,amsthm,mathtools}
\usepackage{microtype}
\usepackage{needspace}
\usepackage{xcolor}
\usepackage[unicode,breaklinks=true,colorlinks=false,
  pdfborder={0 0 0.5},
  linkbordercolor={0.75 0.05 0.05},
  citebordercolor={0 0.5 0.15},
  urlbordercolor={0 0.3 0.75}]{hyperref}
\hypersetup{
  pdftitle={Linearized contact homology and closed Reeb orbits on prequantization bundles},
  pdfauthor={Dong Xie},
  pdfsubject={Linearized contact homology and multiplicity of closed Reeb orbits},
  pdfkeywords={Linearized contact homology, prequantization bundles, closed Reeb orbits, Morse-Bott theory, index iteration}}
\newtheorem{theorem}{Theorem}[section]
\newtheorem{lemma}[theorem]{Lemma}
\newtheorem{proposition}[theorem]{Proposition}
\newtheorem{corollary}[theorem]{Corollary}
\newtheorem*{conjecture*}{Conjecture}
\numberwithin{equation}{section}
\theoremstyle{remark}
\newtheorem{remark}[theorem]{Remark}
\newcommand{\LCH}{\mathrm{LCH}}
\newcommand{\LCC}{\mathrm{LCC}}
\renewcommand{\theenumi}{\roman{enumi}}
\renewcommand{\labelenumi}{\textup{(\theenumi)}}
\allowdisplaybreaks[2]
\title[Linearized contact homology and Reeb orbits]{Linearized contact homology and closed Reeb orbits on prequantization bundles}
\author[Dong Xie]{Dong Xie}
\address{School of Mathematical Sciences, Nankai University,
Tianjin 300071, China; Centre de recherches math\'ematiques (CRM),
Universit\'e de Montr\'eal, Montr\'eal, QC H3C 3J7, Canada}
\email{\href{mailto:dong.xie@nankai.edu.cn}{dong.xie@nankai.edu.cn}}
\date{}
\begin{document}
\begin{abstract}
We study closed Reeb orbits on prequantization
circle bundles $(M,\xi)\to(B,\omega)$ with
$c_1(\xi)|_{\pi_2(M)}=0$.
When $\omega|_{\pi_2(B)}\ne0$, we construct augmentations and compute
contractible linearized contact homology as a direct sum of shifted
copies of $H_*(B;\mathbb Q)$, under a condition on the rational
Betti numbers of $B$.
The construction uses Morse--Bott methods and requires neither a
filling nor index-admissibility.
When $c_1(TB)|_{\pi_2(B)}\ne0$, index recurrence and Morse inequalities
yield multiplicity and non-hyperbolicity bounds for nondegenerate
contact forms under additional index assumptions.
If $c_1(TB)|_{\pi_2(B)}=0$ and either
$H_{\mathrm{odd}}(B;\mathbb Q)=0$ or $\omega|_{\pi_2(B)}=0$,
every nondegenerate supporting contact form has infinitely many
simple closed Reeb orbits.
We also establish multiplicity and non-hyperbolicity results for
smooth prequantizations of symplectic orbifolds admitting Hamiltonian
circle actions with isolated fixed points.
\end{abstract}
\maketitle
\markboth{LINEARIZED CONTACT HOMOLOGY}{CLOSED REEB ORBITS ON PREQUANTIZATION BUNDLES}
\begingroup
\hypersetup{hidelinks}
% The contents entry links to the proof, including its theorem number.
\RenewDocumentCommand{\ref}{s m}{\getrefnumber{#2}}
\tableofcontents
\endgroup

\section{Introduction}

\subsection{Background}
\label{subsec:intro-background}

We study linearized contact homology and the multiplicity of
closed Reeb orbits on prequantization $S^1$-bundles.

On a closed contact manifold $(M^{2n+1},\xi)$, a contact form
$\alpha$ supporting $\xi$ determines the Reeb vector field
$R_\alpha$ by
\[
 \alpha(R_\alpha)=1,
 \qquad
 d\alpha(R_\alpha,\cdot)=0.
\]
We denote its flow by $\varphi_{R_\alpha}^t$.
A \emph{closed Reeb orbit} of period $T>0$ is a smooth map
$\gamma:\mathbb R/T\mathbb Z\to M$ satisfying
$\dot\gamma(t)=R_\alpha(\gamma(t))$.
For each integer $k\ge1$, the \emph{$k$th iterate} $\gamma^k$
is the $k$-fold cover of $\gamma$ and has period $kT$.
A closed Reeb orbit is \emph{simple} if it is not a
nontrivial multiple cover of another closed Reeb orbit.
Two closed Reeb orbits
$\gamma_1:\mathbb R/T_1\mathbb Z\to M$ and
$\gamma_2:\mathbb R/T_2\mathbb Z\to M$
are \emph{geometrically distinct} if their images are different
subsets of $M$.

A closed Reeb orbit is \emph{contractible} if it is homotopic
to a constant loop in $M$.
An iterate $\gamma^q$ of a simple closed Reeb orbit $\gamma$
is called the \emph{least contractible iterate} of $\gamma$,
or a \emph{least contractible closed Reeb orbit}, if $q$ is
the smallest positive integer such that $\gamma^q$ is contractible.

The \emph{linearized return map} of a closed Reeb orbit
$\gamma$ of period $T$ is
\[
 d\varphi_{R_\alpha}^{T}\big|_{\xi_{\gamma(0)}}:
 \xi_{\gamma(0)}\longrightarrow\xi_{\gamma(0)}.
\]
The orbit $\gamma$ is \emph{nondegenerate} if this map has
no eigenvalue $1$, \emph{elliptic} if all its eigenvalues
lie on the unit circle in $\mathbb C$, and \emph{hyperbolic}
if none of its eigenvalues lies on the unit circle.
The contact form $\alpha$ is \emph{nondegenerate} if every
closed Reeb orbit, including all its iterates, is nondegenerate.

The existence and multiplicity of closed Reeb orbits
have long been central questions in contact dynamics.
Linearized contact homology and symplectic homology, including
its positive $S^1$-equivariant version, provide tools for studying
these questions. Contact homology was introduced by Eliashberg,
Givental, and Hofer as part of symplectic field theory
\cite{EGH00}. The $S^1$-equivariant version of symplectic homology
was introduced by Viterbo \cite[Section~5]{Vit99} and developed
by Bourgeois and Oancea \cite{BO13,BO17}; see also \cite{Gu17}.
For suitable Liouville fillings, Bourgeois and Oancea \cite{BO17}
identify positive $S^1$-equivariant symplectic homology over
$\mathbb Q$ with linearized contact homology defined using the filling.

A closed contact manifold $(M,\xi)$ is called
\emph{homologically unbounded} if, for some augmentation, the
ranks of its linearized contact homology are unbounded along
a sequence of degrees tending to $+\infty$ or $-\infty$.
For such a manifold with $c_1(\xi)=0$, Hryniewicz and Macarini
\cite{HM15} proved that every contact form
supporting $\xi$ has infinitely many geometrically distinct
simple closed Reeb orbits.
McLean \cite{McL12} proved a corresponding result
using symplectic homology of Liouville domains with trivialized
canonical bundle.
For a survey of multiplicity results for Reeb flows, see
\cite{Mac25}.

Among \emph{homologically bounded} contact manifolds, prequantization
$S^1$-bundles provide an interesting and natural class of examples.
(Here infinite rank in finitely many degrees is allowed.)
Let $(B^{2n},\omega)$ be a closed connected integral symplectic
manifold. A \emph{prequantization $S^1$-bundle}
$\pi:(M^{2n+1},\xi)\to(B^{2n},\omega)$ is a principal
$S^1$-bundle with first Chern class $-[\omega]$, equipped
with a connection form $\beta$ satisfying $\ker\beta=\xi$
and $d\beta=\pi^*\omega$ \cite{BoW58}.
With the convention $S^1=\mathbb R/\mathbb Z$, the fibers
are the simple closed Reeb orbits of $\beta$, each of period $1$.
A closed Reeb orbit whose image is a fiber is called a
\emph{fiber orbit}.

A classical example of a prequantization $S^1$-bundle is the
Hopf fibration $S^{2n+1}\to\mathbb{CP}^n$.
Irrational ellipsoids give nondegenerate contact forms supporting
the standard contact structure on $S^{2n+1}$, with exactly
$n+1$ geometrically distinct simple closed Reeb orbits,
all elliptic; see \cite[Section~1.1]{AbM24}.
Another class of examples consists of the unit cotangent bundles
$S^*N$ of compact rank-one symmetric spaces $N$ of dimension
at least two. The geodesic flow of the standard metric generates
a free $S^1$-action on $S^*N$. The quotient $B=S^*N/S^1$ is
the space of oriented geodesics of this metric, and $S^*N\to B$
is a prequantization $S^1$-bundle.
The Katok--Ziller Finsler metrics on $N$ induce nondegenerate
contact forms supporting the same contact structure on $S^*N$,
with finitely many geometrically distinct simple closed Reeb
orbits, all elliptic \cite{Zil83}.
In both examples, the number of these orbits equals the total
rational Betti number of the corresponding symplectic base manifolds; see
\cite[Section~2.1]{GGM18}.

For convenience, we write
\[
 b_j(B)=\dim_{\mathbb Q}H_j(B;\mathbb Q),
 \qquad
 b(B)=\sum_{j=0}^{2n}b_j(B),
\]
with $b_j(B)=0$ for $j<0$ or $j>2n$.
Regarding the multiplicity of closed Reeb orbits on
prequantization $S^1$-bundles, the following conjecture is
stated in \cite[Conjecture~1.1]{LX26}.
For related questions and results, see
\cite{AbM17} and \cite{AbM26}.

\begin{conjecture*}
Let $\pi:(M^{2n+1},\xi)\to(B^{2n},\omega)$, $n\ge1$, be a
prequantization circle bundle over a closed connected integral
symplectic manifold. Every contact form supporting $\xi$ has
at least $b(B)$ geometrically distinct simple closed Reeb orbits.
\end{conjecture*}

The conjecture holds in dimension three; see also
\cite[Remark~1.2]{LX26}. If $B=S^2$,
then $M$ is a lens space (possibly $S^3$) and $c_1(\xi)$
is torsion. By \cite{CGHHL26}, every contact form supporting
$\xi$ has either two or infinitely many geometrically distinct
simple closed Reeb orbits. Since $b(S^2)=2$, the conjecture
follows in this case. If $B$ has positive genus, its universal
cover is contractible. Thus $B$ is aspherical, that is,
$\pi_j(B)=0$ for every $j\ge2$. Moreover, $c_1(TB)$ is a
rational multiple of $[\omega]$, since $H^2(B;\mathbb Q)$
is one-dimensional. The identity $d\beta=\pi^*\omega$ gives
$\pi^*[\omega]=0$, and hence $c_1(\xi)=\pi^*c_1(TB)=0$
in $H^2(M;\mathbb Q)$. Therefore
\cite[Theorem~6.4]{GGM18} applies and gives infinitely many
geometrically distinct simple closed Reeb orbits for every
contact form supporting $\xi$.

In higher dimensions, the conjecture remains open for arbitrary
supporting contact forms.
On the standard contact sphere $S^{2n+1}$,
\c{C}ineli, Ginzburg, and G\"urel \cite[Theorem~A]{CGG24}
proved that every dynamically convex contact form has at least
$n+1$ geometrically distinct simple closed Reeb orbits,
allowing degenerate closed Reeb orbits.
For nondegenerate contact forms supporting the standard contact
structure on $S^{2n+1}$, \c{C}ineli, Ginzburg, and G\"urel
\cite[Theorem~B]{CGG26} proved that a closed Reeb orbit with
negative mean index forces infinitely many geometrically
distinct simple closed Reeb orbits.
For further multiplicity and stability results for closed
characteristics on star-shaped hypersurfaces, see
\cite{DLLW24,DLLQW25,DuQ26b}.

For the unit cotangent bundles $S^*N$ described above,
a Finsler metric on $N$ induces a supporting
contact form whose Reeb flow corresponds to the geodesic
flow of that metric.
For multiplicity and stability results for closed geodesics
on these spaces, see \cite{DLW16,DuX23,DuQ25,DuX26}.
For related results on restricted
contact type hypersurfaces in $T^*S^m$ under dynamical convexity,
see \cite{DuQ26a}.

For general prequantization bundles, the known multiplicity
results involve additional geometric or index hypotheses.
When $\omega|_{\pi_2(B)}=0$,
infinitely many geometrically distinct simple closed Reeb orbits
are known to exist under additional hypotheses
\cite{GGM15,GS18}.
When $\omega|_{\pi_2(B)}\ne0$, Abreu and Macarini
proved the lower bound $b(B)$ for nondegenerate
dynamically convex supporting forms under the hypotheses of
\cite[Theorem~2.2]{AbM17}.
For multiplicity results under pinching and other geometric
hypotheses, see \cite{AGH18,Ker17}.

Under additional geometric and index hypotheses,
Ginzburg, G\"urel, and Macarini
\cite[Theorems~2.1 and~6.2]{GGM18} obtained lower bounds for
geometrically distinct contractible closed Reeb orbits of
nondegenerate contact forms, assuming that all contractible
closed Reeb orbits have nonzero mean indices of the same sign.
These lower bounds equal $b(B)$ when
$H_{\mathrm{odd}}(B;\mathbb Q)=0$.
In the symplectically aspherical filling case,
Duan, Liu, Long, and Wang \cite[Theorem~1.3]{DLLW24}
obtained the same lower bounds without requiring the
nonzero mean indices to have the same sign.

In \cite{AbM24}, Abreu and Macarini studied nondegenerate
lacunary contact forms on prequantization bundles and certain
smooth prequantizations of symplectic orbifolds.
Here lacunarity means that all contractible closed Reeb orbits
have Conley--Zehnder indices of the same parity.
Under suitable geometric and index hypotheses, they obtained
exact multiplicity results.
For prequantization bundles, these exact multiplicity results
allow a weaker restriction on the minimal Chern number of the
base than that in \cite{GGM18}.
In the non-filling case, they also allow the first Chern class
of the contact structure to be torsion rather than zero.

In \cite{AbM26}, they established multiplicity and
stability results for nondegenerate contact forms on contact
manifolds with periodic positive $S^1$-equivariant symplectic
homology. Their applications include prequantization bundles
and certain smooth prequantizations of symplectic orbifolds.
For positive monotone bases with
$H_{\mathrm{odd}}(B;\mathbb Q)=0$, they obtain at least
$b(B)$ geometrically distinct simple closed Reeb orbits
under their remaining hypotheses, assuming that every closed
orbit has nonzero mean index and that no good closed orbit
has Conley--Zehnder index $\pm p$ for some integer $p\ge0$
with $p\not\equiv n\pmod2$.
These results further relax the index restrictions in the
corresponding multiplicity results of \cite{GGM18,AbM24}.
The proofs of these prequantization results still rely largely
on suitable symplectic fillings or on index-admissibility
in the non-filling setting.

Liu \cite[Theorem~1.6]{LX26} obtained multiplicity bounds
for nondegenerate supporting contact forms when $n$ is even,
under related geometric and topological hypotheses and
different index assumptions. Under the hypotheses of
\cite[Theorem~1.11]{LX26}, he also proved the lower bound
$b(B)$ in dimension five.

For prequantization bundles, the periodic Reeb flow of the
connection form provides the starting point for Morse--Bott calculations.
For each $m\ge1$, the $m$th iterates of the simple closed
Reeb orbits of $\beta$ form a family parametrized by $B$.
Bourgeois's Morse--Bott approach to contact homology uses
Morse functions on these orbit spaces \cite{Bou}; see also
\cite{BO09} for Morse--Bott methods in symplectic homology
and \cite{Nel20} for computations of cylindrical contact
homology on three-dimensional prequantization bundles.

Under the hypotheses of \cite[Proposition~3.1]{GGM18},
positive $S^1$-equivariant symplectic homology in the
contractible class is a direct sum of shifted copies of
$H_*(B;\mathbb Q)$.
A corresponding computation with Novikov coefficients is
given in \cite[Proposition~2.1]{AbM24}.
For positive monotone bases, Abreu and Macarini also obtain
multiplicity bounds without assuming
$H_{\mathrm{odd}}(B;\mathbb Q)=0$, provided that the
analogous homology formula with all free homotopy classes
included holds \cite[equation~(8.1) and Remark~8.1]{AbM26}.

Linearized contact homology depends on a choice of augmentation
of the contact algebra.
Pardon constructed the contact algebra using virtual
fundamental cycles \cite{Par16,Pardon}; see \cite{BaH23} for an
alternative construction.
Chaidez proved that the set of weak equivalence classes
of augmentations is a contact invariant and that the linearized
homology associated with each class is well defined up to isomorphism
\cite{Chaidez}.

Ginzburg, G\"urel, and Macarini
\cite[Section~6.1]{GGM18}
construct an augmentation for any nondegenerate supporting
contact form by pulling back the trivial augmentation of a nondegenerate
index-admissible contact form.
Here index-admissibility excludes contractible closed orbits with
Conley--Zehnder indices $1-n$, $2-n$, or $3-n$
\cite[Section~2.1]{GGM18}.

For prequantization bundles, they also pointed out that
the Morse--Bott calculation needed for the homology formula remained
to be justified within Pardon's framework
\cite[Section~6.2]{GGM18}.
Zhang subsequently developed Morse--Bott contact homology
using Pardon's implicit atlases \cite{Zhang19}.

In this paper, we use Morse--Bott methods to construct
augmentations and compute contractible linearized contact homology
for prequantization circle bundles under the topological hypotheses
stated below.
The construction requires neither a filling nor index-admissibility.
We then establish multiplicity and non-hyperbolicity results for
closed Reeb orbits on these bundles and on smooth prequantizations
of symplectic orbifolds.

\subsection{Main results}
\label{subsec:intro-main-results}

Let $\pi:(M^{2n+1},\xi)\to(B^{2n},\omega)$, $n\ge1$,
be a prequantization circle bundle over a closed connected integral
symplectic manifold.
In all our main results we assume that $c_1(\xi)|_{\pi_2(M)}=0$.

We denote by $\chi(B)$ the Euler characteristic of $B$ and by $c_B$
its minimal Chern number, defined by
\begin{equation}\label{eq:minimal-chern-number}
 c_B :=
 \begin{cases}
  \inf\{k\in\mathbb{Z}_{>0} \mid
        \langle c_1(TB),A\rangle=k \text{ for some } A\in\pi_2(B)\},
        & \text{if } c_1(TB)|_{\pi_2(B)}\ne0,\\[2pt]
  0,    & \text{if } c_1(TB)|_{\pi_2(B)}=0.
 \end{cases}
\end{equation}

When $c_B>0$, Lemma~\ref{lem:topology} implies that
$\omega|_{\pi_2(B)}\ne0$ and that $(B,\omega)$ is
\emph{spherically monotone},~i.e.,
\[
 c_1(TB)|_{\pi_2(B)}=\lambda[\omega]|_{\pi_2(B)}
\]
for some nonzero $\lambda\in\mathbb R$. We call $B$
\emph{spherically positive monotone} if $\lambda>0$ and
\emph{spherically negative monotone} if $\lambda<0$.
We set $\sigma=\operatorname{sign}(\lambda)$ when $c_B>0$ and
$\sigma=0$ when $c_B=0$.

Let $\gamma$ be a closed Reeb orbit and let $\tau$ be a symplectic
trivialization of $\xi$ along $\gamma$. We write
$\widehat\mu(\gamma;\tau)$ for its mean index and, when $\gamma$ is
nondegenerate, $\mu(\gamma;\tau)$ for its Conley--Zehnder index.
If $\gamma$ is contractible in $M$, we choose $\tau$ to extend over
a capping disk. Since $c_1(\xi)|_{\pi_2(M)}=0$, the resulting indices
are independent of the choice of capping disk, and we write
$\widehat\mu(\gamma)$ and $\mu(\gamma)$, respectively
(see Section~\ref{subsec:index-iteration}).

{Let $\alpha$ be a nondegenerate contact form supporting
$\xi$, and let $\gamma$ be a simple closed Reeb orbit of $\alpha$.}
Choose a symplectic trivialization $\tau$
of $\xi$ along $\gamma$. The iterate $\gamma^k$ is \emph{good} if
\[
 \mu(\gamma^k;\tau^k)\equiv\mu(\gamma;\tau)\pmod 2
\]
and \emph{bad} otherwise, where $\tau^k$ is the trivialization
along $\gamma^k$ induced by $\tau$. This condition is independent
of the choice of $\tau$.

{Let $\mathcal C(\alpha)$ denote the contact algebra of
$\alpha$ over $\mathbb Q$,
generated by the good closed Reeb orbits; see
\cite[Section~1.2]{Pardon}.
Its \emph{contractible quotient} $(\mathcal B(\alpha),d)$ is the
quotient by the ideal generated by the noncontractible orbits.}
A \emph{graded augmentation} of
$\mathcal B(\alpha)$ is a degree-preserving unital algebra map
$\bar\varepsilon_\alpha:\mathcal B(\alpha)\to\mathbb Q$ satisfying
$\bar\varepsilon_\alpha\circ d=0$, where $\mathbb Q$ is concentrated
in degree zero. Its pullback $\varepsilon_\alpha$ to
$\mathcal C(\alpha)$ vanishes on every noncontractible generator.

A generator $\gamma$ of $\mathcal B(\alpha)$ has degree
$|\gamma|=\mu(\gamma)+n-2$. In the linearized complex at
$\bar\varepsilon_\alpha$, we grade the corresponding generator by
$\mu(\gamma)$ and denote {its degree-$k$ homology} by
$\LCH^0_k(M,\alpha;\varepsilon_\alpha)$
(see Section~\ref{subsec:contact-algebras}).
{The following theorem gives a Morse--Bott computation
of the linearized contact homology in the trivial free homotopy
class of $M$.}

\begin{theorem}\label{thm:homology}
Let $\pi:(M^{2n+1},\xi)\to(B^{2n},\omega)$, $n\ge1$, be a
prequantization circle bundle over a closed connected integral
symplectic manifold, with Euler class $-[\omega]$. Assume that
\[
 c_1(\xi)|_{\pi_2(M)}=0,
 \qquad
 \omega|_{\pi_2(B)}\ne0,
\]
and that
\[
 b_j(B)b_{j+2\sigma c_B-1}(B)=0
 \qquad\text{for every }j\in\mathbb Z.
\]
For every nondegenerate contact form $\alpha$ supporting $\xi$,
there exists a graded augmentation $\bar\varepsilon_\alpha$ of
$\mathcal B(\alpha)$ such that, for every $k\in\mathbb Z$, there is
an isomorphism of rational vector spaces
\begin{equation}\label{eq:BWpattern}
 \LCH^0_k(M,\alpha;\varepsilon_\alpha)
 \cong \bigoplus_{m\ge1} H_{k-2\sigma m c_B+n}(B;\mathbb Q).
\end{equation}
\end{theorem}

\begin{remark}\label{rem:zero-chern-rank}
Theorem~\ref{thm:homology} also applies when $c_B=0$.
In this case, taking $k=-n$ in \eqref{eq:BWpattern} gives
\[
 \LCH^0_{-n}(M,\alpha;\varepsilon_\alpha)
 \cong \bigoplus_{m\ge1} H_0(B;\mathbb Q)
 \cong \bigoplus_{m\ge1}\mathbb Q.
\]
We use the resulting infinite rank in the proof of
Theorem~\ref{thm:zero-chern}\textup{(ii)}.
\end{remark}

\begin{remark}\label{rem:betti-hypothesis}
The positive $S^1$-equivariant symplectic homology computations
in \cite[Proposition~3.1]{GGM18} and
\cite[Proposition~2.1]{AbM24} assume either
\[
 H_{\mathrm{odd}}(B;\mathbb Q)=0
 \qquad\text{or}\qquad c_B>n,
\]
in addition to their respective filling or index-admissibility
hypotheses.
Each condition implies
\[
 b_j(B)b_{j+2\sigma c_B-1}(B)=0
 \qquad\text{for every }j\in\mathbb Z.
\]
Indeed, the two degrees have opposite parity, so the first
condition forces one of the Betti numbers to vanish.
If $c_B>n$, then $|2\sigma c_B-1|>2n$, and the two degrees
cannot both lie in $[0,2n]$.
\end{remark}

\begin{remark}\label{rem:quartic-k3}
A concrete example of the case $c_B=0$ is provided by the
Fermat quartic
\[
 B=\bigl\{[z_0:z_1:z_2:z_3]\in\mathbb{CP}^3
 \,\big|\,
 z_0^4+z_1^4+z_2^4+z_3^4=0\bigr\}.
\]
This is a K3 surface, that is, a compact simply connected
complex surface with trivial canonical bundle.
Equip $B$ with the restriction $\omega$ of the Fubini--Study
form, normalized so that $[\omega]$ is the integral hyperplane
class, and let $\pi:(M,\xi)\to(B,\omega)$ be the associated
prequantization circle bundle.
Then $c_1(TB)=0$, $H_{\mathrm{odd}}(B;\mathbb Q)=0$, and
$\omega|_{\pi_2(B)}\ne0$.
Thus Theorem~\ref{thm:homology} applies with $n=2$ and $c_B=0$,
and Remark~\ref{rem:zero-chern-rank} gives infinite rank in
degree $-2$.

For a Liouville filling of this prequantization bundle,
Kwon and van Koert \cite[Remark~C.1(i)]{KvK16} also observed
that the positive $S^1$-equivariant symplectic homology group
in degree zero is infinite-dimensional.
\end{remark}

For a nondegenerate contact form $\alpha$ supporting $\xi$,
let $r_\alpha$ denote the number of geometrically distinct
contractible closed Reeb orbits of $\alpha$, and let
$r_\alpha^{\mathrm{nh}}$ denote the number of non-hyperbolic
orbits among them.

\begin{theorem}\label{thm:main}
Let $\pi:(M^{2n+1},\xi)\to(B^{2n},\omega)$, $n\ge2$, be a
prequantization circle bundle over a closed connected integral
symplectic manifold, with Euler class $-[\omega]$. Assume that
$c_1(\xi)|_{\pi_2(M)}=0$, $c_B>0$, and
\[
 b_j(B)b_{j+2\sigma c_B-1}(B)=0
 \qquad\text{for every }j\in\mathbb Z.
\]

Let $\alpha$ be a nondegenerate contact form supporting $\xi$,
and fix an integer $p\ge0$. Assume that every good contractible
closed Reeb orbit $\gamma$ of $\alpha$ satisfies
$\widehat\mu(\gamma)\ne0$ and
\[
 \mu(\gamma)\notin\{-p,p\}.
\]
Then the following hold.
\begin{enumerate}
\item
If $p\not\equiv n\pmod2$, then
\begin{equation}\label{eq:main-bound}
 r_\alpha\ge
 \begin{cases}
  \displaystyle\max\left\{
   \chi(B)+2b_n(B),\,
   \frac{b_n(B)-\chi(B)}{2}+1
  \right\},&p=0,\\[6pt]
  \displaystyle\max\left\{
   \chi(B)+4b_{n-p}(B),\,
   \left\lceil-\frac{\chi(B)}{2}\right\rceil+b_{n-p}(B)+1
  \right\},&p>0.
 \end{cases}
\end{equation}

\item
Suppose $p\equiv n\pmod2$. If $p=0$ or $n=2$, then
\begin{equation}\label{eq:main-same-zero}
 r_\alpha\ge
 \max\left\{
  \frac{\chi(B)+b_n(B)}{2}+1,\,
  -\chi(B)+2b_n(B)
 \right\}.
\end{equation}
If $n=p=2$, we also have
\begin{equation}\label{eq:main-same-positive-lowdim}
 r_\alpha\ge4-\chi(B).
\end{equation}

If $p=1$, then
\begin{equation}\label{eq:main-odd-one}
 r_\alpha\ge
 \max\left\{
  \chi(B)+2b_n(B),\,
  \left\lceil\frac{\chi(B)}2\right\rceil+b_{n-1}(B)+1,\,
  -\chi(B)+4b_{n-1}(B)
 \right\}.
\end{equation}

If $p\ge2$ and $n\ge3$, then
\begin{equation}\label{eq:main-same-positive}
 r_\alpha\ge
 \max\left\{
  \left\lceil\frac{\chi(B)}2\right\rceil+b_{n-p}(B)+1,\,
  -\chi(B)+4b_{n-p}(B)
 \right\}.
\end{equation}
\end{enumerate}
\end{theorem}

\begin{remark}\label{rem:same-parity-special}
If $H_{\mathrm{odd}}(B;\mathbb Q)=0$, every maximum in
Theorem~\ref{thm:main} equals its first entry.
The bounds in \eqref{eq:main-bound} and \eqref{eq:main-odd-one},
and in \eqref{eq:main-same-zero} when $n=2$, then give
\[
 r_\alpha\ge\chi(B)=b(B).
\]
The additional bound \eqref{eq:main-same-positive-lowdim}
is then redundant.
\end{remark}

\begin{remark}\label{rem:counting-scope}
Under their respective hypotheses, Ginzburg, G\"urel,
and Macarini \cite[Theorems~2.1 and~6.2]{GGM18} obtained
\[
 r_\alpha\ge
 \begin{cases}
  \chi(B)+2b_n(B),&n\text{ odd},\\
  \chi(B)+4b_{n-1}(B),&n\text{ even}.
 \end{cases}
\]
Theorem~\ref{thm:main}\textup{(i)} recovers these bounds with
$p=0$ and $p=1$, respectively.
Following Abreu and Macarini \cite[Theorem~5]{AbM26},
we exclude good contractible closed Reeb orbits with indices
$\pm p$. Our proof adapts their counting argument in
\cite[Section~6]{AbM26}, using the homology formula of
Theorem~\ref{thm:homology}.

For positive monotone $B$ with
$H_{\mathrm{odd}}(B;\mathbb Q)=0$, Abreu and Macarini
\cite[Theorems~1 and~2, case~\textup{(PM)}]{AbM26} obtain at least
$b(B)$ simple closed Reeb orbits under their filling or non-filling
hypotheses. In their non-filling case, $\alpha$ must be
index-admissible. Theorem~\ref{thm:main}\textup{(i)} gives
$r_\alpha\ge b(B)$ when $H_{\mathrm{odd}}(B;\mathbb Q)=0$ for
both $\sigma=1$ and $\sigma=-1$.

For even $n$, Liu \cite[Theorem~1.6]{LX26} obtained
$r_\alpha\ge(\chi(B)+b_n(B))/2+1$ when $\widehat\mu(\gamma)\ne0$
for every contractible closed Reeb orbit $\gamma$ and
$\mu(\gamma)\ne0$ for every good contractible closed Reeb orbit
$\gamma$. Liu also assumes $c_B>n/2$, either
$H_{\mathrm{odd}}(B;\mathbb Q)=0$ or $c_B>n$, and the additional
hypotheses stated there in the filling and non-filling cases.
Theorem~\ref{thm:main}\textup{(ii)} recovers this estimate for
$p=0$ and also gives estimates for $p>0$ with $p\equiv n\pmod2$.

The homology construction in Theorem~\ref{thm:homology} allows us
to obtain the bounds in Theorem~\ref{thm:main} under
$c_1(\xi)|_{\pi_2(M)}=0$, $c_B>0$, and
$b_j(B)b_{j+2\sigma c_B-1}(B)=0$ for every $j\in\mathbb Z$,
without a filling or index-admissibility of $\alpha$.
Theorem~6.2 of \cite{GGM18} also requires neither a filling nor
index-admissibility of $\alpha$,
but assumes $c_1(\xi)=0$ in $H^2(M;\mathbb Z)$,
$H_{\mathrm{odd}}(B;\mathbb Q)=0$ or $c_B>n$, and
$c_B>n/2$ for $\sigma=1$ or $c_B\ge n$ for $\sigma=-1$.
\end{remark}

\begin{theorem}\label{thm:main-nonhyperbolic}
Assume the hypotheses of Theorem~\ref{thm:main}, and suppose
that $r_\alpha<\infty$.
If $p=0$, then
\begin{equation}\label{eq:main-nonhyperbolic}
 r_\alpha^{\mathrm{nh}}\ge
 \max\left\{0,\,b_n(B)-(-1)^n\chi(B)\right\}.
\end{equation}
{For $p>0$, set}
\begin{equation}\label{eq:main-nonhyperbolic-alternating}
 A_p=-2(-1)^{n+p}\sum_{j=0}^{n-p-1}(-1)^j b_j(B),
\end{equation}
with empty sums understood to be zero.
\begin{enumerate}
\item
Suppose $n=2$. If $p$ is even, then
\begin{equation}\label{eq:main-nonhyperbolic-five-even}
 r_\alpha^{\mathrm{nh}}\ge
 \max\left\{2b_{2-p}(B),\,
 4b_{2-p}(B)-\chi(B),\,b_2(B)-\chi(B)\right\}.
\end{equation}
If $p$ is odd, then
\begin{equation}\label{eq:main-nonhyperbolic-five-odd}
 r_\alpha^{\mathrm{nh}}\ge2b_1(B)+2.
\end{equation}
\par

\item
Suppose $n\ge4$ is even. If $p$ is even, then
\begin{equation}\label{eq:main-nonhyperbolic-even}
 r_\alpha^{\mathrm{nh}}\ge
 \max\left\{
 \begin{aligned}
  &2b_{n-p}(B)+|A_p|,\\
  &4b_{n-p}(B)+|A_p|-A_p-\chi(B),\\
  &1-\chi(B)
 \end{aligned}
 \right\}.
\end{equation}
If $p$ is odd, then
\begin{equation}\label{eq:main-nonhyperbolic-same}
 r_\alpha^{\mathrm{nh}}\ge
 \max\left\{
 \begin{aligned}
  &2b_{n-p}(B)+|A_p|,\\
  &A_p-\chi(B)+1,\\
  &2A_p-\chi(B),\\
  &b_{n-p}(B)+1-\left\lceil\frac{\chi(B)}2\right\rceil
 \end{aligned}
 \right\}.
\end{equation}
\par

\item
Suppose $n\ge3$ is odd. If $p$ is even, then
\begin{equation}\label{eq:main-nonhyperbolic-odd-even}
 r_\alpha^{\mathrm{nh}}\ge
 \max\left\{2b_{n-p}(B)+|A_p|,\,
 4b_{n-p}(B)+|A_p|-A_p+\chi(B)\right\}.
\end{equation}
If $p<n$ is odd, then
\begin{equation}\label{eq:main-nonhyperbolic-odd-low}
 r_\alpha^{\mathrm{nh}}\ge
 \max\left\{
 \begin{aligned}
  &2b_{n-p}(B)+|A_p|,\\
  &A_p+\chi(B)+1,\\
  &2A_p+\chi(B),\\
  &b_{n-p}(B)+1+\left\lfloor\frac{\chi(B)}2\right\rfloor
 \end{aligned}
 \right\}.
\end{equation}
If $p\ge n$ is odd, then
\begin{equation}\label{eq:main-nonhyperbolic-odd-high}
 r_\alpha^{\mathrm{nh}}\ge
 \max\left\{2b_{n-p}(B),\,\chi(B),\,
 b_{n-p}(B)+1+\left\lfloor\frac{\chi(B)}2\right\rfloor\right\}.
\end{equation}
\par
\end{enumerate}
\end{theorem}

\begin{remark}\label{rem:nonhyperbolic-even-homology}
Under the hypotheses of Theorem~\ref{thm:main-nonhyperbolic},
suppose additionally that $H_{\mathrm{odd}}(B;\mathbb Q)=0$.
Then $\chi(B)=b(B)$.

If $n$ is odd and either $p$ is even, $p=1$, or $p\ge n$, then
\[
 r_\alpha^{\mathrm{nh}}\ge b(B).
\]
For odd $n$ and odd $3\le p<n$,
\eqref{eq:main-nonhyperbolic-odd-low} reduces to
\[
 r_\alpha^{\mathrm{nh}}\ge
 \max\left\{
  2\sum_{j=0}^{n-p}b_j(B),\,
  b(B)+1-2\sum_{j=0}^{n-p-1}b_j(B)
 \right\}.
\]

If $n\ge4$ is even and $p>0$, each maximum in
\eqref{eq:main-nonhyperbolic-even} and
\eqref{eq:main-nonhyperbolic-same} equals its first entry, giving
\[
 r_\alpha^{\mathrm{nh}}\ge2\sum_{j=0}^{n-p}b_j(B),
\]
where empty sums are zero.
If $n$ is even and $p=0$, \eqref{eq:main-nonhyperbolic}
reduces to $r_\alpha^{\mathrm{nh}}\ge0$.
For $n=2$, we obtain
\[
 r_\alpha^{\mathrm{nh}}\ge
 \begin{cases}
  2,&p\text{ odd or }p=2,\\
  0,&p=0\text{ or even }p\ge4.
 \end{cases}
\]
\end{remark}

\begin{remark}\label{rem:nonhyperbolic-comparison}
Related lower bounds for non-hyperbolic closed Reeb orbits were
obtained by Abreu--Macarini \cite{AbM17},
Ginzburg--G\"urel--Macarini \cite{GGM18},
Duan--Liu--Long--Wang \cite{DLLW18,DLLW24}, and
Abreu--Macarini \cite{AbM26}, under different
geometric and index hypotheses.
Using contact homology, Ginzburg--G\"urel--Macarini
\cite[Remark~6.3]{GGM18} also treated positive and negative
monotone bases without assuming a filling.
{Theorem~\ref{thm:main-nonhyperbolic} applies for
either sign of $\sigma$ and every integer $p\ge0$, without
requiring a filling or index-admissibility of $\alpha$.}

For comparison, suppose that \cite[Theorem~7]{AbM26} applies
to the closed orbits contractible in $M${, as in
\cite[Remark~2.11]{AbM26}}, with homology ranks
given by the right-hand side of \eqref{eq:BWpattern} for $\sigma=1$.
{For positive even $p$, their estimate becomes}
\[
 r_\alpha^{\mathrm{nh}}\ge2b_{n-p}(B)-(-1)^n\chi(B),
\]
whereas Theorem~\ref{thm:main-nonhyperbolic} gives
\[
 r_\alpha^{\mathrm{nh}}\ge4b_{n-p}(B)-(-1)^n\chi(B).
\]
For $p=0$ and $p=1$, the estimates of \cite[Theorem~7]{AbM26}
also follow from Theorem~\ref{thm:main-nonhyperbolic}.
Theorem~\ref{thm:main-nonhyperbolic} also applies to odd $p\ge3$.
\par
\end{remark}

We now turn to the existence of infinitely many closed
Reeb orbits when $c_B=0$. In this case,
$c_1(TB)|_{\pi_2(B)}=0$,
so $\xi\simeq\pi^*TB$ implies $c_1(\xi)|_{\pi_2(M)}=0$.
When $\omega|_{\pi_2(B)}=0$, we compute linearized
contact homology in the noncontractible fiber classes.
If $\omega|_{\pi_2(B)}\ne0$, the condition on the Betti numbers
in Theorem~\ref{thm:homology} is equivalent to
$H_{\mathrm{odd}}(B;\mathbb Q)=0$;
see Lemma~\ref{lem:betti-support}\textup{(i)}.

\begin{theorem}\label{thm:zero-chern}
Let $\pi:(M^{2n+1},\xi)\to(B^{2n},\omega)$, $n\ge1$, be a
prequantization circle bundle over a closed connected integral
symplectic manifold, with Euler class $-[\omega]$.
Assume that $c_B=0$, and let $\alpha$ be a nondegenerate
contact form supporting $\xi$.
\begin{enumerate}
\renewcommand{\labelenumi}{\textup{(\theenumi)}}
\item If $\omega|_{\pi_2(B)}=0$, then $\alpha$ has
infinitely many geometrically distinct simple closed Reeb orbits.
\item If $\omega|_{\pi_2(B)}\ne0$ and
$H_{\mathrm{odd}}(B;\mathbb Q)=0$, then $\alpha$ has infinitely
many geometrically distinct contractible closed Reeb orbits.
\end{enumerate}
\end{theorem}

\begin{remark}\label{rem:zero-chern-comparison}
Ginzburg, G\"urel, and Macarini
\cite[Theorem~6.4]{GGM18} proved the existence of infinitely many
simple closed Reeb orbits with contractible projections to $B$
when $B$ is aspherical and $c_1(\xi)$ is atoroidal.
For nondegenerate contact forms, Ginzburg and Shon
\cite[Theorem~6.5]{GS18} obtained the same conclusion when
$B$ is symplectically aspherical and $\pi_1(B)$ is torsion free.
Theorem~\ref{thm:zero-chern}\textup{(i)} does not require
$\pi_1(B)$ to be torsion free, but does not assert that the
simple orbits have contractible projections to $B$.
\end{remark}

{We also obtain multiplicity and non-hyperbolicity bounds
for prequantization $S^1$-orbibundles}
$\pi:(M^{2n+1},\xi)\to(B^{2n},\omega)$ over closed connected
symplectic orbifolds, {with smooth total space $M$};
see \cite{KeL21}.
A circle action on $(B,\omega)$ is \emph{Hamiltonian} if its
generating vector field $X$ satisfies $\iota_X\omega=-dH$
for some smooth function $H:B\to\mathbb R$; see \cite{LeT97}.

{Suppose that $B$ admits a Hamiltonian circle action
with isolated fixed points, and write $B^{S^1}$ for its
fixed-point set.}
By Lemma~\ref{lem:orbifold-topology},
$\pi_1(M)$ is finite cyclic and $\#B^{S^1}=b(B)$.
Every simple closed Reeb orbit therefore has a contractible iterate.
{Consequently, $r_\alpha$ is the total number of
geometrically distinct simple closed Reeb orbits of $\alpha$,
and $r_\alpha^{\mathrm{nh}}$ is the number of those that are
non-hyperbolic.}

\Needspace{6\baselineskip}
\begin{theorem}\label{thm:main-orbifold}
Let $\pi:(M^{2n+1},\xi)\to(B^{2n},\omega)$, $n\ge1$, be a
prequantization $S^1$-orbibundle over a closed connected
symplectic orbifold, with $M$ smooth.
Fix a connection form $\beta$ with $\ker\beta=\xi$
and $d\beta=\pi^*\omega$.
Assume that $c_1(\xi)|_{\pi_2(M)}=0$ and that $B$ admits a
Hamiltonian circle action with isolated fixed points.
For each $x\in B^{S^1}$, let $\gamma_{x,\beta}$ denote the least
contractible closed Reeb orbit of $\beta$ along the fiber over $x$.
Assume that
\[
 \widehat\mu(\gamma_{x,\beta})\ne0
 \qquad\text{for every }x\in B^{S^1}.
\]

Let $\alpha$ be a nondegenerate contact form supporting $\xi$,
and fix an integer $p\ge0$ with $p\not\equiv n\pmod2$.
Assume that every good contractible closed Reeb orbit $\gamma$
of $\alpha$ satisfies
\[
 \widehat\mu(\gamma)\ne0,
 \qquad \mu(\gamma)\notin\{-p,p\}.
\]
Then $r_\alpha\ge b(B)$.
If $n$ is odd and $r_\alpha<\infty$, then also
$r_\alpha^{\mathrm{nh}}\ge b(B)$.
\end{theorem}

\begin{remark}\label{rem:orbifold-comparison}
Abreu--Macarini
{\cite[Corollary~7]{AbM26}}
obtain the lower bound $b(B)$ under their filling or non-filling
hypotheses {and index assumptions on $\alpha$}.
{They assume that the connection form $\beta$ is
index-positive, which implies that its contractible closed
Reeb orbits have positive mean index.}
Theorem~\ref{thm:main-orbifold} requires neither a filling nor
index-admissibility {of $\alpha$}.
{We assume that the mean index
$\widehat\mu(\gamma_{x,\beta})$ is nonzero for every
$x\in B^{S^1}$, but we do not require it to be positive.}
\end{remark}

\subsection{Outline of the proofs}
\label{subsec:intro-proof-outline}

The proof of Theorem~\ref{thm:homology} begins with a Morse--Bott
calculation below a fixed action cutoff.
In Section~\ref{sec:local-bw-proof}, we use Bourgeois's local
comparison \cite{Bou} to identify the cylinder differential
for each covering multiplicity with the Morse differential
on $B$, up to a grading shift.
In Section~\ref{sec:reference-tower}, we then construct
augmentations that are compatible with the cobordism maps
as the action cutoff increases.
The hypothesis on the Betti numbers allows us to apply
Lemma~\ref{lem:layer-augmentation-extension} to extend each
augmentation as generators of higher covering multiplicity
are adjoined.
By Lemma~\ref{lem:betti-support}, this hypothesis also implies
that the spectral sequence associated with the filtration
by covering multiplicity collapses at $E^1$.

Proposition~\ref{prop:reference-homology} identifies the homology
of the direct limit of these linearized complexes with the
direct sum in \eqref{eq:BWpattern}.
For $c>0$ sufficiently large, we use the cobordism maps to
construct an augmentation of $\mathcal B(c\alpha)$.
The induced map from the corresponding linearized complex
to this direct limit is a quasi-isomorphism.
The rescaling isomorphism \eqref{eq:contact-rescaling} then
gives the augmentation and homology formula for $\alpha$,
completing the proof of Theorem~\ref{thm:homology}.

In Section~\ref{sec:complete-mixed-sign-count}, we prove
the multiplicity and non-hyperbolicity bounds for $c_B>0$
by adapting the arguments in \cite[Section~6]{AbM26}.
We combine the homology formula \eqref{eq:BWpattern} with
index recurrence, the resonance relations, and Morse
inequalities on finite intervals.
In particular, we count good contractible orbits with
Conley--Zehnder indices $\pm(d-p)$ and $\pm(d+p)$, where $d$
is the even integer provided by Lemma~\ref{lem:index-recurrence}.
The homology formula gives lower bounds for these counts
in terms of $b_{n-p}(B)$.
For $c_B=0$,
$|\mu(\gamma^k)-k\widehat\mu(\gamma)|<n$
shows that only finitely many iterates of a fixed orbit
can have index $n$ or $-n$.
The homology computations, however, give nonzero homology
in degree $n$ in the free homotopy class of the $k$th iterate of a
simple fiber orbit, for every $k\ge1$, or infinite rank
in degree $-n$ in the contractible class.
These facts imply that there are infinitely many
geometrically distinct closed orbits.
In the orbifold case, we perturb the connection form to a
nondegenerate form with only elliptic closed orbits and
compare its linearized contact homology with that of $\alpha$.

\section{Preliminaries}
\label{sec:nd-foundations}

\subsection{Index conventions and iteration}
\label{subsec:index-iteration}

We use the conventions of \cite[Section~3.1]{GGM18} for
the mean and Conley--Zehnder indices.
Let $\operatorname{Sp}(2n)$ denote the symplectic group
of $(\mathbb R^{2n},\omega_0)$, where $\omega_0$ is the
standard symplectic form.
Its universal cover $\widetilde{\operatorname{Sp}}(2n)$
consists of homotopy classes of continuous paths
$\Phi:[0,1]\to\operatorname{Sp}(2n)$ with $\Phi(0)=I$,
where homotopies fix both endpoints.
The powers $\Phi^k$, $k\in\mathbb Z$,
are represented by $t\mapsto\Phi(t)^k$.

For such a path $\Phi$, let $\widehat\mu(\Phi)\in\mathbb R$
denote its mean index, which is invariant under homotopies
with fixed endpoints.
The resulting function
$\widehat\mu:\widetilde{\operatorname{Sp}}(2n)\to\mathbb R$
is the unique continuous homogeneous quasimorphism
satisfying the normalization
\[
 \widehat\mu(\Phi_0)=2,
 \qquad
 \Phi_0(t)=e^{2\pi i t}\oplus I_{2n-2},
 \qquad t\in[0,1];
\]
see \cite{BG92}.

A path $\Phi$ is \emph{nondegenerate} if
$\det(\Phi(1)-I)\ne0$.
For such a path, let $\mu(\Phi)\in\mathbb Z$
denote its Conley--Zehnder index.
We normalize the index so that the path
$t\mapsto e^{i\theta t}$ in $\operatorname{Sp}(2)$
has index $1$ for $0<\theta<2\pi$.
The Conley--Zehnder index is invariant under homotopies
through nondegenerate paths starting at $I$.
We also have
\begin{gather}
 \widehat\mu(\Phi^k)=k\widehat\mu(\Phi)
 \qquad(k\ge1),
 \label{eq:index-homogeneity}\\
 \widehat\mu(\Phi^{-1})=-\widehat\mu(\Phi),\qquad
 \mu(\Phi^{-1})=-\mu(\Phi).
 \label{eq:index-homogeneity-inversion}
\end{gather}

For every nondegenerate path $\Phi$, we have
\begin{equation}\label{eq:index-mean-bound}
 \bigl|\mu(\Phi)-\widehat\mu(\Phi)\bigr|<n.
\end{equation}
A path $\Phi$ is \emph{strongly nondegenerate} if
$\Phi^k$ is nondegenerate for every $k\ge1$.
For such a path, \eqref{eq:index-mean-bound} and
\eqref{eq:index-homogeneity} give
\begin{equation}\label{eq:count-mean-estimate}
 \bigl|\mu(\Phi^k)-k\widehat\mu(\Phi)\bigr|<n
 \qquad(k\ge1),
\end{equation}
and hence
\[
 \widehat\mu(\Phi)
 =\lim_{k\to\infty}\frac{\mu(\Phi^k)}{k}.
\]
If $\Phi(1)$ is hyperbolic, then
\begin{equation}\label{eq:index-hyperbolic}
 \mu(\Phi^k)=k\widehat\mu(\Phi)
 \qquad(k\ge1).
\end{equation}

We also use the Robbin--Salamon index
$\mu_{\mathrm{RS}}(\Phi)\in\tfrac12\mathbb Z$, which
extends the Conley--Zehnder index to paths with degenerate
endpoints.
It is defined as the Maslov index of the path
$t\mapsto\operatorname{graph}\Phi(t)$ relative to the
diagonal in
$(\mathbb R^{2n}\oplus\mathbb R^{2n},-\omega_0\oplus\omega_0)$;
see \cite{RS93}.
It is invariant under homotopies with fixed endpoints
and vanishes on the constant path $\Phi(t)\equiv I$.

Let $\alpha$ be a contact form supporting $\xi$, and let
$\gamma$ be a closed Reeb orbit of period $T$.
Fix a symplectic trivialization
$\tau_t:\xi_{\gamma(Tt)}\to\mathbb R^{2n}$,
$t\in\mathbb R/\mathbb Z$.
In this trivialization, the linearized Reeb flow is
represented by the symplectic path
\[
 \Phi_\gamma^\tau(t)
 =\tau_t\circ
 d\varphi_{R_\alpha}^{Tt}\big|_{\xi_{\gamma(0)}}
 \circ\tau_0^{-1},
 \qquad 0\le t\le1.
\]
We define the mean and Robbin--Salamon indices of $\gamma$
with respect to $\tau$ by
\[
 \widehat\mu(\gamma;\tau)
 =\widehat\mu(\Phi_\gamma^\tau),
 \qquad
 \mu_{\mathrm{RS}}(\gamma;\tau)
 =\mu_{\mathrm{RS}}(\Phi_\gamma^\tau).
\]
When $\gamma$ is nondegenerate, its Conley--Zehnder index
with respect to $\tau$ is
$\mu(\gamma;\tau)=\mu(\Phi_\gamma^\tau)$.
Along each iterate $\gamma^k$, we use the trivialization
$\tau^k$ induced by $\tau$.

If $\tau'$ is another symplectic trivialization along $\gamma$,
then $t\mapsto\tau'_t\tau_t^{-1}$ is a loop in
$\operatorname{Sp}(2n)$.
Denote its winding number by
$\operatorname{wind}(\tau'\tau^{-1})\in\mathbb Z$,
normalized so that a positive full rotation in one
symplectic plane has winding number $1$.
For a nondegenerate orbit $\gamma$,
\cite[Section~3.2]{GGM18} gives
\begin{equation}\label{eq:index-trivialization-change}
 \begin{aligned}
 \mu(\gamma;\tau')-\mu(\gamma;\tau)
 &=2\operatorname{wind}(\tau'\tau^{-1}),\\
 \widehat\mu(\gamma;\tau')-\widehat\mu(\gamma;\tau)
 &=2\operatorname{wind}(\tau'\tau^{-1}).
 \end{aligned}
\end{equation}
The same formula holds for the mean and Robbin--Salamon
indices even when $\gamma$ is degenerate.

For a contractible orbit $\gamma$, choose a capping disk
$u:D^2\to M$ with
\[
 u(e^{2\pi it})=\gamma(Tt),
 \qquad t\in\mathbb R/\mathbb Z.
\]
A symplectic trivialization of $u^*\xi$ induces a
trivialization along $\gamma$ by restriction to $\partial D^2$.
Since $D^2$ is contractible, the homotopy class of this
induced trivialization is independent of the choice of
symplectic trivialization of $u^*\xi$.
If $c_1(\xi)|_{\pi_2(M)}=0$, the resulting indices are
also independent of the capping disk.
Under this assumption, we write $\widehat\mu(\gamma)$,
$\mu_{\mathrm{RS}}(\gamma)$, and, when $\gamma$ is
nondegenerate, $\mu(\gamma)$.

By \eqref{eq:index-trivialization-change}, the parity of
the Conley--Zehnder index is independent of the choice of
trivialization.
A closed Reeb orbit $\gamma$ is \emph{strongly nondegenerate}
if $\gamma^k$ is nondegenerate for every $k\ge1$.
For such an orbit,
\cite[Section~3.2]{GGM18} gives
\[
 \mu(\gamma^{k+2};\tau^{k+2})
 \equiv\mu(\gamma^k;\tau^k)\pmod2
 \qquad(k\ge1).
\]
For a strongly nondegenerate simple orbit, either all iterates
are good or precisely the even iterates are bad.
If $\gamma$ is elliptic, every nondegenerate iterate
$\gamma^k$ satisfies
\begin{equation}\label{eq:index-elliptic-parity}
 \mu(\gamma^k;\tau^k)\equiv n\pmod2,
\end{equation}
where $\tau$ is any symplectic trivialization of $\xi$
along $\gamma$.
In particular, every such iterate is good.

\Needspace{5\baselineskip}
Index recurrence theorems play an important role in proving
multiplicity results for periodic orbits, and several versions
have been developed in the literature.
We use the following version, stated in
\cite[Theorem~4]{AbM26}.
It also follows from Theorem~3.6 of \cite{DLLW24}
and Step~2 in the proof of their Theorem~1.2.
For related versions, see \cite{LZ02,DLW16,GGM18,HW22}.

\begin{lemma}\label{lem:index-recurrence}
Let $\Phi_1,\ldots,\Phi_r$ be a finite collection of strongly
nondegenerate elements of $\widetilde{\operatorname{Sp}}(2n)$ with
$\widehat\mu(\Phi_i)\ne0$ for all $1\le i\le r$.
Denote by $\sigma_i=\widehat\mu(\Phi_i)/|\widehat\mu(\Phi_i)|$
the sign of $\widehat\mu(\Phi_i)$.
Then, for any $\eta>0$ and any $\ell_0\in\mathbb N$, there exist
two integer sequences $d_j^\pm\to\infty$ and two sequences of
integer vectors $\vec k_j^\pm=(k_{1j}^\pm,\ldots,k_{rj}^\pm)$,
with all components going to infinity as $j\to\infty$, such that
for all $i$ and $j$, and all $\ell\in\mathbb Z$ in the range
$1\le|\ell|\le\ell_0$, we have
\begin{enumerate}
\item\label{eq:count-recurrence-mean}
$\bigl|\widehat\mu(\Phi_i^{k_{ij}^\pm})-\sigma_i d_j^\pm\bigr|<\eta$
with the equality
$\widehat\mu(\Phi_i^{k_{ij}^\pm})
 =\mu(\Phi_i^{k_{ij}^\pm})=\sigma_i d_j^\pm$
whenever $\Phi_i(1)$ is hyperbolic,
\item\label{eq:count-recurrence-offset}
$\mu(\Phi_i^{k_{ij}^\pm+\ell})
 =\sigma_i d_j^\pm+\mu(\Phi_i^\ell)$, and
\item\label{eq:count-recurrence-companion}
$\mu(\Phi_i^{k_{ij}^-})-\sigma_i d_j^-
 =-\bigl(\mu(\Phi_i^{k_{ij}^+})-\sigma_i d_j^+\bigr)$.
\end{enumerate}
Furthermore, for any $N\in\mathbb N$, we can make all $d_j^\pm$
and $k_{ij}^\pm$ divisible by $N$.
\end{lemma}

\subsection{Holomorphic curves and Morse trajectories}
\label{subsec:holomorphic-morse}

Let $(M^{2n+1},\xi)$ be a closed cooriented contact manifold with
$c_1(\xi)|_{\pi_2(M)}=0$, and let $\alpha$ be a contact form
supporting $\xi$.
The action of a closed Reeb orbit $\gamma$ of period $T$ is
\[
 \mathcal A_\alpha(\gamma)=\int_\gamma\alpha=T.
\]
The \emph{action spectrum} of $\alpha$ is
\[
 \operatorname{Spec}(\alpha)
 =\bigl\{\mathcal A_\alpha(\gamma)
   \mid \gamma \text{ is a closed Reeb orbit of }\alpha\bigr\}.
\]

The symplectization of $\alpha$ is
$(\mathbb R\times M,d(e^a\alpha))$, where $a$ denotes the
coordinate on the $\mathbb R$-factor.
An almost complex structure $J$ on this symplectization is
\emph{cylindrical} if it is $\mathbb R$-invariant, satisfies
$J\partial_a=R_\alpha$ and $J\xi=\xi$, and its restriction
to $\xi$ is $d\alpha$-compatible; see
\cite[Section~2.1]{Bou}.

For such pairs $(\alpha,J)$, we use the topology defined by
\begin{equation}\label{eq:contact-pair-topology}
 (\alpha_\nu,J_\nu)\longrightarrow(\alpha,J)
 \quad\Longleftrightarrow\quad
 \alpha_\nu\longrightarrow\alpha
 \ \text{and}\ J_\nu|_\xi\longrightarrow J|_\xi,
\end{equation}
where both convergences on the right are in the $C^\infty$
topology on $M$.

A smooth map $u:(\Sigma,j)\to(\mathbb R\times M,J)$ from a
Riemann surface is $J$-holomorphic if $du\circ j=J\circ du$.
On the cylinder $\mathbb R\times(\mathbb R/\mathbb Z)$,
with coordinates $(s,t)$ and complex structure
$j\partial_s=\partial_t$, the Cauchy--Riemann equation
takes the form
\begin{equation}\label{eq:holomorphic-cylinder-equation}
 u_s+J(u)u_t=0.
\end{equation}

For a punctured $J$-holomorphic curve $u=(a,v)$, where
$a:\Sigma\to\mathbb R$ and $v:\Sigma\to M$, we define
its Hofer energy following \cite[(2.28)--(2.30)]{Pardon}:
\[
 E(u)=\int_\Sigma v^*d\alpha+
 \sup_{\varphi}\int_\Sigma
 \varphi'(a)\,da\wedge v^*\alpha.
\]
The supremum is taken over smooth nondecreasing functions
$\varphi:\mathbb R\to[0,1]$ equal to $0$ near $-\infty$
and to $1$ near $+\infty$.
The curves considered below have finite energy and are
asymptotic to closed Reeb orbits at their punctures.
A puncture is positive or negative according as
$a\to+\infty$ or $a\to-\infty$.
If $u$ has one positive end $\gamma^+$ and negative ends
$\gamma_1^-,\ldots,\gamma_r^-$, Stokes' theorem gives
\begin{equation}\label{eq:symplectization-action-energy}
 \begin{gathered}
 0\le\int_\Sigma v^*d\alpha
 =\mathcal A_\alpha(\gamma^+)
   -\sum_{i=1}^r\mathcal A_\alpha(\gamma_i^-),\\
 E(u)\le 2\mathcal A_\alpha(\gamma^+).
 \end{gathered}
\end{equation}
The integral $\int_\Sigma v^*d\alpha$ is the $d\alpha$-area
of $u$.
If this area vanishes, the $\xi$-component of $dv$ vanishes.
For a connected curve with positive end $\gamma^+$, the image
of $u$ then lies in $\mathbb R\times\gamma^+$.
A trivial cylinder over $\gamma$ has image
$\mathbb R\times\gamma$ and has zero $d\alpha$-area
and positive Hofer energy, even when multiply covered.

For curves with nondegenerate ends, we follow the marking
conventions of \cite[Section~1.2]{Pardon}.
An \emph{asymptotic marker} is a ray in the complex tangent
line at a puncture of the compactified domain.
We fix a basepoint on each underlying simple Reeb orbit
and require the asymptotic parametrization at each puncture
to send its marker to the corresponding basepoint.
Moduli spaces are taken modulo biholomorphic
reparametrizations preserving the labeled punctures and
their asymptotic markers and, in the symplectization,
modulo translations in the $\mathbb R$-coordinate.

Let $D_u$ denote the linearized Cauchy--Riemann operator
at $u$, as in \cite[(2.42)]{Pardon}.
We write $\operatorname{ind}(u)=\operatorname{ind}D_u$
for its Fredholm index and call $u$ \emph{regular} if
$D_u$ is surjective.

Let $u$ be a genus-zero $J$-holomorphic curve with one
positive end $\gamma^+$ and negative ends
$\gamma_1^-,\ldots,\gamma_r^-$.
Assume that all asymptotic orbits are nondegenerate and
contractible.
For each orbit, choose a capping disk and compute its
Conley--Zehnder index using the induced trivialization.

Capping the projection of $u$ to $M$ with these disks
gives a sphere.
Since $c_1(\xi)|_{\pi_2(M)}=0$, the relative Chern term
in the index formula vanishes.
Thus \cite[(2.50)]{Pardon} gives
\begin{equation}\label{eq:capped-punctured-index}
 \operatorname{ind}(u)
 =(n-2)(1-r)+\mu(\gamma^+)-\sum_{i=1}^r\mu(\gamma_i^-).
\end{equation}

Let $u$ be a nontrivial regular cylinder with nondegenerate
contractible ends $\gamma^+$ and $\gamma^-$.
With the marking and quotient conventions specified above,
the moduli space is smooth near the class of $u$ and has
dimension
\begin{equation}\label{eq:cylinder-moduli-dimension}
 \operatorname{ind}(u)-1
 =\mu(\gamma^+)-\mu(\gamma^-)-1.
\end{equation}

We use coherent orientations compatible with gluing;
see \cite[Section~2.13]{Pardon}.
For a compact oriented zero-dimensional marked moduli space
$\mathcal M$, we write $\#\mathcal M$ for its signed count.

In the local Morse--Bott calculation of
Section~\ref{sec:local-bw-proof}, we use Morse trajectories
on the base to compute signed counts of rigid holomorphic
cylinders.
For background on Morse theory, see \cite{Sch93}.

Let $f:B\to\mathbb R$ be a Morse function on a closed
Riemannian manifold $(B,g)$.
We denote its critical set by $\operatorname{Crit}(f)$
and its gradient by $\nabla_g f$.
For $x\in\operatorname{Crit}(f)$, write
$\operatorname{ind}_f(x)$ for its Morse index.

For $x,y\in\operatorname{Crit}(f)$, a parametrized
\emph{Morse trajectory} from $x$ to $y$ is a smooth map
$z:\mathbb R\to B$ satisfying
\begin{equation}\label{eq:morse-trajectory}
 \frac{dz}{dr}=-\nabla_g f(z(r)),
 \qquad
 \lim_{r\to-\infty}z(r)=x,
 \qquad
 \lim_{r\to+\infty}z(r)=y.
\end{equation}
For $x\ne y$, let $\mathcal M_f(x,y)$ denote the space of
\emph{unparametrized Morse trajectories} from $x$ to $y$, that is,
solutions of \eqref{eq:morse-trajectory} modulo time translation.

The pair $(f,g)$ is \emph{Morse--Smale} if the unstable
and stable manifolds $W^u(x)$ and $W^s(y)$ of the negative
gradient flow intersect transversely for all critical
points $x,y$.
In this case, for distinct critical points $x,y$,
$\mathcal M_f(x,y)$ is a smooth manifold of dimension
$\operatorname{ind}_f(x)-\operatorname{ind}_f(y)-1$
whenever it is nonempty.
We choose orientations of the unstable manifolds and
equip the spaces $\mathcal M_f(x,y)$ with the induced
orientations.

For a Morse--Smale pair $(f,g)$ and critical points $x,y$ with
$\operatorname{ind}_f(x)-\operatorname{ind}_f(y)=1$,
$\mathcal M_f(x,y)$ is compact and zero-dimensional.
Its elements are called \emph{rigid Morse trajectories}.
We denote their signed count by
\[
 n_f(x,y)=\#\mathcal M_f(x,y).
\]
The Morse differential defined by these counts lowers
the Morse index by one.

A \emph{broken Morse trajectory} is a finite sequence of
nonconstant Morse trajectories in which each trajectory
ends at the critical point where the next one begins.

\subsection{Contact algebras and linearized contact homology}
\label{subsec:contact-algebras}
\label{subsec:free-homotopy}

We work over $\mathbb Q$.
A differential graded algebra (DGA) is a graded algebra
$A$ equipped with a differential $d$ of degree $-1$
such that $d^2=0$ and
\[
 d(ab)=d(a)b+(-1)^{|a|}a\,d(b)
\]
for homogeneous elements $a,b\in A$, where $|a|$
denotes the degree of $a$.
All DGAs considered below have a unit and are
graded-commutative, meaning that
$ab=(-1)^{|a||b|}ba$.

A Sullivan algebra is a free graded-commutative DGA
admitting a well-ordered set of homogeneous generators
such that the differential of each generator belongs
to the subalgebra generated by the preceding generators;
see \cite[Definition~2.4]{Chaidez}.
For integer gradings, we impose no upper or lower bound
on the generator degrees.

Let $\alpha$ be a nondegenerate supporting contact form.
Choose a cylindrical almost complex structure $J$ and
the auxiliary data required in Pardon's construction,
and denote the resulting contact DGA by $\mathcal C(\alpha)$.
After choosing generators of the orbit orientation lines,
its underlying algebra is the free unital graded-commutative
algebra over $\mathbb Q$ generated by the good Reeb orbits.
We denote each generator by the same symbol as its orbit.
Elements of $\mathcal C(\alpha)$ are finite rational
linear combinations of monomials in these generators,
including the empty monomial $1$.

The differential is determined on generators by
\[
 d\gamma^+
 =\sum_{\Gamma}c_\alpha(\gamma^+;\Gamma)\,\Gamma,
\]
where we fix an ordering of the good orbits and sum over
the corresponding monomial basis.
For $\Gamma=\gamma_1^-\cdots\gamma_r^-$, the coefficient
$c_\alpha(\gamma^+;\Gamma)\in\mathbb Q$ is defined by
virtual counts of finite-energy genus-zero
$J$-holomorphic curves in the symplectization with
positive end $\gamma^+$ and negative ends
$\gamma_1^-,\ldots,\gamma_r^-$.
We use the marking and quotient conventions of
Section~\ref{subsec:holomorphic-morse}.
Each coefficient includes the automorphism factors
specified in \cite[(1.5)]{Pardon}.
Only moduli spaces of expected dimension zero contribute.
The case $r=0$ gives the constant term.
The differential extends by the graded Leibniz rule
and satisfies $d^2=0$; see \cite[Section~1.2]{Pardon}.
When the compactified moduli spaces consist of regular
isolated curves, the virtual counts agree with the ordinary
signed counts \cite[Proposition~4.33]{Pardon}.

The full contact DGA $\mathcal C(\alpha)$ is naturally
$\mathbb Z/2$-graded.
If $c_1(\xi)|_{\pi_2(M)}=0$, each contractible generator
has the integer degree
\begin{equation}\label{eq:contact-grading}
 |\gamma|=\mu(\gamma)+n-2.
\end{equation}

If $c_1(\xi)=0$ in $H^2(M;\mathbb Z)$, fix a
trivialization $\tau$ of $\det_{\mathbb C}\xi$.
For every closed Reeb orbit $\gamma$, we compute $\mu(\gamma)$
using a symplectic trivialization of $\xi$ along $\gamma$
whose determinant agrees with $\tau$. On each iterate
$\gamma^k$, we use the induced trivialization.
The full contact algebra then has the integer grading
\eqref{eq:contact-grading}, and its differential has degree $-1$.
For contractible orbits, these indices agree with those
computed using capping disks.

For contact forms $\alpha^\pm$ defining the same cooriented
contact structure $\xi$, write $\alpha^+=g\alpha^-$ with $g>0$.
We define
\begin{equation}\label{eq:contact-form-order}
 \alpha^+>\alpha^-
 \quad\Longleftrightarrow\quad
 g(x)>1\quad\text{for every }x\in M.
\end{equation}
If $\alpha^\pm$ are nondegenerate and $\alpha^+>\alpha^-$,
the region
\[
 W=\{(s,x)\in\mathbb R\times M:0\le s\le\log g(x)\},
 \qquad \lambda=e^s\alpha^-,
\]
with symplectic form $d\lambda$, is an exact symplectic
cobordism from $(M,\alpha^+)$ to $(M,\alpha^-)$.
Under the natural identifications of its boundary
components with $M$, the restrictions of $\lambda$
are $\alpha^+$ and $\alpha^-$, respectively.
Its completion is diffeomorphic to $\mathbb R\times M$.
Fix the auxiliary data defining
$\mathcal C(\alpha^+)$ and $\mathcal C(\alpha^-)$.
Choose a tame almost complex structure on the completed
cobordism that agrees with the cylindrical structures at both ends.
The prescribed perturbation data at the ends extend over the
cobordism by the surjectivity of Pardon's forgetful map
\cite[Section~1.3]{Pardon}. Each such extension defines a DGA map
\[
 F:\mathcal C(\alpha^+)\longrightarrow\mathcal C(\alpha^-).
\]
For a generator $\gamma^+$, the coefficients of $F(\gamma^+)$
are defined by virtual counts of genus-zero punctured
holomorphic curves in the completed cobordism.
These curves have one positive end at $\gamma^+$ and
negative ends at the orbits appearing in the output monomial.
The map extends multiplicatively to the contact algebra.
Throughout, we use the cobordisms and cylindrical end
identifications chosen in Lemma~\ref{lem:found-action}.

Contractibility of the positive end does not force each
negative end to be contractible.
Let $\mathcal I_{\mathrm{nc}}(\alpha)$ be the ideal in
$\mathcal C(\alpha)$ generated by the good noncontractible
Reeb orbits, and set
\[
 \mathcal B(\alpha)
 =\mathcal C(\alpha)/\mathcal I_{\mathrm{nc}}(\alpha).
\]
Pardon describes the contact algebra generated by
contractible orbits in \cite[Section~1.8]{Pardon}.

\begin{lemma}\label{lem:found-contractible-quotient}
Assume that $c_1(\xi)|_{\pi_2(M)}=0$.
The ideal $\mathcal I_{\mathrm{nc}}(\alpha)$ is preserved
by the differential.
The quotient $\mathcal B(\alpha)$ is Sullivan, with the
integer grading \eqref{eq:contact-grading} on its
contractible generators and differential of degree $-1$.
Each cobordism map described above descends to a map
between the quotients that preserves the integer grading
and induces an isomorphism on unfiltered homology.
\end{lemma}

\begin{proof}
Consider a connected genus-zero curve with one positive
puncture.
If all its negative asymptotic orbits are contractible,
capping the negative ends gives a capping disk for the
positive orbit.
Consequently, every monomial in the differential of a
noncontractible generator contains at least one generator
corresponding to a noncontractible orbit.
The derivation rule then gives
$d\mathcal I_{\mathrm{nc}}\subset\mathcal I_{\mathrm{nc}}$.
The same argument applies to the cobordism maps on
generators, since the completed cobordisms are
diffeomorphic to $\mathbb R\times M$.

Capping all ends of a genus-zero curve gives a sphere
in the symplectization or the completed cobordism.
The inclusion of a slice $\{s\}\times M$ in a cylindrical end
into the completed cobordism is a homotopy equivalence.
On this slice the complex
tangent bundle of the cobordism restricts to $\xi\oplus\mathbb C$,
so its first Chern class also vanishes on spheres.
The relative Chern term therefore vanishes in the
Fredholm index formula
\eqref{eq:capped-punctured-index}.
The differential counts curves of Fredholm index one,
whereas the cobordism map counts curves of index zero.
Thus the differential lowers the degree by one and the
cobordism map preserves it.

Order the generators by action, breaking ties arbitrarily.
There are finitely many generators below any fixed
action bound.
Every generator appearing in the differential of a
generator has strictly smaller action.
By \eqref{eq:symplectization-action-energy}, the sum of the
negative actions is at most the positive action.
If there are several negative ends, each individual
negative action is smaller than the positive action.
With one negative end, equal actions give zero
$d\alpha$-area. The curve is then a trivial cylinder, possibly
multiply covered. Its Fredholm index is zero, so it does not
contribute to the differential.
This ordering makes the full contact algebra Sullivan
\cite[Lemma~3.1]{Chaidez}. Its restriction to the contractible
generators shows that the quotient is also Sullivan.

To prove the quasi-isomorphism assertion, we use Pardon's
deformation and composition homotopies
\cite[Sections~1.4--1.5]{Pardon}. In Setup III, the ordinary
chain homotopy is defined on each monomial using disconnected
configurations, with one connected genus-zero curve for each
orbit in the monomial.
If one of these positive orbits is noncontractible, its component
must have a noncontractible negative end.
Thus the homotopy preserves $\mathcal I_{\mathrm{nc}}$.
The same argument applies to the composition homotopy
of Setup IV, so the deformation and composition
identities descend to the quotient.
The completion of $W$ is the symplectization of $\xi$
with end identifications determined by $\alpha^+$ and $\alpha^-$.
It is therefore the trivial exact cobordism of \cite[Lemma~1.2]{Pardon}.
By the composition identity on the quotient, it suffices, as in
that lemma, to consider equal contact forms and almost complex
structures at the two ends. For the translation-invariant almost
complex structure, only trivial cylinders contribute to the associated
graded map for the filtration used there. The trivial cylinder over
each contractible orbit identifies its orientation lines at the two
ends, so the associated graded map is an isomorphism on the quotient
as well. The induced map on unfiltered homology is therefore an
isomorphism.
\end{proof}

For $L>0$, let $\mathcal B(\alpha)^{<L}$ be the unital
subalgebra generated by the good contractible orbits
of action $<L$.
A monomial in $\mathcal B(\alpha)$ belongs to this subalgebra
whenever each of its factors has action $<L$, even if their
total action is at least $L$.
It is a Sullivan sub-DGA because the differential
of each generator involves only generators of
strictly smaller action.

\begin{lemma}\label{lem:found-action}
Let $\alpha^\pm$ be nondegenerate supporting contact
forms with $\alpha^+>\alpha^-$.
Fix cylindrical almost complex structures at both ends.
One can choose an exact symplectic cobordism of the type
described above, with fixed identifications of its
cylindrical ends and an admissible almost complex
structure extending the prescribed end structures,
so that the induced DGA map $F$ has the following property.
For each generator $\gamma^+$ and each monomial
$\gamma_1^-\cdots\gamma_r^-$ with nonzero coefficient in
$F(\gamma^+)$, we have
\begin{equation}\label{eq:found-total-output-action}
 \sum_i\mathcal A_{\alpha^-}(\gamma_i^-)
 \leq\mathcal A_{\alpha^+}(\gamma^+).
\end{equation}
The sum is zero for the constant monomial.
Consequently,
\[
 F\bigl(\mathcal B(\alpha^+)^{<L}\bigr)
 \subset\mathcal B(\alpha^-)^{<L}
 \qquad (L>0).
\]
\end{lemma}

\begin{proof}
Fix a supporting contact form $\beta$ and write
$\alpha^\pm=h_\pm\beta$.
Choose a monotone interpolation $h_s$ from $h_-$ to $h_+$,
constant for sufficiently large $|s|$.
On $X=\mathbb R\times M$, use the primitive
$\Lambda=e^s h_s\beta$ and an almost complex structure
satisfying
$\widehat J\partial_s=R_{h_s\beta}$,
with compatible restriction to $\xi$ and the prescribed
structures at the ends.
This gives an admissible exact cobordism in Pardon's convention,
with $\Lambda=e^s\alpha^\pm$ on its cylindrical ends.
The change of coordinates $(s,x)\mapsto(t,x)$, where
$t=s+\log h_s(x)$, pulls $e^t\beta$ back to $\Lambda$.

The closed auxiliary form $\eta=d(h_s\beta)$ is
nonnegative on complex lines.
Its restriction to $\xi$ is $h_s d\beta$.
For every $Y\in\xi$, both $\eta(\partial_s,Y)$ and
$\eta(R_{h_s\beta},Y)$ vanish. On the plane spanned by $\partial_s$
and $R_{h_s\beta}$, we have
\[
 \eta(\partial_s,R_{h_s\beta})
 =\frac{\partial_s h_s}{h_s}\ge0.
\]
For a curve with positive end $\gamma^+$ and negative
ends $\gamma_1^-,\ldots,\gamma_r^-$, Stokes' theorem
therefore gives \eqref{eq:found-total-output-action}.
The prescribed end identifications ensure that the actions in
\eqref{eq:found-total-output-action} are computed with $\alpha^\pm$.

For a holomorphic building, apply this estimate on the cobordism
level and \eqref{eq:symplectization-action-energy} on the
symplectization levels. Summing over the levels cancels the
actions at matching ends. If the asymptotic orbits violate
\eqref{eq:found-total-output-action}, the compactified moduli
space is therefore empty and its virtual count vanishes by
\cite[Theorem~1.1(iv)]{Pardon}.
Since the orbit at each negative end has positive action, every generator
occurring in the image of a generator of action $<L$
also has action $<L$.
The inclusion follows because $F$ is multiplicative.
\end{proof}

The contact differential may have nonzero constant terms,
so its linear part need not square to zero.
An augmentation allows us to eliminate the constant terms
before taking the linear part.
An augmentation is a grading-preserving unital DGA map
$\varepsilon:A\to\mathbb Q$, where $\mathbb Q$ is
concentrated in degree zero and has zero differential.

For an augmented DGA $(A,\varepsilon)$, set
\begin{equation}\label{eq:indecomposables}
 I_A=\ker\varepsilon,
 \qquad
 Q(A,\varepsilon)=I_A/I_A^2.
\end{equation}
Here $I_A^2$ is the ideal generated by products $xy$
with $x,y\in I_A$.
Since $\varepsilon d=0$ and $d$ is a derivation,
$dI_A\subset I_A$ and $d(I_A^2)\subset I_A^2$.
Thus $d$ induces a differential on the space of
indecomposables $Q(A,\varepsilon)$.
The resulting complex is the \emph{linearized complex}
at $\varepsilon$ \cite[Definition~2.10]{Chaidez}.

If $A$ is freely generated by homogeneous elements $v$,
then the classes $[v-\varepsilon(v)]$ form a basis
of the graded vector space $Q(A,\varepsilon)$.
To compute its differential, rewrite $dv$ in terms
of the generators $w-\varepsilon(w)$.
The constant term is $\varepsilon(dv)=0$, and terms
of word length at least two vanish modulo $I_A^2$.
Thus the linearized differential is obtained by
retaining the terms of word length one.

An augmented DGA is called \emph{pointed}
\cite[Definition~2.2]{Chaidez}.
We write $Q(A)$ for $Q(A,\varepsilon_A)$ when the
augmentation is understood.
A DGA map $g:A\to D$ is pointed if
$\varepsilon_D\circ g=\varepsilon_A$. Such a map sends $I_A$
into $I_D$ and therefore induces a chain map
$Qg:Q(A)\to Q(D)$.

\begin{lemma}
\label{lem:found-forward-sullivan}
A pointed quasi-isomorphism $g:S\to T$ between Sullivan algebras induces a
quasi-isomorphism $Qg:Q(S)\to Q(T)$.
\end{lemma}

\begin{proof}
For a $\mathbb Z/2$-grading, this is the forward implication
of \cite[Lemma~2.15]{Chaidez}. For an integer-graded complex $C$,
reducing the grading modulo two gives homology
$\bigoplus_{k\equiv i\bmod 2}H_k(C)$ in parity $i$.
Thus $g$ remains a quasi-isomorphism
after reduction. The ordering of the generators is unchanged, and reduction
commutes with taking the augmentation ideal and its square.
The $\mathbb Z/2$-graded case therefore shows that $Qg$ induces an
isomorphism on homology in each parity. Since $Qg$ preserves the
integer grading, it induces an isomorphism in each degree.
\end{proof}

For a graded augmentation $\varepsilon$ of
$\mathcal C(\alpha)$, let $\LCC_*(M,\alpha;\varepsilon)$ denote
$Q(\mathcal C(\alpha),\varepsilon)$, with the grading shifted
by subtracting $n-2$.
This gives the Conley--Zehnder grading modulo two.
For linearized contact homology of the full contact algebra
at an augmentation class, see \cite[Definition~3.3]{Chaidez}.

Assume that $c_1(\xi)|_{\pi_2(M)}=0$.
Let $\alpha$ be a nondegenerate supporting contact form,
and let $\bar\varepsilon_\alpha$ be a graded augmentation
of $\mathcal B(\alpha)$.
By \eqref{eq:contact-grading},
$\bar\varepsilon_\alpha(\gamma)$ can be nonzero only if
$\mu(\gamma)=2-n$.
Denote its pullback along the quotient map by
\[
 \varepsilon_\alpha:
 \mathcal C(\alpha)\longrightarrow\mathcal B(\alpha)
 \xrightarrow{\ \bar\varepsilon_\alpha\ }\mathbb Q.
\]
Here the quotient map preserves the grading after we reduce
the integer grading on $\mathcal B(\alpha)$ modulo two.
This augmentation vanishes on noncontractible generators.

Using the integer grading on $\mathcal B(\alpha)$, we
define the contractible linearized contact complex in
Conley--Zehnder grading by
\begin{equation}\label{eq:contractible-linearized-complex}
 \LCC^0_k(M,\alpha;\varepsilon_\alpha)
 :=Q(\mathcal B(\alpha),\bar\varepsilon_\alpha)_{k+n-2},
\end{equation}
with differential $\partial_{\varepsilon_\alpha}$
induced by the contact differential.
Here $0$ denotes the contractible free homotopy class in $M$.
The shift in \eqref{eq:contractible-linearized-complex}
converts the contact-algebra grading
\eqref{eq:contact-grading} to the Conley--Zehnder grading.

Its homology is the contractible linearized contact
homology, denoted
\begin{equation}\label{eq:linearized-contact-homology}
 \LCH^0_k(M,\alpha;\varepsilon_\alpha)
 :=H_k\bigl(
 \LCC^0_*(M,\alpha;\varepsilon_\alpha),
 \partial_{\varepsilon_\alpha}
 \bigr).
\end{equation}
We retain $\alpha$ and $\varepsilon_\alpha$ in the
notation; the homology may depend on the augmentation.

\begin{lemma}\label{lem:free-homotopy-splitting}
Let $\varepsilon$ be a graded augmentation of
$\mathcal C(\alpha)$ that vanishes on every noncontractible
generator. Its linearization decomposes as
\[
 \LCC_*(M,\alpha;\varepsilon)
 =\bigoplus_{a\in[S^1,M]}\LCC_*^a(M,\alpha;\varepsilon),
\]
where the classes $[\gamma-\varepsilon(\gamma)]$
for the good orbits $\gamma$ in the free homotopy class $a$
form a basis of the corresponding summand.
The linearizations of pointed DGA maps between such
augmented contact algebras, induced by the cobordisms of
Lemma~\ref{lem:found-action}, preserve this decomposition.
\end{lemma}

\begin{proof}
Each nonzero contribution to the linearized
differential is obtained by retaining one negative orbit
and evaluating the augmentation on all the others.
Since the augmentation vanishes on noncontractible
generators, these other orbits are contractible.
Capping their ends with disks and projecting to $M$
gives a free homotopy between the positive orbit and
the retained negative orbit.
The same argument applies to the linearization of a
pointed cobordism map.
\end{proof}

We denote the homology of
$\LCC_*^a(M,\alpha;\varepsilon)$ by
$\LCH_*^a(M,\alpha;\varepsilon)$.
When an integer grading compatible with the differential
and the linearized cobordism maps is available, we use it
on the corresponding summands.
For any subset $\Gamma\subset[S^1,M]$, let
$\LCC_*^\Gamma=\bigoplus_{a\in\Gamma}\LCC_*^a$ and denote
its homology by $\LCH_*^\Gamma$.

For the augmentation $\varepsilon_\alpha$ pulled back
from $\mathcal B(\alpha)$, the quotient map identifies the
contractible summand with the integer-graded complex
in \eqref{eq:contractible-linearized-complex}.

\subsection{Local contact homology, resonance relations, and Morse inequalities}
\label{subsec:morse-mean-euler}

Let $\gamma$ be an isolated closed Reeb orbit of $\alpha$, possibly
degenerate or multiply covered. Choose an isolating neighborhood
$U$ of its image and a sufficiently small nondegenerate perturbation
$\widetilde\alpha$ of $\alpha$. The closed Reeb orbits of
$\widetilde\alpha$ that are freely homotopic to $\gamma$ in $U$
form a finite set. The good orbits in this set generate the local
chain complex over $\mathbb Q$. Fix a cylindrical almost complex
structure $J$ adapted to $\alpha$. For regular local data
$(\widetilde\alpha,\widetilde J)$ sufficiently
$C^\infty$-close to $(\alpha,J)$, the differential counts
rigid marked $\widetilde J$-holomorphic cylinders in
$\mathbb R\times U$, with each signed count divided by the covering
multiplicity of the negative orbit as in \cite[(16)]{HM15}.
The homology is independent of the choice of such regular data up to isomorphism
\cite[Corollary~4.11]{HM15}. We call it the \emph{local contact
homology} of $\gamma$ and denote it by $\mathrm{HC}_*(\alpha,\gamma)$.

Fix a symplectic trivialization of $\xi$ along $\gamma$. We use
the induced trivializations along the nearby orbits to grade each
generator by its Conley--Zehnder index. For contractible orbits we use a capping
trivialization. When $c_1(\xi)|_{\pi_2(M)}=0$, the resulting local grading
agrees with that in \eqref{eq:contractible-linearized-complex}.

If $\gamma$ is nondegenerate, the computation in
\cite[Example~3.9]{GG15} gives
\[
 \mathrm{HC}_j(\alpha,\gamma)\cong
 \begin{cases}
  \mathbb Q, & \text{if $\gamma$ is good and $j=\mu(\gamma)$},\\
  0,         & \text{otherwise}.
 \end{cases}
\]

For an isolated orbit $\gamma$, the local chain complex is
finite-dimensional. We define the \emph{local Euler characteristic}
of $\gamma$ by
\[
 \chi(\gamma)
 =\sum_{j\in\mathbb Z}(-1)^j
   \dim_{\mathbb Q}\mathrm{HC}_j(\alpha,\gamma).
\]
If $\gamma$ is nondegenerate, then $\chi(\gamma)=(-1)^{\mu(\gamma)}$
when $\gamma$ is good and $\chi(\gamma)=0$ when it is bad.

If every iterate of $\gamma$ is isolated, we define its
\emph{local mean Euler characteristic} by
\[
 \widehat\chi(\gamma)
 =\lim_{N\to\infty}\frac{1}{N}\sum_{k=1}^N\chi(\gamma^k).
\]
The limit exists and is rational by
\cite[Remarks~3.13 and~3.15]{GG15}.

Suppose that $\gamma$ is strongly nondegenerate. For every $k\ge1$,
the iterates $\gamma^k$ and $\gamma^{k+2}$ have the same index
parity and are either both good or both bad
(Section~\ref{subsec:index-iteration}). Hence
$\chi(\gamma^{k+2})=\chi(\gamma^k)$, and the definition of
$\widehat\chi(\gamma)$ gives
\[
 \widehat\chi(\gamma)
 =\frac{\chi(\gamma)+\chi(\gamma^2)}{2}
 =
 \begin{cases}
  0, & \text{if $\gamma$ is bad},\\
  (-1)^{\mu(\gamma)}, & \text{if $\gamma$ and $\gamma^2$ are good},\\
  \tfrac12(-1)^{\mu(\gamma)},
     & \text{if $\gamma$ is good and $\gamma^2$ is bad}.
 \end{cases}
\]
For every positive even integer $N$, it follows that
\begin{equation}\label{eq:count-period-two}
 \sum_{k=1}^N\chi(\gamma^k)=N\widehat\chi(\gamma),\qquad
 \chi(\gamma^{N\pm\ell})=\chi(\gamma^\ell),
\end{equation}
where $\ell\ge1$ is an integer and the minus sign requires $\ell<N$.

Let $\alpha$ be a nondegenerate contact form and fix an augmentation
$\varepsilon_\alpha$ as in Section~\ref{subsec:contact-algebras}.
Suppose that $\alpha$ has finitely many geometrically distinct
good least contractible closed Reeb orbits.
Denote them by $\gamma_1,\ldots,\gamma_r$ and assume that
$\widehat\mu(\gamma_i)\ne0$ for every $i$.

Let
\[
 b_m=\dim_{\mathbb Q}\LCH^0_m(M,\alpha;\varepsilon_\alpha)
\]
be the $m$th Betti number of the linearized contact homology
in the trivial free homotopy class.

Following \cite[p.~173]{HM15}, we define the positive and
negative mean Euler characteristics by
\begin{equation}\label{eq:mean-euler-definition}
 \chi_\pm(M,\alpha;\varepsilon_\alpha)=
 \lim_{N\to\infty}\frac{1}{N}\sum_{m=0}^{N}(-1)^m b_{\pm m}.
\end{equation}
We adapt the proof of \cite[Theorem~1.5]{HM15} to the trivial
free homotopy class.
Both limits in \eqref{eq:mean-euler-definition} then exist and
satisfy the resonance relations
\begin{equation}\label{eq:orbit-resonance}
 \chi_\pm(M,\alpha;\varepsilon_\alpha)=
 \sum_{\pm\widehat\mu(\gamma_i)>0}
 \frac{\widehat\chi(\gamma_i)}{|\widehat\mu(\gamma_i)|}.
\end{equation}
For the corrected sign in the negative relation, see
\cite[(3.1) and Remark~3.1]{AbM26}.

The $m$th Morse-type number
$c_m=\dim_{\mathbb Q}\LCC^0_m(M,\alpha;\varepsilon_\alpha)$
counts the good contractible orbits of Conley--Zehnder index $m$.
Grouping these orbits by their least contractible iterates gives
\[
 c_m=\sum_{i=1}^{r}\sum_{k=1}^{\infty}
     \dim_{\mathbb Q}\mathrm{HC}_m(\alpha,\gamma_i^k),
\]
as in \cite[p.~218]{HM15}.
By \eqref{eq:count-mean-estimate}, $c_m$ is bounded uniformly in $m$.

\begin{proposition}[Morse inequality]\label{prop:finite-morse}
For $a\le b$ with $a\equiv b\pmod 2$,
\begin{equation}\label{eq:count-finite-morse}
 \sum_{k=a}^{b}(-1)^{b-k}c_k
 \ge \sum_{k=a}^{b}(-1)^{b-k}b_k.
\end{equation}
\end{proposition}

\begin{proof}
Let
\[
 u_k=\operatorname{rank}\!\left(
 \partial_{\varepsilon_\alpha}\colon
 \LCC^0_k(M,\alpha;\varepsilon_\alpha)
 \longrightarrow
 \LCC^0_{k-1}(M,\alpha;\varepsilon_\alpha)
 \right).
\]
By rank-nullity,
\begin{equation}\label{eq:count-rank-nullity}
 c_k-b_k=u_k+u_{k+1}.
\end{equation}
Summing \eqref{eq:count-rank-nullity} with alternating signs over
any finite interval $a\le k\le b$ gives
\begin{equation}\label{eq:finite-morse-general}
 \sum_{k=a}^{b}(-1)^k(c_k-b_k)
 =(-1)^a u_a+(-1)^b u_{b+1}.
\end{equation}
When $a\equiv b\pmod 2$, multiplying
\eqref{eq:finite-morse-general} by $(-1)^b$ gives
$u_a+u_{b+1}\ge0$ on the right, proving \eqref{eq:count-finite-morse}.
\end{proof}

\section{Morse--Bott calculation}
\label{sec:local-bw-proof}

Throughout this section, we assume that
$c_1(\xi)|_{\pi_2(M)}=0$.
We perturb the connection form to construct the sub-DGAs
in Lemma~\ref{lem:local-bw-homology} from contractible
iterates $\gamma^k$ of simple fiber orbits $\gamma$
below a fixed action cutoff.
For each fixed $k$, the differential counting holomorphic
cylinders between the resulting orbits agrees with the
Morse differential on $B$.
The calculation also applies when $c_B=0$.

\Needspace{8\baselineskip}

\subsection{Perturbations of the connection form}
\label{subsec:finite-action-perturbations}

Ginzburg--G\"urel--Macarini use perturbations of the connection
form in their Morse--Bott computation for prequantization
bundles \cite[Section~3.4]{GGM18}; see also
\cite[Sections~6.1--6.2]{GGM18} for applications to linearized
contact homology. We recall the construction in the present
setting, following Bourgeois \cite[Section~2.2]{Bou}.

We use the prequantization circle bundle
\[
 \pi:(M^{2n+1},\xi)\longrightarrow(B^{2n},\omega)
\]
and the connection form $\beta$ fixed in
Section~\ref{subsec:intro-main-results}.

The restriction
\[
 d\pi|_{\xi_x}:(\xi_x,d\beta|_{\xi_x})
 \longrightarrow(T_{\pi(x)}B,\omega)
\]
is a symplectic isomorphism.

The horizontal lift $\widetilde Y$
of a vector field $Y$ on $B$ is determined by
$\beta(\widetilde Y)=0$ and $d\pi(\widetilde Y)=Y$.
Through the symplectic identification $d\pi|_\xi$,
an $\omega$-compatible almost complex structure on $TB$ induces
a $d\beta|_\xi$-compatible almost complex structure on $\xi$.

A horizontal lift of a smooth path $z:\mathbb R\to B$ is
a path $\widetilde{z}:\mathbb R\to M$ satisfying
\[
 \pi\circ\widetilde{z}=z,
 \qquad \beta(\widetilde{z}'(s))=0.
\]
A choice of $\widetilde{z}(0)\in\pi^{-1}(z(0))$ determines
a unique horizontal lift.
\par

A contact form $\alpha$ is of \emph{Morse--Bott type} if
$\operatorname{Spec}(\alpha)$ is discrete and, for every
$T\in\operatorname{Spec}(\alpha)$, the fixed-point set
$N_T=\{x\in M:\varphi_{R_\alpha}^T(x)=x\}$ is a closed smooth
submanifold, $d\alpha|_{TN_T}$ has locally constant rank, and
\[
 T_xN_T=\ker\bigl(d\varphi_{R_\alpha}^T(x)-\operatorname{id}\bigr)
 \qquad(x\in N_T);
\]
see \cite[Definition~1.7]{Bou}.
For the connection form $\beta$, we have
$\operatorname{Spec}(\beta)=\mathbb N$ and
$\varphi_{R_\beta}^\ell=\operatorname{id}_M$ for every
$\ell\in\mathbb N$. Thus $N_\ell=M$ and the Morse--Bott
condition holds. For each $\ell\in\mathbb N$, the $\ell$th
iterates of the simple fiber orbits have action $\ell$ and,
modulo time translation, are parametrized by $B$.

We use the minimal Chern number $c_B$ defined in
\eqref{eq:minimal-chern-number}.

\begin{lemma}\label{lem:topology}
The following assertions hold.
\begin{enumerate}
\item If $c_B>0$, then $\omega|_{\pi_2(B)}\ne0$.

\item Suppose that $\omega|_{\pi_2(B)}\ne0$, and let $h$
be the positive generator of $\omega(\pi_2(B))$. There is a unique
$\lambda\in\mathbb R$ such that
\[
 c_1(TB)[A]=\lambda\omega[A]
 \qquad\text{for all }A\in\pi_2(B).
\]
Writing $\sigma=\operatorname{sign}(\lambda)$, we have
$\lambda=\sigma c_B/h$, and $\sigma=0$ if and only if $c_B=0$.

Every simple closed Reeb orbit $\gamma$ of $\beta$ represents
an element of order $h$ in $\pi_1(M)$. For every $m\in\mathbb N$,
the iterate $\gamma^{mh}$ is contractible and satisfies
\[
 \mu_{\mathrm{RS}}(\gamma^{mh})=2\sigma m c_B,
\]
where the index is computed using a capping disk in $M$.
\end{enumerate}
\end{lemma}

\begin{proof}
Let $i:S^1\to M$ be the inclusion of a fiber.
The homotopy exact sequence of the circle bundle contains
\[
 \pi_2(M)\xrightarrow{\pi_*}\pi_2(B)
 \xrightarrow{\partial}\pi_1(S^1)\cong\mathbb Z
 \xrightarrow{i_*}\pi_1(M).
\]
The connecting map $\partial$ is evaluation on the Euler
class $-[\omega]$, up to the choice of generator of $\pi_1(S^1)$.
By exactness, a class in $\pi_2(B)$ lies in the image of $\pi_*$
precisely when its symplectic area is zero.
The identification $\xi\simeq\pi^*TB$, together with
$c_1(\xi)|_{\pi_2(M)}=0$, shows that $c_1(TB)$ vanishes on
every such class. If $\omega|_{\pi_2(B)}=0$, then
$c_1(TB)$ vanishes on all of $\pi_2(B)$ and hence $c_B=0$.
This proves \textup{(i)}.

For \textup{(ii)}, assume that $\omega|_{\pi_2(B)}\ne0$.
Since $[\omega]$ is integral, $\omega(\pi_2(B))=h\mathbb Z$ for
some $h\in\mathbb N$. Evaluation by $c_1(TB)$ factors through
$\omega$ and induces a homomorphism $h\mathbb Z\to\mathbb Z$ with
image $c_B\mathbb Z$. Its value on $h$ therefore has absolute value
$c_B$. Dividing this value by $h$ gives the unique proportionality
constant $\lambda$ and yields $\lambda=\sigma c_B/h$.

All simple fiber orbits of $\beta$ are freely homotopic.
Exactness gives $\ker i_*=h\mathbb Z$, so their classes
have order $h$ in $\pi_1(M)$.

Let $\gamma$ be a simple fiber orbit of $\beta$.
A capping disk $u$ for $\gamma^{mh}$ projects to
a sphere $A$ of area
$\int u^*d\beta=\int_{\partial u}\beta=mh$. In the horizontal
trivialization the linearized fiber flow is constant. Passing to the
capping trivialization adds $2c_1(TB)[A]=2\sigma m c_B$.
\end{proof}

For the remainder of this section, assume that
$\omega|_{\pi_2(B)}\ne0$ and let $h$ be as in
Lemma~\ref{lem:topology}.
The $m$th summand in \eqref{eq:BWpattern} corresponds
to the iterates $\gamma^{mh}$, where $\gamma$ ranges
over the simple fiber orbits of $\beta$.

Choose a smooth Morse function $f:B\to\mathbb R$ with pairwise
distinct critical values and $1\le f\le2$.
Let $J_B$ be an $\omega$-compatible almost complex structure
on $B$, and define the Riemannian metric
\[
 g_B(v,w)=\omega(v,J_Bw),
 \qquad v,w\in T_bB,\quad b\in B.
\]
By \cite[Theorem~8.1]{SZ92}, we can choose $J_B$ so that the
negative gradient flow of $f$ with respect to $g_B$ is
Morse--Smale, as defined in
Section~\ref{subsec:holomorphic-morse}.

For $\epsilon>0$, define $F_\epsilon:B\to\mathbb R_{>0}$ by
$F_\epsilon(b)=1+\epsilon f(b)$,
and set
\begin{equation}\label{eq:reference-perturbation}
 \beta_\epsilon=(F_\epsilon\circ\pi)\beta.
\end{equation}
Since $F_\epsilon>0$, $\beta_\epsilon$ is a contact form
with $\ker\beta_\epsilon=\xi$.
The computation in \cite[proof of Lemma~2.3]{Bou} gives
\begin{equation}\label{eq:reference-reeb-field}
 R_{\beta_\epsilon}
 =\frac{1}{F_\epsilon\circ\pi}\,R_\beta
  +\widetilde Y_\epsilon,
\end{equation}
where $Y_\epsilon$ is the vector field on $B$ determined by
\[
 \iota_{Y_\epsilon}\omega
 =\frac{1}{F_\epsilon^2}\,dF_\epsilon,
\]
and $\widetilde Y_\epsilon$ is its horizontal lift
with respect to $\beta$.
Equation~\eqref{eq:reference-reeb-field} shows that
$R_{\beta_\epsilon}$ projects to $Y_\epsilon$.
Since $dF_\epsilon=\epsilon\,df$ and $\omega$ is nondegenerate,
$Y_\epsilon$ vanishes precisely at the critical points of $f$.
Consequently, a fiber $\pi^{-1}(b)$ is invariant under the
Reeb flow of $\beta_\epsilon$ if and only if
$b\in\operatorname{Crit}(f)$.

Let $x\in\operatorname{Crit}(f)$. Since
$\left.\widetilde Y_\epsilon\right|_{\pi^{-1}(x)}=0$,
\eqref{eq:reference-reeb-field} gives
\[
 \left.R_{\beta_\epsilon}\right|_{\pi^{-1}(x)}
 =\frac{1}{1+\epsilon f(x)}
  \left.R_\beta\right|_{\pi^{-1}(x)}.
\]
Let $\gamma_{x,\beta_\epsilon}$ denote the simple
closed Reeb orbit of $\beta_\epsilon$ that traverses the
fiber $\pi^{-1}(x)$ once, considered up to time translation.
Since the flow of $R_\beta$ traverses each fiber once in
time $1$, the minimal period of $\gamma_{x,\beta_\epsilon}$
is $1+\epsilon f(x)$.
When only $\beta_\epsilon$ is under consideration, we write
$\gamma_x^\ell=\gamma_{x,\beta_\epsilon}^\ell$ for $\ell\ge1$.

Since $\gamma_{x}$ traverses its fiber once in the direction
of $R_\beta$, its class also has order $h$.
Thus $\gamma_{x}^\ell$ is contractible in $M$ if and only if
$h\mid\ell$.
The least contractible iterate of $\gamma_{x}$ is therefore
$\gamma_{x}^h$, and its contractible iterates are precisely
\[
 \gamma_{x}^{mh}=(\gamma_{x}^h)^m,\qquad m\ge1.
\]

\begin{lemma}\label{lem:local-bw-finite}
Fix a positive integer $N$. There is $\epsilon_0(N)>0$ such that,
for every $0<\epsilon<\epsilon_0(N)$,
$N+\tfrac12\notin\operatorname{Spec}(\beta_\epsilon)$,
and the closed Reeb orbits of $\beta_\epsilon$ with action
less than $N+\tfrac12$ are precisely
\[
 \gamma_{x}^\ell,\qquad
 1\le\ell\le N,\quad x\in\operatorname{Crit}(f).
\]
These orbits are nondegenerate and good, and their actions are
\begin{equation}\label{eq:reference-orbit-action}
 \mathcal A_{\beta_\epsilon}(\gamma_{x}^\ell)
 =\ell(1+\epsilon f(x)).
\end{equation}
For $x\in\operatorname{Crit}(f)$ and each positive integer $m$
with $mh\le N$, the Conley--Zehnder index of the contractible
orbit $\gamma_{x}^{mh}$, computed using a capping disk as in
Section~\ref{subsec:index-iteration}, is
\begin{equation}\label{eq:BWindex}
 \mu(\gamma_{x}^{mh})
 =2\sigma m c_B-n+\operatorname{ind}_f(x).
\end{equation}
Here $\sigma$ is as in Lemma~\ref{lem:topology}, and
$\operatorname{ind}_f(x)$ is the Morse index defined in
Section~\ref{subsec:holomorphic-morse}.
\end{lemma}

\begin{proof}
For each $1\le\ell\le N$, the space of $\ell$th iterates of
the simple fiber orbits of $\beta$, modulo time translation,
is $B$.
We use the same Morse function $f$ on each of these orbit
spaces. By \cite[Lemma~2.3]{Bou}, there is
$\epsilon_0(N)>0$ such that, for every
$0<\epsilon<\epsilon_0(N)$, the closed Reeb orbits of
$\beta_\epsilon$ with action at most $N+\tfrac12$ are
precisely $\gamma_{x}^\ell$, where $1\le\ell\le N$ and
$x\in\operatorname{Crit}(f)$, and all these orbits are
nondegenerate.
Decrease $\epsilon_0(N)$ so that
$2N\epsilon_0(N)<\tfrac12$.
Since $f\le2$, \eqref{eq:reference-orbit-action} gives
\[
 \mathcal A_{\beta_\epsilon}(\gamma_{x}^\ell)
 \le N(1+2\epsilon)<N+\tfrac12.
\]
Thus $N+\tfrac12\notin\operatorname{Spec}(\beta_\epsilon)$.

For each $x\in\operatorname{Crit}(f)$, choose a symplectic
basis of $T_{x}B$. The symplectic identification
$d\pi|_\xi$ gives a trivialization of $\xi$ along $\gamma_{x}$.
Use the induced trivializations along its iterates.
In these trivializations, the linearized unperturbed
fiber flow is the constant identity path, whose
Robbin--Salamon index is zero.
Since the orbit space $B$ has dimension $2n$,
\cite[Lemma~2.4]{Bou} gives the Conley--Zehnder index
$-n+\operatorname{ind}_f(x)$ for $\gamma_{x}^\ell$ in
these trivializations, for every $1\le\ell\le N$,
after decreasing $\epsilon_0(N)$ if necessary.
For each fixed $x$, these indices are independent of
$\ell$. Hence every $\gamma_{x}^\ell$ in this range has
the same index parity as $\gamma_{x}$ and is good by
the criterion in
Section~\ref{subsec:index-iteration};
see also \cite[proof of Lemma~7.5]{Bou}.

For a contractible iterate $\gamma_{x}^{mh}$ with
$mh\le N$, changing from the horizontal trivialization
to the trivialization induced by a capping disk
adds $2\sigma m c_B$ to the index, as computed in
Lemma~\ref{lem:topology}.
This proves \eqref{eq:BWindex}.
\end{proof}

For a contact form $\widetilde\beta$ supporting $\xi$
sufficiently $C^\infty$-close to $\beta_\epsilon$ and for
$x\in\operatorname{Crit}(f)$, let $\gamma_{x,\widetilde\beta}$
be the unique nearby simple closed orbit of $R_{\widetilde\beta}$,
up to time translation, whose period is close to that of
$\gamma_{x,\beta_\epsilon}$.
We call $\gamma_{x,\widetilde\beta}$ the \emph{continuation}
of $\gamma_{x,\beta_\epsilon}$ to $\widetilde\beta$;
cf. \cite[Section~1.3, p.~8]{dPHKS24}.
Here $x$ labels the original fiber orbit;
$\gamma_{x,\widetilde\beta}$ need not lie in $\pi^{-1}(x)$.

\begin{corollary}\label{cor:finite-orbit-persistence}
Fix $N$ and $0<\epsilon<\epsilon_0(N)$ as in
Lemma~\ref{lem:local-bw-finite}. There is a $C^\infty$
neighborhood $\mathcal U$ of $\beta_\epsilon$ among contact
forms supporting $\xi$ such that the following
statements hold for every $\widetilde\beta\in\mathcal U$.
We have $N+\tfrac12\notin\operatorname{Spec}({\widetilde\beta})$,
and the closed Reeb orbits of ${\widetilde\beta}$ with action less than
$N+\tfrac12$ are precisely the nondegenerate iterates
$\gamma_{x,{\widetilde\beta}}^\ell$, where $1\le\ell\le N$ and
$x\in\operatorname{Crit}(f)$.

For each $1\le\ell\le N$, the orbit $\gamma_{x,{\widetilde\beta}}^\ell$
is the continuation of $\gamma_{x,\beta_\epsilon}^\ell$ and is freely homotopic
to it. All these orbits are good.
The Conley--Zehnder indices of the contractible orbits,
computed with capping disks, are still given by
\eqref{eq:BWindex}.
The actions of $\gamma_{x,\widetilde\beta}^\ell$ have
the same strict ordering as those of $\gamma_{x,\beta_\epsilon}^\ell$.
\end{corollary}

\begin{proof}
By Lemma~\ref{lem:local-bw-finite}, all the orbits
$\gamma_{x,\beta_\epsilon}^\ell$ in the stated range
are nondegenerate.
Applying the implicit function theorem to the corresponding
local Poincar\'e return maps gives unique nearby nondegenerate
closed orbits, up to time translation, with periods close
to the original ones.
Their periods depend continuously on $\widetilde\beta$.
After choosing time origins and parametrizing the orbits
on $\mathbb R/\mathbb Z$, the orbits also depend continuously
on $\widetilde\beta$ in the $C^\infty$ topology.
Since there are only finitely many such orbits, we can
choose a common neighborhood $\mathcal U$ in which their
continuations remain nondegenerate and have action less
than $N+\tfrac12$.
By local uniqueness, the continuation of
$\gamma_{x,\beta_\epsilon}^\ell$ is
$\gamma_{x,\widetilde\beta}^\ell$, the $\ell$th iterate
of the continued simple orbit.

We claim that, after shrinking $\mathcal U$, there are no
other closed orbits of action at most $N+\tfrac12$.
Otherwise, choose contact forms ${\widetilde\beta}_k\to\beta_\epsilon$
and additional closed orbits with periods $T_k\le N+\tfrac12$.
The vector fields $R_{{\widetilde\beta}_k}$ converge smoothly to the
nonvanishing vector field $R_{\beta_\epsilon}$, so their
closed-orbit periods have a uniform positive lower bound.
After passing to a subsequence, the periods converge to
$T\in(0,N+\tfrac12]$, and the periodic trajectories converge
to a closed orbit of $R_{\beta_\epsilon}$ of period $T$.
By Lemma~\ref{lem:local-bw-finite}, this limiting orbit is one
of the $\gamma_{x,\beta_\epsilon}^\ell$ and has period less than
$N+\tfrac12$. Local uniqueness then implies that, for all
sufficiently large $k$, the chosen orbits are its
continuations, contrary to their choice.

Each orbit is freely homotopic to its continuation.
Transport a symplectic trivialization of $\xi$ along each
$\gamma_{x,\beta_\epsilon}$ to $\gamma_{x,{\widetilde\beta}}$ and use the induced
trivializations along their iterates.
In these trivializations, the Conley--Zehnder indices remain
constant because the orbits remain nondegenerate.
Their parities are therefore unchanged, so the continued
iterates remain good.
For a contractible orbit, transport its capping disk along
the continuation. The Conley--Zehnder index computed in the
induced trivialization is also constant and hence is still
given by \eqref{eq:BWindex}.

At ${\widetilde\beta}=\beta_\epsilon$, the actions are pairwise distinct
by \eqref{eq:reference-orbit-action}, the distinct critical
values of $f$, and $2N\epsilon<\tfrac12$.
There are only finitely many such actions, so their strict
ordering persists by continuity after shrinking $\mathcal U$.
\end{proof}

Given $0<\delta<\tfrac12$, choose $0<\epsilon<\epsilon_0(N)$
with $2N\epsilon<\delta$. Since $1\le f\le2$,
\eqref{eq:reference-orbit-action} gives
$\ell<\mathcal A_{\beta_\epsilon}(\gamma_{x,\beta_\epsilon}^\ell)<\ell+\delta$
for $1\le\ell\le N$.
After shrinking the neighborhood in
Corollary~\ref{cor:finite-orbit-persistence}, we obtain
\begin{equation}\label{eq:finite-layer-action-bounds}
 \ell<\mathcal A_{\widetilde\beta}(\gamma_{x,{\widetilde\beta}}^\ell)
 <\ell+\delta,
 \qquad 1\le\ell\le N,\quad x\in\operatorname{Crit}(f),
 \quad{\widetilde\beta}\in\mathcal U.
\end{equation}
We take $\delta=\tfrac14$ in the subsequent arguments.
The choices of $\epsilon$ and $\mathcal U$ may depend on $N$.

\subsection{Holomorphic cylinders and the Morse differential}

Fix $N$, and let $\beta_\epsilon$ be the contact form of
Lemma~\ref{lem:local-bw-finite}, with $\epsilon>0$ sufficiently
small. Let $J_\epsilon$ be the cylindrical almost complex
structure for $\beta_\epsilon$, in the sense of
Section~\ref{subsec:holomorphic-morse}, with
$J_\epsilon|_\xi$ equal to the horizontal lift of $J_B$
with respect to $\beta$.
Since $J_B$ is $\omega$-compatible and
$d\beta_\epsilon|_\xi=(F_\epsilon\circ\pi)d\beta|_\xi$
with $F_\epsilon>0$, $J_\epsilon|_\xi$ is compatible
with $d\beta_\epsilon|_\xi$.

Fix $1\le\ell\le N$, and let $z:\mathbb R\to B$ be a
nonconstant Morse trajectory from $x$ to $y$ for $(f,g_B)$,
as in \eqref{eq:morse-trajectory}.
Following \cite[Section~5.3.2, pp.~61--62]{Bou}, we construct
from $z$ a finite-energy $J_\epsilon$-holomorphic cylinder
\[
 u=(a,v):\mathbb R\times(\mathbb R/\mathbb Z)
 \longrightarrow \mathbb R\times M
\]
with positive asymptotic orbit $\gamma_{x}^\ell$ and negative
asymptotic orbit $\gamma_{y}^\ell$ as follows.

We use a decreasing change of parameter $r=r(s)$ along the
Morse trajectory $z$. We choose the speed
so that the cylinder defined below satisfies the holomorphic
equation:
\begin{equation}\label{eq:lg-parameter}
 r'(s)=-\frac{\ell\epsilon}{F_\epsilon(z(r(s)))},
 \qquad r(0)=0.
\end{equation}
For fixed $\epsilon>0$, the bounds $1\le f\le2$ imply that
$-r'(s)$ is bounded above and below by positive constants.
Hence $r$ is defined on all of $\mathbb R$ and is a decreasing
diffeomorphism $\mathbb R\to\mathbb R$, with
$r(s)\to\mp\infty$ as $s\to\pm\infty$.

Choose a horizontal lift $\widetilde{z}:\mathbb R\to M$
of $z$ with respect to $\beta$, as in
Section~\ref{subsec:finite-action-perturbations}.
Define $u(s,t)=(a(s),v(s,t))$ by
\begin{equation}\label{eq:lg-cylinder}
 \begin{aligned}
 a(s)&=\ell\int_0^s F_\epsilon(z(r(\sigma)))\,d\sigma,\\
 v(s,t)&=\varphi_{R_\beta}^{\ell t}
          (\widetilde{z}(r(s))),
 \qquad t\in\mathbb R/\mathbb Z.
 \end{aligned}
\end{equation}
Since the Reeb flow of $\beta$ traverses each fiber once
in time $1$, $v$ is well defined for
$t\in\mathbb R/\mathbb Z$.
It satisfies $\pi(v(s,t))=z(r(s))$, and each loop
$v(s,\cdot)$ traverses the fiber over $z(r(s))$
exactly $\ell$ times.

We verify \eqref{eq:holomorphic-cylinder-equation} with
$J=J_\epsilon$.
The flow of $R_\beta$ preserves $\beta$, and
$\widetilde{z}$ is horizontal, so $v_s\in\xi$.
Also, $v_t=\ell R_\beta(v(s,t))$.
With respect to the splitting
\[
 T(\mathbb R\times M)
 =\mathbb R\partial_a\oplus\mathbb R R_{\beta_\epsilon}
   \oplus\xi,
\]
the equations in the $\partial_a$ and $R_{\beta_\epsilon}$
directions follow from \eqref{eq:lg-cylinder}:
\[
 a_s=\ell F_\epsilon(z(r(s)))=\beta_\epsilon(v_t),
 \qquad a_t=0=-\beta_\epsilon(v_s).
\]

By \eqref{eq:reference-reeb-field}, the $\xi$-component
of $v_t$ in this splitting is
$-\ell F_\epsilon(z(r(s)))\widetilde Y_\epsilon(v(s,t))$.
The definitions of $Y_\epsilon$ and $g_B$ give
$J_BY_\epsilon=F_\epsilon^{-2}\nabla_{g_B}F_\epsilon$.
Together with \eqref{eq:morse-trajectory} and
\eqref{eq:lg-parameter}, this yields
\begin{equation}\label{eq:lg-gradient}
\begin{aligned}
 \frac{d}{ds}z(r(s))
 &=-r'(s)\nabla_{g_B}f(z(r(s)))\\
 &=\frac{\ell\epsilon}{F_\epsilon(z(r(s)))}
   \nabla_{g_B}f(z(r(s)))\\
 &=\ell F_\epsilon(z(r(s)))
   J_B\bigl(Y_\epsilon(z(r(s)))\bigr).
\end{aligned}
\end{equation}
We have $d\pi(v_s)=\frac{d}{ds}z(r(s))$.
Since $d\pi|_\xi$ is an isomorphism intertwining
$J_\epsilon|_\xi$ and $J_B$, \eqref{eq:lg-gradient} gives
\[
 v_s-\ell F_\epsilon(z(r(s)))
 J_\epsilon\bigl(\widetilde Y_\epsilon(v(s,t))\bigr)=0.
\]
This verifies the remaining $\xi$-component of the
holomorphic equation, so $u$ is $J_\epsilon$-holomorphic.

The Morse trajectory $z$ converges exponentially to $x$ and
$y$ as $r\to-\infty$ and $r\to+\infty$, respectively.
It therefore has finite length near both ends, so its
horizontal lift $\widetilde{z}$ converges to points in the
corresponding fibers.
Since $r(s)\to\mp\infty$ and $a(s)\to\pm\infty$ as
$s\to\pm\infty$, it follows from \eqref{eq:lg-cylinder}
that $u$ has positive asymptotic orbit $\gamma_{x}^\ell$
and negative asymptotic orbit $\gamma_{y}^\ell$.
By \eqref{eq:symplectization-action-energy}, its Hofer
energy satisfies $E(u)\le2\ell F_\epsilon(x)<\infty$.
If $z(r)\equiv x$, \eqref{eq:lg-cylinder} gives a trivial
cylinder over $\gamma_{x}^\ell$.

We use the following consequence of \cite[Proposition~4.16]{Bou}.
Convergence is understood in the sense of
\cite[Definitions~4.12--4.15]{Bou}.

\Needspace{7\baselineskip}
\begin{lemma}\label{lem:local-bw-low-index}
Fix a positive integer $N$, and let $\epsilon_0(N)>0$ be as in
Lemma~\ref{lem:local-bw-finite}. There exists
$0<\epsilon_*(N)\le\epsilon_0(N)$ such that the following hold.
\begin{enumerate}
\item Fix $1\le\ell\le N$ and distinct critical points
$x,y$ of $f$. Let $\epsilon_\nu\in(0,\epsilon_*(N))$ tend to
zero, and let $u_\nu$ be finite-energy
$J_{\epsilon_\nu}$-holomorphic cylinders with positive end
$\gamma_{x}^\ell$ and negative end $\gamma_{y}^\ell$.
After passing to a subsequence, $u_\nu$ converges to a limit
whose holomorphic components are trivial cylinders over the fibers.

The gradient-flow segments concatenate to a possibly broken
Morse trajectory from $x$ to $y$, with $k\ge1$ nonconstant
components. If $\operatorname{ind}_f(x)-\operatorname{ind}_f(y)=1$,
then $k=1$.

\item Let $0<\epsilon<\epsilon_*(N)$,
$1\le\ell\le N$, and $x,y\in\operatorname{Crit}(f)$.
If there exists a nontrivial finite-energy
$J_\epsilon$-holomorphic cylinder with positive end
$\gamma_{x}^\ell$ and negative end $\gamma_{y}^\ell$, then
\[
 \operatorname{ind}_f(x)-\operatorname{ind}_f(y)\ge1.
\]
If, moreover, $\ell=mh$ for some positive integer $m$,
then the asymptotic orbits are contractible and
\[
 \mu(\gamma_{x}^{mh})-\mu(\gamma_{y}^{mh})
 =\operatorname{ind}_f(x)-\operatorname{ind}_f(y)\ge1.
\]
\end{enumerate}
\end{lemma}

\begin{proof}
We first prove \textup{(i)} for any sequence
$\epsilon_\nu\in(0,\epsilon_0(N))$ tending to zero.
The actions
$\mathcal A_{\beta_{\epsilon_\nu}}(\gamma_{x}^\ell)
=\ell(1+\epsilon_\nu f(x))$
are uniformly bounded.
Thus \eqref{eq:symplectization-action-energy} gives a
uniform bound on the Hofer energies $E(u_\nu)$.
Moreover,
\[
 0\le\int u_\nu^*d\beta_{\epsilon_\nu}
   =\ell\epsilon_\nu\bigl(f(x)-f(y)\bigr)\longrightarrow0.
\]
By \cite[Proposition~4.16]{Bou}, a subsequence converges
in the sense specified above.
Every holomorphic component of the limit has zero
$d\beta$-area. Since $d\beta=\pi^*\omega$ and its
projection to $B$ is $J_B$-holomorphic, this projection
is constant.
Planes are excluded because their $d\beta$-area equals
the positive action of their asymptotic orbit.
The symplectization is exact, so nonconstant sphere
bubbles cannot occur.
After omitting constant components, a component with
more than two punctures would force a plane component,
since the limit has genus zero and only one positive
and one negative end.
Thus the holomorphic components are trivial cylinders
over the fibers, connected in a chain by gradient segments.
The covering multiplicity is preserved along this chain,
so every cylinder has multiplicity $\ell$.

Following the chain from the positive end to the negative
end, we parametrize each gradient segment in $B$ so that
$z'=-\nabla_{g_B}f(z)$.
The two ends of each trivial cylinder project to the
same point of $B$. If adjacent gradient segments meet
there at finite values of their parameters, uniqueness
for the gradient equation joins them into a single
trajectory. Breaks can occur only at critical points
approached in infinite time.
The gradient trajectories attached to the positive and
negative ends of the chain converge to $x$ and $y$,
respectively. Including these trajectories, we obtain
a possibly broken Morse trajectory from $x$ to $y$;
see \cite[proof of Proposition~4.16, p.~45]{Bou}.

Let $k$ be the number of its nonconstant components.
Since $x\ne y$, we have $k\ge1$.
Since $(f,g_B)$ is Morse--Smale, each nonconstant
Morse trajectory decreases the Morse index by at
least one; see
Section~\ref{subsec:holomorphic-morse}.
Hence
\[
 1\le k\le
 \operatorname{ind}_f(x)-\operatorname{ind}_f(y).
\]
If the index difference is one, then $k=1$, so the
trajectory is unbroken. This proves \textup{(i)}.

For \textup{(ii)}, a cylinder $u$ with $x=y$ satisfies
$\int u^*d\beta_\epsilon=0$.
Hence its image lies in $\mathbb R\times\gamma_{x}$,
and $u$ is a trivial cylinder, possibly multiply covered.

If no uniform $\epsilon_*(N)>0$ existed, we could
choose $\epsilon_\nu\in(0,\epsilon_0(N))$ tending to
zero and nontrivial cylinders $u_\nu$ violating the
Morse index inequality in \textup{(ii)}.
Since $1\le\ell\le N$ and $\operatorname{Crit}(f)$ is
finite, after passing to a subsequence we may assume
that their positive and negative ends are
$\gamma_{x}^\ell$ and $\gamma_{y}^\ell$, respectively,
for fixed $\ell,x,y$ with $x\ne y$ and
\[
 \operatorname{ind}_f(x)\le\operatorname{ind}_f(y).
\]
The argument for \textup{(i)} gives
$\operatorname{ind}_f(x)>\operatorname{ind}_f(y)$,
a contradiction.
Thus there exists $0<\epsilon_*(N)\le\epsilon_0(N)$
such that the Morse index inequality in \textup{(ii)}
holds for every $0<\epsilon<\epsilon_*(N)$.
When $\ell=mh$, the asymptotic orbits are contractible,
and \eqref{eq:BWindex} gives the claimed
Conley--Zehnder index difference.
\end{proof}

Let $\ell=mh\le N$, where $m$ is a positive integer,
and let $x,y\in\operatorname{Crit}(f)$ be distinct.
Denote by $\mathcal M_{\ell,\epsilon}(x,y)$ the moduli
space of finite-energy $J_\epsilon$-holomorphic cylinders
from $\gamma_{x}^\ell$ to $\gamma_{y}^\ell$.
We include all relative homotopy classes and allow
multiply covered cylinders.
We use the asymptotic markers and equivalence relation
specified in Section~\ref{subsec:holomorphic-morse}.
By \eqref{eq:cylinder-moduli-dimension} and
\eqref{eq:BWindex}, the expected dimension is
\[
 \operatorname{ind}_f(x)-\operatorname{ind}_f(y)-1.
\]
We use the Morse moduli space $\mathcal M_f(x,y)$ and
its orientation from
Section~\ref{subsec:holomorphic-morse}, with $g=g_B$.

\Needspace{12\baselineskip}
\begin{proposition}\label{prop:local-bw-counts}
Fix $N$. For all sufficiently small $\epsilon>0$,
if $\ell=mh\le N$ and $x,y\in\operatorname{Crit}(f)$
satisfy
\[
 \operatorname{ind}_f(x)-\operatorname{ind}_f(y)=1,
\]
then every cylinder in $\mathcal M_{\ell,\epsilon}(x,y)$
is regular.
With coherent orientations compatible with the Morse
orientations, this moduli space is a compact oriented
zero-dimensional manifold, and
\begin{equation}\label{eq:local-cylinder-count}
 \#\mathcal M_{\ell,\epsilon}(x,y)
 =\ell\,\#\mathcal M_f(x,y)
 =\ell\,n_f(x,y).
\end{equation}
Here the cylinder moduli space uses the asymptotic markers
specified in Section~\ref{subsec:holomorphic-morse}.
\end{proposition}

\begin{proof}
Fix $\ell=mh\le N$ and critical points $x,y$ satisfying
\[
 \operatorname{ind}_f(x)-\operatorname{ind}_f(y)=1.
\]
We use the comparison argument in
\cite[Section~8.1]{Bou}, restricted to the
action range of Lemma~\ref{lem:local-bw-finite}.
Throughout the proof, we take
$0<\epsilon<\epsilon_*(N)$ and decrease the upper bound
on $\epsilon$ as needed.

For every covering multiplicity, the unperturbed orbit
space is $B$. The restriction $J_\epsilon|_\xi$ is
independent of $\epsilon$, and the limiting almost complex
structure is invariant under fiber rotation.
Thus our family has the form considered in
\cite[Sections~2.1 and~4.2]{Bou}.
By Lemma~\ref{lem:local-bw-low-index}\textup{(i)}, the
holomorphic components of any limit as $\epsilon\to0$
are trivial fiber cylinders, and the gradient segments
concatenate to an unbroken Morse trajectory from $x$ to $y$.
These trajectories represent elements of
$\mathcal M_f(x,y)$, which is compact and zero-dimensional
and hence finite.

We first omit the asymptotic markers and retain the
finite domain automorphism groups.
At $\epsilon=0$, linearize at a trivial cylinder over
the fiber at $b\in B$, allowing the asymptotic orbits
at the two ends to vary independently in $B$.
Using $d\pi|_\xi$, identify horizontal sections with
$T_bB$-valued functions. The horizontal linearized
operator is
\[
 D_b=\partial_s+J_B(b)\partial_t.
\]
Choose a positive exponential weight smaller than the
absolute value of every nonzero eigenvalue of
$J_B(b)\partial_t$, and allow sections to converge
exponentially to independently varying constants
at the two ends.
On the nonzero Fourier modes, $D_b$ is invertible.
On the zero mode, it is $\partial_s$.
Integrating an exponentially decaying right-hand side
gives a section with the required asymptotics, so
$\partial_s$ is surjective, with kernel consisting
of the constant sections, identified with $T_bB$.
The operator in the directions spanned by $\partial_a$
and $R_\beta$ is likewise surjective when asymptotic
translations and phases are included.
The same calculation applies to covered fiber cylinders,
since their horizontal linearized Reeb flow is the identity.

Modulo reparametrization, the family of trivial fiber
cylinders is $B$, and the pair of asymptotic evaluation
maps is
\[
 B\longrightarrow B\times B,\qquad b\longmapsto(b,b).
\]
Its transversality to $W^u(x)\times W^s(y)$ is precisely
the Morse--Smale condition
\[
 T_bW^u(x)+T_bW^s(y)=T_bB,
 \qquad b\in W^u(x)\cap W^s(y).
\]
Each individual endpoint evaluation is the identity on
$B$. Inserting a gradient segment composes an evaluation
map with the negative gradient flow, which preserves the
relevant stable and unstable manifolds.
Thus the matching and endpoint conditions are transverse,
as required in \cite[Section~6.1]{Bou}.

By
\cite[Proposition~5.14 and Corollary~5.15]{Bou},
the preceding surjectivity and transversality give right
inverses for the linearized operators in the gluing
construction that are bounded uniformly as $\epsilon\to0$,
in the norms used there.
The implicit function argument in
\cite[Section~5.3.3]{Bou} then shows that all solutions
near the approximate maps are regular.
Together with the charts in
\cite[Lemma~6.3]{Bou}, this gives neighborhoods
of the limiting configurations in the parameterized
moduli space.

These neighborhoods contain every cylinder under
consideration for sufficiently small $\epsilon$.
Otherwise, a sequence outside them with $\epsilon\to0$
would, by Lemma~\ref{lem:local-bw-low-index}\textup{(i)},
have a subsequence converging to one of the unbroken
Morse trajectories, contrary to the local description.
Hence every cylinder in
$\mathcal M_{\ell,\epsilon}(x,y)$ is regular.

We next check compactness of the parameterized family,
using \cite[proof of Theorem~2.27]{Pardon}.
For $2N\epsilon<\tfrac12$,
\eqref{eq:symplectization-action-energy} gives the uniform
bound $E(u)<2N+1$. Moreover,
\[
 0\le\int u^*d\beta_\epsilon
   =\ell\epsilon\bigl(f(x)-f(y)\bigr)
   \le N\epsilon<1.
\]
Suppose that $\epsilon_\nu\to\epsilon_\infty>0$.
Every asymptotic orbit in the limiting building is an orbit
of $\beta_{\epsilon_\infty}$ with action less than $N+\tfrac12$,
so Lemma~\ref{lem:local-bw-finite} applies.
In particular, its action is greater than $1$.
A plane would have $d\beta_{\epsilon_\infty}$-area equal to its
positive asymptotic action, so no plane can occur.
Exactness excludes nonconstant sphere bubbles.
As in the proof of
Lemma~\ref{lem:local-bw-low-index}, a genus-zero limit
with one positive and one negative end must therefore
be a chain of cylinders.

By \eqref{eq:reference-orbit-action}, the actions of
orbits of multiplicity $k\le N$ lie in
$(k,k+\tfrac12)$.
The action of every intermediate orbit lies between
those of the original ends, so its multiplicity is
$\ell$.
Each nontrivial component strictly decreases the Morse
index by Lemma~\ref{lem:local-bw-low-index}\textup{(ii)}.
Since the total decrease is one, there is only one
nontrivial component.
Thus the compactification adds no broken cylinders.

At parameter zero, the limits are those described
in Lemma~\ref{lem:local-bw-low-index}\textup{(i)}.
Thus, for every sufficiently small $\epsilon>0$,
the family over $[0,\epsilon]$, compactified by these
limits at zero, is compact and has boundary only
at zero and $\epsilon$.

We choose coherent orientations compatible with the
Morse orientations fixed in
Section~\ref{subsec:holomorphic-morse}.
The orbit space $B$ is symplectically oriented.
The iterates of the simple fiber orbits of $\beta$ are good
in the sense of \cite[Definition~7.1]{Bou}: their horizontal
linearized Reeb flows are constant identity paths,
and the orbit space is $B$ for every iterate.
Along any loop $c:S^1\to B$, choose a unitary trivialization
of $c^*TB$. Pulling this trivialization back along
$c^*M\to S^1$ makes the horizontal asymptotic operators
constant, so their orientation lines have trivial monodromy.
The orientation hypotheses of
\cite[Lemma~7.2]{Bou} are therefore satisfied.
By \cite[Lemma~7.6]{Bou}, the comparison preserves the
resulting Morse orientations; the additional sign
exponents there vanish for one positive and one
negative puncture.

We now restore the asymptotic markers, using the conventions
of Section~\ref{subsec:holomorphic-morse}.
For the gradient segments, we use the matching
condition in \cite[Definition~4.15(iii)]{Bou},
which in our setting is given by horizontal lifts
with respect to $\beta$.

In a local orbifold chart $U/G$ of the parameterized
moduli space, let $S$ denote the finite set of pairs of
asymptotic markers over one point, with both markers
mapping to the prescribed basepoints.
After shrinking $U$, the chart with markers is
\[
 (U\times S)/G,
\]
with the action on the rays induced by derivatives
at the punctures.
The gluing maps and changes of chart respect these
actions and the markers at the external punctures.
Thus these charts define the cobordism with both
asymptotic markers retained.

The marker choices are discrete and do not change the
linearized operators. They also eliminate finite domain
automorphisms.
All relevant perturbed orbits are good by
Lemma~\ref{lem:local-bw-finite}.
Consequently, the $\mathbb Z/\ell$-action on the marker
choices at either end preserves the orientation; see
\cite[Definition~2.49]{Pardon}.
Thus $\mathcal M_{\ell,\epsilon}(x,y)$ is a compact
oriented zero-dimensional manifold, and the cobordism
with markers is oriented.

It remains to count the marker choices over each
limiting Morse trajectory.
Consider the corresponding cylinder
\eqref{eq:lg-cylinder}.
There are $\ell$ choices of marker at each end,
hence $\ell^2$ pairs.
The deck group is $\mathbb Z/\ell$, generated by
\[
 (s,t)\longmapsto(s,t+1/\ell).
\]
It acts freely on the pairs of markers, leaving
$\ell$ inequivalent choices.
These are the marker choices at the corresponding
$\epsilon=0$ boundary of the cobordism, and each has
the sign of the Morse trajectory.
The oriented cobordism therefore gives
\[
 \#\mathcal M_{\ell,\epsilon}(x,y)
 =\ell\,\#\mathcal M_f(x,y)
 =\ell\,n_f(x,y),
\]
which proves \eqref{eq:local-cylinder-count}.

Since there are only finitely many triples $(\ell,x,y)$ satisfying
the hypotheses, the upper bound on $\epsilon$ can be chosen uniformly.
\end{proof}

\Needspace{11\baselineskip}
\subsection{The contact differential}
\label{subsec:layer-complexes}

To apply the preceding cylinder calculation to the
contact algebra of Section~\ref{subsec:contact-algebras},
we first perturb $\beta_\epsilon$ to a globally
nondegenerate contact form.
Lemma~\ref{lem:local-bw-finite} gives nondegeneracy only
below the action bound $N+\tfrac12$.
The following lemma shows that the signed counts in
Proposition~\ref{prop:local-bw-counts} are preserved under
sufficiently small perturbations.

\begin{lemma}\label{lem:local-bw-perturbation-persistence}
Fix $N$ and a positive $\epsilon$ for which
Lemma~\ref{lem:local-bw-low-index} holds. Suppose that every
finite-energy $J_\epsilon$-holomorphic cylinder from
$\gamma_{x,\beta_\epsilon}^{mh}$ to $\gamma_{y,\beta_\epsilon}^{mh}$ with $mh\le N$ and
Conley--Zehnder index difference one is regular.

Let $\mathcal U$ be the neighborhood of $\beta_\epsilon$
in Corollary~\ref{cor:finite-orbit-persistence}.
For $\widetilde\beta\in\mathcal U$, let $J$ be a cylindrical
almost complex structure for $\widetilde\beta$.
If $(\widetilde\beta,J)$ is sufficiently $C^\infty$-close to
$(\beta_\epsilon,J_\epsilon)$ in the sense of
\eqref{eq:contact-pair-topology}, then every nontrivial
finite-energy $J$-holomorphic cylinder from
$\gamma_{x,\widetilde\beta}^{mh}$ to
$\gamma_{y,\widetilde\beta}^{mh}$, with $mh\le N$, has
positive Conley--Zehnder index difference.
When this difference is one, the cylinders are regular,
and their moduli spaces have the same signed counts as
those for $(\beta_\epsilon,J_\epsilon)$.
\end{lemma}

\begin{proof}
By Corollary~\ref{cor:finite-orbit-persistence}, the continued
orbits have the same covering multiplicities and Conley--Zehnder
indices, and their actions retain the same strict ordering.
If $x=y$, any such cylinder has zero $d\widetilde\beta$-area
and is trivial.

Consider nontrivial finite-energy cylinders $u_\nu$ for
$(\widetilde\beta_\nu,J_\nu)\to(\beta_\epsilon,J_\epsilon)$,
with both ends of the same covering multiplicity $mh\le N$.
After passing to a subsequence, we may fix $m,x,y$ with
$x\ne y$ so that their positive and negative ends are
$\gamma_{x,\widetilde\beta_\nu}^{mh}$ and
$\gamma_{y,\widetilde\beta_\nu}^{mh}$, respectively.
Their actions are less than $N+\tfrac12$, so
\eqref{eq:symplectization-action-energy} gives
$E(u_\nu)<2N+1$.
The compactness argument in \cite[proof of Theorem~2.27]{Pardon}
applies to this convergent sequence of contact forms and
almost complex structures: all asymptotic orbits lie below
the action bound $N+\tfrac12$, where $\beta_\epsilon$ is
nondegenerate.
Thus a subsequence converges to a $J_\epsilon$-holomorphic
building.

As in the proof of Proposition~\ref{prop:local-bw-counts},
the area bound excludes planes, and exactness excludes
nonconstant sphere bubbles. The domains have no nodes;
see \cite[Remark~2.8]{Pardon}.
The limit is therefore a chain of cylinders, and the action
estimates in that proof show that every intermediate orbit
has covering multiplicity $mh$.
We use the stable compactification of \cite[Remark~2.9]{Pardon},
in which these cylinders are nontrivial.

Let $u^1,\ldots,u^k$ be the cylinders in this chain.
Since $x\ne y$, we have $k\ge1$.
Lemma~\ref{lem:local-bw-low-index} gives
$\operatorname{ind}(u^i)\ge1$ for every $i$, and hence
\[
\begin{aligned}
 \mu(\gamma_{x,\widetilde\beta_\nu}^{mh})
 -\mu(\gamma_{y,\widetilde\beta_\nu}^{mh})
 &=\mu(\gamma_{x,\beta_\epsilon}^{mh})-\mu(\gamma_{y,\beta_\epsilon}^{mh})\\
 &=\sum_{i=1}^k\operatorname{ind}(u^i)
 \ge k\ge1.
\end{aligned}
\]
This rules out nontrivial cylinders of nonpositive index
for arbitrarily small perturbations.
If the $u_\nu$ have index one, then $k=1$, and the limit
is a single index-one cylinder, regular by hypothesis.
Applied to the fixed pair $(\beta_\epsilon,J_\epsilon)$,
the same argument shows that its index-one moduli spaces
are compact and, by regularity, finite.

Choose a common positive exponential weight smaller than
the spectral gaps of the finitely many relevant asymptotic
operators and small enough that this remains true under
perturbation.
Transport the basepoints and asymptotic coordinates along
the continued orbits.
On each end, use the asymptotic expansion of
\cite[Definition~2.6]{Pardon}, with the period and
parametrization of the corresponding continued orbit.
Express the difference from this asymptotic cylinder
in the transported coordinates.
After fixing domain parametrizations and target translations,
the exponential estimates give convergence of these differences
in a fixed exponentially weighted Sobolev space;
see \cite[Proposition~2.20 and Section~2.11]{Pardon}.
The implicit function theorem then gives a unique nearby
index-one cylinder for each original cylinder, with the
chosen asymptotic markers, and regularity persists.
Compactness shows that every nearby index-one cylinder
arises in this way.
There are finitely many original cylinders over the
finitely many relevant pairs of ends, so a single
perturbation neighborhood works for all of them.
Coherent orientations extend along the perturbation,
and the signed counts are unchanged.
\end{proof}

Proposition~\ref{prop:local-bw-counts} supplies the
regularity hypothesis of
Lemma~\ref{lem:local-bw-perturbation-persistence} for all
sufficiently small $\epsilon>0$.
Fix such an $\epsilon$, with $2N\epsilon<\tfrac14$,
and choose a pair $(\widetilde\beta,J)$ sufficiently
$C^\infty$-close to $(\beta_\epsilon,J_\epsilon)$,
with $\widetilde\beta$ globally nondegenerate, so that
the conclusions of the lemma hold
and \eqref{eq:finite-layer-action-bounds} holds with
$\delta=\tfrac14$.

For $1\le m\le\lfloor N/h\rfloor$, let $V_m$ be the
$\mathbb Q$-vector space generated by
$\gamma_{x,\widetilde\beta}^{mh}$, $x\in\operatorname{Crit}(f)$,
graded by \eqref{eq:contact-grading}.
Define $\partial_m:V_m\to V_m$ by the signed counts of
index-one $J$-holomorphic cylinders between these orbits,
with normalization factor $1/(mh)$.
We use the marking conventions of
Section~\ref{subsec:holomorphic-morse} and the coherent
orientations transported from $(\beta_\epsilon,J_\epsilon)$.
Let $d$ denote the differential of the contact algebra
$\mathcal B(\widetilde\beta)$ introduced in
Section~\ref{subsec:contact-algebras}.
Define the graded subalgebras
$A_{\le m}\subset\mathcal B(\widetilde\beta)$ by
\[
 A_{\le0}=\mathbb Q,\qquad
 A_{\le m}=\operatorname{Sym}(V_1\oplus\cdots\oplus V_m),
 \qquad A_{<m}=A_{\le m-1},
\]
where $\operatorname{Sym}$ denotes the free
graded-commutative algebra.

\begin{lemma}\label{lem:local-bw-homology}
For $1\le m\le\lfloor N/h\rfloor$, the algebra $A_{\le m}$ is a sub-DGA of
$\mathcal B(\widetilde\beta)$, and
there is a linear map $r_m:V_m\to A_{<m}$
of degree $-1$ such that
\begin{equation}\label{eq:reference-triangular-differential}
dv=\partial_m v+r_m(v)\qquad(v\in V_m).
\end{equation}
The differential $\partial_m$ is given by
\begin{equation}\label{eq:layer-morse-differential}
 \partial_m\gamma_{x,\widetilde\beta}^{mh}
 =\sum_{\operatorname{ind}_f(y)=\operatorname{ind}_f(x)-1}
       n_f(x,y)\,\gamma_{y,\widetilde\beta}^{mh}.
\end{equation}
In particular, $\partial_m^2=0$, and for every $q\in\mathbb Z$,
\begin{equation}\label{eq:reference-layer-dga-homology}
 H_q(V_m,\partial_m)
 \cong H_{q-2\sigma m c_B+2}(B;\mathbb Q).
\end{equation}
\end{lemma}

\begin{proof}
By Proposition~\ref{prop:local-bw-counts} and
Lemma~\ref{lem:local-bw-perturbation-persistence}, the
signed count of index-one cylinders from
$\gamma_{x,\widetilde\beta}^{mh}$ to
$\gamma_{y,\widetilde\beta}^{mh}$ is $mh\,n_f(x,y)$.
Division by $mh$ gives \eqref{eq:layer-morse-differential}.
Since the Conley--Zehnder indices are unchanged,
\eqref{eq:BWindex} and \eqref{eq:contact-grading} give
\[
 |\gamma_{x,\widetilde\beta}^{mh}|
 =2\sigma mc_B-2+\operatorname{ind}_f(x).
\]
Thus $(V_m,\partial_m)$ is the Morse complex of $f$ with
this grading shift. It follows that $\partial_m^2=0$ and
\eqref{eq:reference-layer-dga-homology} holds.

Fix a generator $v=\gamma_{x,\widetilde\beta}^{mh}\in V_m$.
The compactified moduli spaces defining $\partial_m$
consist of regular isolated cylinders.
The asymptotic markers eliminate finite domain automorphisms,
so the virtual counts agree with the ordinary signed counts
by \cite[Proposition~4.33]{Pardon}.
These cylinders therefore contribute $\partial_m v$ to $dv$.

Consider a monomial $\gamma_1^-\cdots\gamma_r^-$ with
nonzero coefficient in $dv$.
By Stokes' theorem and \eqref{eq:finite-layer-action-bounds},
every negative end has action less than $N+\tfrac12$.
All these orbits are contractible, since we work in
$\mathcal B(\widetilde\beta)$.
By Corollary~\ref{cor:finite-orbit-persistence}, each
$\gamma_i^-$ has covering multiplicity $m_i h$ for some
positive integer $m_i$.
Using \eqref{eq:finite-layer-action-bounds} again, we obtain
\[
 h\sum_{i=1}^r m_i
 \le \sum_{i=1}^r
       \mathcal A_{\widetilde\beta}(\gamma_i^-)
 \le \mathcal A_{\widetilde\beta}(v)
 < mh+\tfrac14.
\]
Since $h$ and the $m_i$ are positive integers,
$\sum_i m_i\le m$.
If some $m_i=m$, then $r=1$, so the corresponding
term comes from a cylinder whose two ends have
covering multiplicity $mh$.
Consequently, the terms of $dv$ involving generators
from $V_m$ sum to $\partial_m v$, and all remaining
monomials lie in $A_{<m}$.
This proves \eqref{eq:reference-triangular-differential}.
Since $d$ is a derivation, each $A_{\le m}$ is a sub-DGA.
\end{proof}

\section{Computation of linearized contact homology}
\label{sec:reference-tower}

Throughout this section, we assume that
$\omega|_{\pi_2(B)}\ne0$, $c_1(\xi)|_{\pi_2(M)}=0$, and
\begin{equation}\label{eq:betti-support}
 b_j(B)b_{j+2\sigma c_B-1}(B)=0
 \qquad (j\in\mathbb Z).
\end{equation}
Let $\alpha$ be a nondegenerate contact form supporting $\xi$.
We use the finitely generated sub-DGAs from
Section~\ref{sec:local-bw-proof} to construct augmentations
compatible with the cobordism maps as the action cutoff
increases.
These maps give an augmentation of $\mathcal B(\alpha)$,
whose pullback to $\mathcal C(\alpha)$ vanishes on
noncontractible generators.
We compute the homology of the direct limit of the resulting
linearized complexes and compare it with the contractible
linearized contact homology of $\alpha$ at this augmentation,
proving Theorem~\ref{thm:homology}.

\subsection{Extension and transport of augmentations}
\label{subsec:layer-augmentations}

\begin{lemma}\label{lem:betti-support}
The following assertions hold:
\begin{enumerate}
\item If $2c_B\le n$, then
$H_{\mathrm{odd}}(B;\mathbb Q)=0$.

\item For all $j\in\mathbb Z$ and $r\in\mathbb N$,
\[
 b_j(B) b_{j+2\sigma r c_B-1}(B)=0.
\]

\item For every $m\in\mathbb N$,
\[
 H_{3-2\sigma m c_B}(B;\mathbb Q)=0.
\]
\end{enumerate}
\end{lemma}

\begin{proof}
Every even Betti number $b_{2k}(B)$, $0\le k\le n$, is positive.
Indeed, $[\omega]^k\ne0$, since
\[
 \bigl\langle [\omega]^k\smile[\omega]^{n-k},[B]\bigr\rangle
 =\int_B\omega^n>0.
\]

For \textup{(i)}, suppose that $2c_B\le n$ and set
$s=|2\sigma c_B-1|$. Then $s$ is an odd positive integer
with $s\le n+1$.
Let $t$ be an odd integer with $1\le t\le2n-1$.

If $t\ge s$, then $t-s$ is even and lies in $[0,2n]$.
Thus $b_{t-s}(B)>0$, and
$b_{t-s}(B)b_t(B)=0$ by \eqref{eq:betti-support}.
It follows that $b_t(B)=0$.

If $t<s$, then $t\le s-2$, since both $t$ and $s$ are odd.
Hence
\[
 0<t+s\le2s-2\le2n.
\]
Thus $t+s$ is even and lies in $[0,2n]$, so $b_{t+s}(B)>0$.
Condition~\eqref{eq:betti-support} now gives
$b_t(B)b_{t+s}(B)=0$, and again $b_t(B)=0$.
This proves \textup{(i)}.

For \textup{(ii)}, if $2c_B\le n$, then \textup{(i)} applies, and the
assertion follows because the two degrees have opposite
parity.
Suppose instead that $2c_B>n$.
Since $c_B$ is an integer, $2c_B\ge n+1$.
Hence, for $r\ge2$,
\[
 |2\sigma r c_B-1|
 \ge4c_B-1
 \ge2n+1.
\]
The two degrees therefore cannot both lie in $[0,2n]$,
so at least one of the corresponding Betti numbers is zero.
For $r=1$, the assertion is exactly
\eqref{eq:betti-support}.

For \textup{(iii)}, apply \textup{(ii)} with
$r=m$ and $j=3-2\sigma m c_B$ to obtain
\[
 b_{3-2\sigma m c_B}(B)\,b_2(B)=0.
\]
Since $b_2(B)>0$, the assertion follows.
\end{proof}

Fix $N$, set $K=\lfloor N/h\rfloor$, and consider
the sub-DGAs $A_{\le m}$ of Lemma~\ref{lem:local-bw-homology},
for $0\le m\le K$.
With the grading in \eqref{eq:contact-grading},
\eqref{eq:reference-layer-dga-homology} and
Lemma~\ref{lem:betti-support}\textup{(iii)} give
\begin{equation}\label{eq:reference-layer-obstruction-vanishing}
 H_1(V_m,\partial_m)
 \cong H_{3-2\sigma m c_B}(B;\mathbb Q)=0
 \qquad (1\le m\le K).
\end{equation}

Recall that $r_m:V_m\to A_{<m}$ is defined by
\eqref{eq:reference-triangular-differential}.
Since $d^2=0$ by Section~\ref{subsec:contact-algebras}
and $\partial_m^2=0$ by Lemma~\ref{lem:local-bw-homology},
we have
\begin{equation}\label{eq:layer-differential-identities}
 r_m\partial_m+d r_m=0.
\end{equation}

\begin{lemma}\label{lem:layer-augmentation-extension}
Let $\varepsilon_{<m}:A_{<m}\to\mathbb Q$ be a graded augmentation.
If $H_1(V_m,\partial_m)=0$, then $\varepsilon_{<m}$ extends to a
graded augmentation of $A_{\le m}$.
\end{lemma}

\begin{proof}
Define the linear functional
\[
 \kappa_m=\varepsilon_{<m}\circ r_m:
 V_m\longrightarrow\mathbb Q.
\]
Since $r_m$ has degree $-1$ and $\varepsilon_{<m}$ is graded,
$\kappa_m|_{(V_m)_q}=0$ for $q\ne1$.
By \eqref{eq:layer-differential-identities} and
$\varepsilon_{<m}d=0$,
\[
 \kappa_m\partial_m
 =\varepsilon_{<m}r_m\partial_m
 =-\varepsilon_{<m}d r_m=0.
\]
The assumption $H_1(V_m,\partial_m)=0$ therefore gives
\[
 \ker\bigl(\partial_m:(V_m)_1\to(V_m)_0\bigr)
 =\partial_m((V_m)_2)
 \subset\ker\kappa_m.
\]
We may thus define a linear functional
$\ell_m:\partial_m((V_m)_1)\to\mathbb Q$ by
\begin{equation}\label{eq:layer-augmentation-equation}
 \ell_m(\partial_m v)=-\kappa_m(v)
 \qquad (v\in(V_m)_1).
\end{equation}

Extend $\ell_m$ linearly to $(V_m)_0$.
Let $\varepsilon:A_{\le m}\to\mathbb Q$ be the graded
algebra homomorphism determined by
\[
 \varepsilon|_{A_{<m}}=\varepsilon_{<m},\qquad
 \varepsilon|_{(V_m)_0}=\ell_m,\qquad
 \varepsilon|_{(V_m)_q}=0\quad(q\ne0).
\]
For $v\in(V_m)_1$, equations
\eqref{eq:reference-triangular-differential} and
\eqref{eq:layer-augmentation-equation} give
\[
 \varepsilon(dv)
 =\ell_m(\partial_m v)+\kappa_m(v)=0.
\]
For $v\in(V_m)_q$ with $q\ne1$,
$\varepsilon(dv)=0$ follows from grading.
Since $\varepsilon|_{A_{<m}}=\varepsilon_{<m}$,
we also have $\varepsilon d=0$ on $A_{<m}$.
The derivation rule then gives $\varepsilon d=0$
on $A_{\le m}$, so $\varepsilon$ is the required augmentation.
\end{proof}

By \eqref{eq:reference-layer-obstruction-vanishing} and
Lemma~\ref{lem:layer-augmentation-extension}, successive extensions
of $\mathrm{id}_{\mathbb Q}$ on $A_{\le0}=\mathbb Q$ give a
graded augmentation of $A_{\le K}$.

Let $\varepsilon$ be any graded augmentation of
$A_{\le K}$, and let
\[
 C=Q(A_{\le K},\varepsilon)
\]
be the linearized complex defined in
\eqref{eq:indecomposables}.
For $0\le m\le K$, define
\begin{equation}\label{eq:local-bw-filtration}
 \mathcal F_mC
 =\operatorname{span}_{\mathbb Q}
 \{[v-\varepsilon(v)]:v\in V_r,\ 1\le r\le m\}.
\end{equation}

\begin{lemma}\label{lem:local-bw-linearization}
The subspaces $\mathcal F_mC$ form an increasing filtration of $C$
by subcomplexes, with
\[
 \mathcal F_0C=0,\qquad
 \mathcal F_{K}C=C.
\]
For $1\le m\le K$, the map
\begin{equation}\label{eq:local-bw-filtration-quotient}
 V_m\longrightarrow\mathcal F_mC/\mathcal F_{m-1}C,
 \qquad v\longmapsto
 [v-\varepsilon(v)]+\mathcal F_{m-1}C
\end{equation}
is an isomorphism of chain complexes.
Thus the homology of the $m$-th successive quotient
$\mathcal F_mC/\mathcal F_{m-1}C$ is given by
\eqref{eq:reference-layer-dga-homology}.
\end{lemma}

\begin{proof}
For $v\in V_m$, equation
\eqref{eq:reference-triangular-differential} and $\varepsilon d=0$ give
\[
 [dv]
 =[\partial_m v-\varepsilon(\partial_m v)]
 +[r_m(v)-\varepsilon(r_m(v))].
\]
Since $r_m(v)\in A_{<m}$, the description of the linearized
complex in Section~\ref{subsec:contact-algebras} gives
\[
 [r_m(v)-\varepsilon(r_m(v))]\in\mathcal F_{m-1}C.
\]
The linearized differential sends $[v-\varepsilon(v)]$ to $[dv]$
and therefore preserves each $\mathcal F_mC$.

Since $A_{\le K}$ is the free graded-commutative algebra on
$\bigoplus_{r=1}^K V_r$, the map in
\eqref{eq:local-bw-filtration-quotient} is an isomorphism
of graded vector spaces. Under this identification, the
induced differential on the quotient is $\partial_m$.
\end{proof}

\Needspace{8\baselineskip}
Let
\[
 A'_{\le K}=\operatorname{Sym}(V'_1\oplus\cdots\oplus V'_{K})
\]
be a graded DGA, and let $A'_{<m}$ be the subalgebra
generated by $V'_1,\ldots,V'_{m-1}$, with $A'_{<1}=\mathbb Q$.
Consider a graded DGA map
$F:A_{\le K}\to A'_{\le K}$ such that, for $1\le m\le K$,
\begin{equation}\label{eq:layer-triangular-map}
 F(v)=T_m(v)+P_m(v)\qquad(v\in V_m),
\end{equation}
where
\[
 T_m:V_m\longrightarrow V'_m,
 \qquad
 P_m:V_m\longrightarrow A'_{<m}
\]
are grading-preserving linear maps.

For a graded augmentation $\varepsilon':A'_{\le K}\to\mathbb Q$
satisfying $\varepsilon'\circ F=\varepsilon$, set
\[
 C'=Q(A'_{\le K},\varepsilon'),
\]
with filtration defined as in \eqref{eq:local-bw-filtration}.
As in Section~\ref{subsec:contact-algebras}, let
$QF:C\to C'$ denote the linearization of $F$
at $\varepsilon$ and $\varepsilon'$.

\begin{lemma}\label{lem:layer-augmentation-transport}
Let $F:A_{\le K}\to A'_{\le K}$ be a graded DGA map of the form
\eqref{eq:layer-triangular-map}, and suppose that every $T_m$
is invertible. Then $F$ is a DGA isomorphism.
For every graded augmentation
$\varepsilon:A_{\le K}\to\mathbb Q$, the map
$\varepsilon'=\varepsilon\circ F^{-1}$ is the unique
graded augmentation of $A'_{\le K}$ satisfying
$\varepsilon'\circ F=\varepsilon$.
Moreover, $QF:C\to C'$ is a chain isomorphism satisfying
$QF(\mathcal F_mC)=\mathcal F_mC'$ for $0\le m\le K$.
Under the identifications \eqref{eq:local-bw-filtration-quotient},
$QF$ induces $T_m:V_m\to V'_m$ on the $m$-th quotients
for $1\le m\le K$.
\end{lemma}

\begin{proof}
We construct an inverse to $F$ by induction on $m$.
For $0\le m\le K$, write
\[
 F_m=F|_{A_{\le m}}:A_{\le m}\longrightarrow A'_{\le m}.
\]
Set $G_0=\mathrm{id}_{\mathbb Q}$, so that
$F_0G_0=G_0F_0=\mathrm{id}_{\mathbb Q}$.

For $1\le m\le K$, suppose that a graded algebra map
$G_{m-1}:A'_{\le m-1}\to A_{\le m-1}$ satisfies
\begin{equation}\label{eq:layer-inverse-lower-relations}
 F_{m-1}G_{m-1}=\mathrm{id}_{A'_{\le m-1}},
 \qquad
 G_{m-1}F_{m-1}=\mathrm{id}_{A_{\le m-1}}.
\end{equation}
Extend $G_{m-1}$ to a graded algebra map
$G_m:A'_{\le m}\to A_{\le m}$ by setting
\begin{equation}\label{eq:layer-triangular-inverse}
 G_m(w)=T_m^{-1}(w)
 -G_{m-1}\bigl(P_m(T_m^{-1}(w))\bigr)
 \qquad(w\in V'_m).
\end{equation}
For $v\in V_m$, equations
\eqref{eq:layer-triangular-map},
\eqref{eq:layer-inverse-lower-relations}, and
\eqref{eq:layer-triangular-inverse} give
\[
\begin{aligned}
 G_m(F_m(v))
 &=G_m(T_m(v))+G_{m-1}(P_m(v))\\
 &=v-G_{m-1}(P_m(v))+G_{m-1}(P_m(v))
 =v,
\end{aligned}
\]
and
\[
\begin{aligned}
 F_m(G_m(T_m(v)))
 &=F_m(v)-F_{m-1}(G_{m-1}(P_m(v)))\\
 &=T_m(v)+P_m(v)-P_m(v)
 =T_m(v).
\end{aligned}
\]
Since $T_m(V_m)=V'_m$, we obtain
\[
 F_mG_m=\mathrm{id}_{A'_{\le m}},
 \qquad
 G_mF_m=\mathrm{id}_{A_{\le m}}.
\]
Taking $m=K$ gives $G_{K}=F^{-1}$.

Let $d$ and $d'$ denote the differentials on $A_{\le K}$
and $A'_{\le K}$, respectively.
The relation $d'F=Fd$ implies
$dF^{-1}=F^{-1}d'$, so $F^{-1}$ is a DGA map.
For a graded augmentation $\varepsilon:A_{\le K}\to\mathbb Q$,
the map $\varepsilon'=\varepsilon\circ F^{-1}$ is a
graded unital algebra homomorphism satisfying
\[
 \varepsilon'd'
 =\varepsilon F^{-1}d'
 =\varepsilon dF^{-1}=0,
 \qquad
 \varepsilon'F=\varepsilon.
\]
Thus $\varepsilon'$ is a graded augmentation.
Since $F$ is invertible, the relation
$\varepsilon'\circ F=\varepsilon$ determines $\varepsilon'$ uniquely.

The DGA isomorphism $F$ with $\varepsilon'F=\varepsilon$ induces
a chain isomorphism $QF:C\to C'$.
For $0\le m\le K$, we have
\[
 F(A_{\le m})=A'_{\le m},
 \qquad F^{-1}|_{A'_{\le m}}=G_m.
\]
Thus $A'_{\le m}$ is a sub-DGA and
$QF(\mathcal F_mC)=\mathcal F_mC'$.
For $v\in V_m$, the relations
$\varepsilon'F=\varepsilon$ and $P_m(v)\in A'_{<m}$,
together with \eqref{eq:layer-triangular-map}, give
\begin{equation}\label{eq:layer-linearized-map}
 QF\bigl([v-\varepsilon(v)]\bigr)
 \equiv [T_m(v)-\varepsilon'(T_m(v))]
 \pmod{\mathcal F_{m-1}C'}.
\end{equation}
Under the identifications \eqref{eq:local-bw-filtration-quotient},
\eqref{eq:layer-linearized-map} shows that $QF$ induces $T_m$
on the $m$-th quotients.
\end{proof}

\subsection{Perturbations and cobordism maps}
\label{subsec:reference-construction}

We use the Morse function $f$ and the almost complex
structure $J_B$ fixed in
Section~\ref{subsec:finite-action-perturbations}.

For each positive integer $N$, decrease $\epsilon_0(N)$
from Lemma~\ref{lem:local-bw-finite}, if necessary, so that
Proposition~\ref{prop:local-bw-counts} applies for every
$0<\epsilon<\epsilon_0(N)$ and
\[
 2N\epsilon_0(N)<\tfrac14.
\]
For two pairs $(\widetilde\beta,J)$ and
$(\widetilde\beta',J')$ chosen as in
Section~\ref{subsec:layer-complexes} for the same $N$,
write
\[
 V_m,\ A_{\le m}
 \qquad\text{and}\qquad
 V'_m,\ A'_{\le m}
\]
for the corresponding graded vector spaces and subalgebras.
By Corollary~\ref{cor:finite-orbit-persistence},
\[
 A_{\le\lfloor N/h\rfloor}
 =\mathcal B(\widetilde\beta)^{<N+\frac12},
 \qquad
 A'_{\le\lfloor N/h\rfloor}
 =\mathcal B(\widetilde\beta')^{<N+\frac12}.
\]
Here $\mathcal B(\alpha)^{<L}$ denotes the unital subalgebra
generated by the contractible good orbits of action $<L$,
as in Section~\ref{subsec:contact-algebras}.

For a finite sequence of globally nondegenerate contact
forms supporting $\xi$,
\[
 \eta_0=\widetilde\beta>\eta_1>\cdots>\eta_s
 =\widetilde\beta',
\]
fix the DGAs $\mathcal B(\eta_i)$ as in
Section~\ref{subsec:contact-algebras}.
Lemma~\ref{lem:found-action} gives cobordism maps
\[
 G_i:\mathcal B(\eta_{i-1})
 \longrightarrow\mathcal B(\eta_i),
 \qquad 1\le i\le s,
\]
satisfying
\[
 G_i\bigl(\mathcal B(\eta_{i-1})^{<N+\frac12}\bigr)
 \subset \mathcal B(\eta_i)^{<N+\frac12}.
\]
For each such choice, define the graded DGA map
\begin{equation}\label{eq:reference-cobordism-map}
 F=\left.(G_s\circ\cdots\circ G_1)
 \right|_{A_{\le\lfloor N/h\rfloor}}:
 A_{\le\lfloor N/h\rfloor}
 \longrightarrow A'_{\le\lfloor N/h\rfloor}.
\end{equation}

\begin{lemma}\label{lem:close-reference-triangular}
Fix a positive integer $N$ and
$0<\epsilon'<\epsilon<\epsilon_0(N)$.
For any $C^\infty$-neighborhoods $\mathcal U$ and
$\mathcal U'$ of $(\beta_\epsilon,J_\epsilon)$ and
$(\beta_{\epsilon'},J_{\epsilon'})$, respectively,
there exist pairs
\[
 (\widetilde\beta,J)\in\mathcal U,
 \qquad
 (\widetilde\beta',J')\in\mathcal U',
\]
a sequence $(\eta_i)_{i=0}^s$, and maps $(G_i)_{i=1}^s$
such that $F$ in \eqref{eq:reference-cobordism-map}
has the form \eqref{eq:layer-triangular-map}
and, for every $1\le m\le\lfloor N/h\rfloor$,
the linear map $T_m:V_m\to V'_m$ in
\eqref{eq:layer-triangular-map} satisfies
\begin{equation}\label{eq:reference-cylinder-coefficients}
 T_m(\gamma_{x,\widetilde\beta}^{mh})
 =
 \gamma_{x,\widetilde\beta'}^{mh}
 +
 \sum_{\substack{y\in\operatorname{Crit}(f)\\f(y)<f(x)}}
 a_{yx}^{(m)}\gamma_{y,\widetilde\beta'}^{mh},
 \qquad a_{yx}^{(m)}\in\mathbb Q,
\end{equation}
for every $x\in\operatorname{Crit}(f)$.
In particular, by
Lemma~\ref{lem:layer-augmentation-transport},
$F$ is a graded DGA isomorphism.
\end{lemma}

\begin{proof}
For every $\tau\in[\epsilon',\epsilon]$, the orbits of
$\beta_\tau$ below $N+\tfrac12$ are those in
Lemma~\ref{lem:local-bw-finite}, and their actions
satisfy \eqref{eq:finite-layer-action-bounds} with
$\delta=\tfrac14$.
By \eqref{eq:reference-orbit-action},
\begin{equation}\label{eq:reference-uniform-action-gap}
\begin{aligned}
 \Delta
 &:=h\epsilon'
   \min_{\substack{x,y\in\operatorname{Crit}(f)\\x\ne y}}
   |f(x)-f(y)|>0,\\
 \left|
 \mathcal A_{\beta_\tau}(\gamma_{x,\beta_\tau}^{mh})
 -
 \mathcal A_{\beta_\tau}(\gamma_{y,\beta_\tau}^{mh})
 \right|
 &=mh\tau|f(x)-f(y)|\ge\Delta
\end{aligned}
\end{equation}
for $\tau\in[\epsilon',\epsilon]$,
$1\le m\le\lfloor N/h\rfloor$, and
distinct $x,y\in\operatorname{Crit}(f)$.

We first compare perturbation pairs
$(\widetilde\beta,J)$ and $(\widetilde\beta',J')$
close to $(\beta_{\tau^+},J_{\tau^+})$ and
$(\beta_{\tau^-},J_{\tau^-})$, respectively, where
\[
 \epsilon'\le\tau^-<\tau^+\le\epsilon
\]
and $\tau^+-\tau^-$ is sufficiently small.
Choose globally nondegenerate forms satisfying
$\widetilde\beta>\widetilde\beta'$.
Choose the perturbations within the neighborhoods supplied by
Corollary~\ref{cor:finite-orbit-persistence} and
Lemma~\ref{lem:local-bw-perturbation-persistence}, so that
the orbit correspondence, action bounds, and cylinder
differentials are preserved.

Write $\widetilde\beta=h_+\beta$ and
$\widetilde\beta'=h_-\beta$.
Fix a smooth nondecreasing function
$\chi:\mathbb R\to[0,1]$ with $\chi(s)=0$ for $s\le-1$
and $\chi(s)=1$ for $s\ge1$, and set
\[
 h_s=(1-\chi(s))h_-+\chi(s)h_+,\qquad
 \Lambda=e^s h_s\beta.
\]
Choose an almost complex structure $\widehat J$ as in
Lemma~\ref{lem:found-action}, equal to $J$ and $J'$ on
the positive and negative ends, respectively.
These choices can be made so that
\begin{equation}\label{eq:reference-cobordism-convergence}
 (\Lambda,\widehat J)
 \longrightarrow(e^s\beta_\tau,J_\tau)
 \quad\text{in }C^\infty_{\mathrm{loc}}(\mathbb R\times M)
\end{equation}
whenever both perturbation pairs converge to
$(\beta_\tau,J_\tau)$.
For sufficiently close pairs, $d\Lambda$ uniformly tames
$\widehat J$ on compact subsets.

Let
\[
 G:\mathcal B(\widetilde\beta)
 \longrightarrow\mathcal B(\widetilde\beta')
\]
be the DGA map induced by the exact symplectic cobordism
$(\mathbb R\times M,d\Lambda)$, with the almost complex
structure $\widehat J$ and auxiliary data chosen as in
Lemma~\ref{lem:found-action}.
By \eqref{eq:found-total-output-action}, every curve
contributing to $G(\gamma)$ satisfies
\begin{equation}\label{eq:continuation-action}
 \sum_i\mathcal A_{\widetilde\beta'}(\delta_i)
 \le\mathcal A_{\widetilde\beta}(\gamma),
\end{equation}
where the sum is zero for a constant term.

Suppose that $\gamma\in V_m$ is a generator
and that $\delta_1\cdots\delta_k$ has nonzero
coefficient in $G(\gamma)$, with $k=0$ allowed.
Every $\delta_i$ is contractible and has action below
$N+\tfrac12$, so its covering multiplicity is
$m_i h$ for some positive integer $m_i$.
Equations \eqref{eq:finite-layer-action-bounds} and
\eqref{eq:continuation-action} give
\[
 h\sum_{i=1}^k m_i
 \le\sum_{i=1}^k
      \mathcal A_{\widetilde\beta'}(\delta_i)
 \le\mathcal A_{\widetilde\beta}(\gamma)
 <mh+\tfrac14.
\]
Hence $\sum_i m_i\le m$.
If some $m_i=m$, then $k=1$.
Thus the restriction of $G$ to
$A_{\le\lfloor N/h\rfloor}$ has the form
\eqref{eq:layer-triangular-map}:
\begin{equation}\label{eq:reference-local-map}
\begin{gathered}
 G(v)=S_m(v)+P_m(v)\qquad(v\in V_m),\\
 S_m:V_m\longrightarrow V'_m,\qquad
 P_m:V_m\longrightarrow A'_{<m},
\end{gathered}
\end{equation}
where $S_m$ and $P_m$ are grading-preserving linear
maps, and $S_m$ is given by the cylinder coefficients.

Write
\[
 S_m(\gamma_{x,\widetilde\beta}^{mh})
 =
 \sum_{y\in\operatorname{Crit}(f)}
 b_{yx}^{(m)}\gamma_{y,\widetilde\beta'}^{mh},
 \qquad b_{yx}^{(m)}\in\mathbb Q.
\]
By \eqref{eq:reference-uniform-action-gap}, we can
choose the perturbations so that
\[
 \mathcal A_{\widetilde\beta}
   (\gamma_{y,\widetilde\beta}^{mh})
 -
 \mathcal A_{\widetilde\beta}
   (\gamma_{x,\widetilde\beta}^{mh})
 \ge\frac{\Delta}{2}
 \qquad\text{if }f(y)>f(x).
\]
Taking $\tau^+-\tau^-$ and the perturbations smaller
if necessary, we also have
\[
 \left|
 \mathcal A_{\widetilde\beta'}
   (\gamma_{y,\widetilde\beta'}^{mh})
 -
 \mathcal A_{\widetilde\beta}
   (\gamma_{y,\widetilde\beta}^{mh})
 \right|
 <\frac{\Delta}{4}
\]
for every $y\in\operatorname{Crit}(f)$ and
$1\le m\le\lfloor N/h\rfloor$.
Consequently, if $f(y)>f(x)$, then
\[
 \mathcal A_{\widetilde\beta'}
   (\gamma_{y,\widetilde\beta'}^{mh})
 >
 \mathcal A_{\widetilde\beta}
   (\gamma_{x,\widetilde\beta}^{mh}),
\]
and \eqref{eq:continuation-action} gives
$b_{yx}^{(m)}=0$.

We next prove that $b_{xx}^{(m)}=1$ when
$\tau^+-\tau^-$ and the perturbations are sufficiently
small, uniformly for $\tau^\pm\in[\epsilon',\epsilon]$.
By Lemma~\ref{lem:local-bw-finite} and compactness of
$[\epsilon',\epsilon]$, the asymptotic operators of
the orbits below $N+\tfrac12$ have a uniform positive
spectral gap.
Choose a common exponential weight smaller than
this gap.

In the Fredholm setup of \cite[Section~2.11]{Pardon},
the linearized operators at the trivial cylinders are
invertible by \cite[Lemma~2.40]{Pardon}.
The parameterized implicit function theorem gives
a unique nearby regular cylinder for sufficiently
small changes in the two perturbation pairs and
in $(\Lambda,\widehat J)$.
The inverses are locally bounded in $\tau$, and a
finite cover of $[\epsilon',\epsilon]$ gives a
uniform bound.

Suppose that, for arbitrarily small
$\tau^+-\tau^-$ and perturbations, there is a
cylinder from $\gamma_{x,\widetilde\beta}^{mh}$ to
$\gamma_{x,\widetilde\beta'}^{mh}$ distinct from the
cylinder given by the implicit function theorem.
Choose a sequence of such cylinders $u_\nu$.
After passing to a subsequence, we may assume that
$m$ and $x$ are fixed and that
\[
 \tau_\nu^\pm\longrightarrow\tau_*
 \in[\epsilon',\epsilon],
 \qquad
 (\widetilde\beta_\nu^\pm,J_\nu^\pm)
 \longrightarrow(\beta_{\tau_*},J_{\tau_*}).
\]
By \eqref{eq:reference-cobordism-convergence},
\[
 (\Lambda_\nu,\widehat J_\nu)
 \longrightarrow(e^s\beta_{\tau_*},J_{\tau_*})
 \quad\text{in }C^\infty_{\mathrm{loc}}(\mathbb R\times M).
\]
Moreover,
\[
 \mathcal A_{\widetilde\beta_\nu^+}
   (\gamma_{x,\widetilde\beta_\nu^+}^{mh})
 -
 \mathcal A_{\widetilde\beta_\nu^-}
   (\gamma_{x,\widetilde\beta_\nu^-}^{mh})
 \longrightarrow0.
\]

The positive asymptotic actions give a uniform
Hofer-energy bound.
By \eqref{eq:continuation-action}, every intermediate
orbit lies below $N+\tfrac12$.
Compactness therefore gives a limiting holomorphic
building with zero $d\beta_{\tau_*}$-area.
Zero $d\beta_{\tau_*}$-area excludes planes, and exactness
excludes nonconstant sphere bubbles.
A genus-zero building with one positive and one
negative external end cannot branch without a plane.
Every remaining component is therefore a trivial
cylinder over $\gamma_{x,\beta_{\tau_*}}^{mh}$.

After removing additional trivial symplectization levels,
choose domain parametrizations compatible with the fixed
asymptotic markers.
Identify the asymptotic coordinates along the continued
orbits, as in the proof of
Lemma~\ref{lem:local-bw-perturbation-persistence}.
With each end expressed relative to its corresponding
asymptotic cylinder, the uniform exponential estimates and
convergence on compact subsets give convergence to the
trivial cobordism cylinder in the weighted Sobolev topology
for the common weight chosen above.
For sufficiently large $\nu$, the cylinders $u_\nu$ lie in the
neighborhood where the implicit function theorem
gives the unique local continuation, a contradiction.

Transport coherent orientations along the matched
orbits.
With the counting conventions of
Section~\ref{subsec:contact-algebras}, the local cylinder
contributes $+1$ to $b_{xx}^{(m)}$, as in the proof of
\cite[Lemma~1.2]{Pardon}.
The agreement of ordinary and virtual counts follows
from \cite[Proposition~4.33]{Pardon}.
We have therefore proved
\begin{equation}\label{eq:reference-local-coefficients}
 b_{xx}^{(m)}=1,\qquad
 b_{yx}^{(m)}=0\quad\text{if }f(y)>f(x).
\end{equation}
No regularity of the off-diagonal moduli spaces is needed.

Choose a finite subdivision
\[
 \epsilon=\tau_0>\tau_1>\cdots>\tau_s=\epsilon'
\]
such that the action comparison and the diagonal
calculation apply to each consecutive pair of parameters.
At each $\tau_i$, choose a globally nondegenerate
perturbation $\eta_i$ of $\beta_{\tau_i}$ and a
cylindrical almost complex structure $J_i$ so that
the conclusions of Lemma~\ref{lem:local-bw-homology} hold.
The pairs $(\eta_i,J_i)$ can be chosen arbitrarily
$C^\infty$-close to $(\beta_{\tau_i},J_{\tau_i})$.
Since
$\beta_{\tau_0}>\beta_{\tau_1}>\cdots>\beta_{\tau_s}$,
we can choose the perturbations small enough that
\[
 \eta_0>\eta_1>\cdots>\eta_s,\qquad
 (\eta_0,J_0)\in\mathcal U,\qquad
 (\eta_s,J_s)\in\mathcal U',
\]
and that each consecutive pair satisfies the
conditions used to prove
\eqref{eq:reference-local-coefficients}.
Set
\[
 (\widetilde\beta,J)=(\eta_0,J_0),\qquad
 (\widetilde\beta',J')=(\eta_s,J_s).
\]

Write $V_m^{(i)}$ and $A_{\le m}^{(i)}$ for the
graded vector spaces and subalgebras corresponding
to $(\eta_i,J_i)$.
At $i=0$ and $i=s$, these are
$V_m,A_{\le m}$ and $V'_m,A'_{\le m}$, respectively.
Fix the DGAs $\mathcal B(\eta_i)$ as in
Section~\ref{subsec:contact-algebras}, and choose
cobordism maps
\[
 G_i:\mathcal B(\eta_{i-1})
 \longrightarrow\mathcal B(\eta_i),
 \qquad 1\le i\le s,
\]
induced by the exact symplectic cobordisms constructed
in Lemma~\ref{lem:found-action}, with each $G_i$ satisfying
\eqref{eq:reference-local-map} and
\eqref{eq:reference-local-coefficients}.
Set
\[
 F_i=
 \left.G_i\right|_{A_{\le\lfloor N/h\rfloor}^{(i-1)}}:
 A_{\le\lfloor N/h\rfloor}^{(i-1)}
 \longrightarrow A_{\le\lfloor N/h\rfloor}^{(i)},
\]
and let
\[
 S_{m,i}:V_m^{(i-1)}\longrightarrow V_m^{(i)}
\]
denote the linear map $S_m$ in
\eqref{eq:reference-local-map} for $G_i$.
Since
\[
 F_i\bigl(A_{<m}^{(i-1)}\bigr)\subset A_{<m}^{(i)},
\]
the map $F$ defined in
\eqref{eq:reference-cobordism-map} satisfies
\begin{equation}\label{eq:reference-step-composition}
 F=F_s\circ\cdots\circ F_1,\qquad
 T_m=S_{m,s}\circ\cdots\circ S_{m,1}.
\end{equation}
By \eqref{eq:reference-local-coefficients}, each
$S_{m,i}$ is triangular with diagonal entries one
when the generators are ordered by the
critical values of $f$.
Equation \eqref{eq:reference-step-composition}
therefore gives
\eqref{eq:reference-cylinder-coefficients}.
The composite $F$ also has the form
\eqref{eq:layer-triangular-map}.
Since each $T_m$ is invertible,
Lemma~\ref{lem:layer-augmentation-transport} implies
that $F$ is a graded DGA isomorphism.
\end{proof}

For the construction of the reference sequence and the comparison
below, assume that $\alpha>2\beta$.
We remove this assumption by constant rescaling in the proof of
Theorem~\ref{thm:homology}.

\begin{proposition}\label{prop:reference-sequence}
There exist increasing positive integers $N_j\to\infty$
and globally nondegenerate contact forms $\beta_j$,
chosen as in Section~\ref{subsec:layer-complexes}
with $N=N_j$, such that
\begin{equation}\label{eq:reference-sequence}
 \alpha>\beta_1>\beta_2>\cdots>\beta.
\end{equation}
The inequalities are understood in the sense of
\eqref{eq:contact-form-order}.
Set
\begin{equation}\label{eq:reference-algebras}
 L_j:=N_j+\tfrac12\qquad(j\ge1).
\end{equation}
For each $j$, write $V_{m,j}$ and $\partial_{m,j}$
for the spaces $V_m$ and differentials $\partial_m$
defined in Section~\ref{subsec:layer-complexes}
for $\beta_j$.
The following assertions hold.

\begin{enumerate}
\item[\textup{(i)}]
For every $j\ge1$, we have
\[
 \mathcal B(\beta_j)^{<L_j}=\operatorname{Sym}\left(
       \bigoplus_{1\le m\le\lfloor N_j/h\rfloor}V_{m,j}
      \right).
\]
Its differential has the form
\eqref{eq:reference-triangular-differential}, and
$H_*(V_{m,j},\partial_{m,j})$ is given by
\eqref{eq:reference-layer-dga-homology}.

\item[\textup{(ii)}]
For every $j\ge1$, there is an integer-graded DGA
quasi-isomorphism
$F_j:\mathcal B(\beta_j)\to\mathcal B(\beta_{j+1})$ satisfying
\begin{equation}\label{eq:reference-transition-action}
 F_j\bigl(\mathcal B(\beta_j)^{<L}\bigr)
 \subset\mathcal B(\beta_{j+1})^{<L}
 \qquad\text{for every }L>0.
\end{equation}
Its restriction
\[
 f_j:=\left.F_j\right|_{\mathcal B(\beta_j)^{<L_j}}:
 \mathcal B(\beta_j)^{<L_j}\longrightarrow\mathcal B(\beta_{j+1})^{<L_{j+1}}
\]
is a DGA isomorphism onto $\mathcal B(\beta_{j+1})^{<L_j}$, and
\begin{equation}\label{eq:reference-transition-image}
 f_j(\mathcal B(\beta_j)^{<L_j})
 =\mathcal B(\beta_{j+1})^{<L_j}
 =\operatorname{Sym}\left(
       \bigoplus_{1\le m\le\lfloor N_j/h\rfloor}V_{m,j+1}
      \right).
\end{equation}
As a map $\mathcal B(\beta_j)^{<L_j}\to\mathcal B(\beta_{j+1})^{<L_j}$,
$f_j$ has the form \eqref{eq:layer-triangular-map}, with
invertible linear maps $T_{m,j}:V_{m,j}\to V_{m,j+1}$
in place of $T_m$ for $1\le m\le\lfloor N_j/h\rfloor$.
\end{enumerate}
\end{proposition}

\begin{proof}
Choose increasing positive integers $N_j\to\infty$
and positive numbers $\epsilon_j$ such that
\begin{equation}\label{eq:reference-parameters}
 \epsilon_j<\epsilon_0(N_j),\qquad
 \epsilon_{j+1}<\epsilon_j/4,\qquad
 2N_j\epsilon_j<1/8.
\end{equation}
For each $j$, choose a finite subdivision of
$[\epsilon_{j+1},\epsilon_j]$ as in the proof of
Lemma~\ref{lem:close-reference-triangular}, with $N=N_j$.

At each subdivision point $\tau$, choose one globally
nondegenerate perturbation of $\beta_\tau$ and a
cylindrical almost complex structure close to $J_\tau$,
so that the construction in the proof of
Lemma~\ref{lem:close-reference-triangular} applies to
every consecutive pair.
At $\epsilon_j$, denote the chosen form by $\beta_j$
and also impose the conditions of
Section~\ref{subsec:layer-complexes} with $N=N_j$.
For $j\ge2$, use the same pair at $\epsilon_j$ in both
adjacent subdivisions.
Only finitely many conditions occur at each subdivision
point, so they can be satisfied by taking sufficiently
small perturbations.
The strict pointwise inequalities between consecutive
forms $\beta_\tau$ allow the perturbed forms to be
chosen in decreasing order and above $\beta$.
Since $\beta_{\epsilon_1}<2\beta<\alpha$, we can also
choose $\beta_1<\alpha$.
Thus \eqref{eq:reference-sequence} holds.

By Corollary~\ref{cor:finite-orbit-persistence}, the
orbits of $\beta_j$ below $L_j$ are precisely the
continuations $\gamma_{x,\beta_j}^{\ell}$ of
$\gamma_{x,\beta_{\epsilon_j}}^{\ell}$, with
$x\in\operatorname{Crit}(f)$ and $1\le\ell\le N_j$.
The bounds \eqref{eq:finite-layer-action-bounds},
with $\delta=\tfrac14$, give
\begin{equation}\label{eq:reference-action-separation}
 \ell<\mathcal A_{\beta_j}(\gamma_{x,\beta_j}^{\ell})
 <\ell+\tfrac14
 \qquad
 (x\in\operatorname{Crit}(f),\ 1\le\ell\le N_j).
\end{equation}
These orbits are contractible exactly when $h\mid\ell$.
Lemma~\ref{lem:local-bw-homology} therefore gives
\textup{(i)}.

Fix the DGAs at all chosen contact forms as in
Section~\ref{subsec:contact-algebras}.
For each consecutive pair, choose a cobordism map as in
the proof of Lemma~\ref{lem:close-reference-triangular},
satisfying \eqref{eq:reference-local-map} and
\eqref{eq:reference-local-coefficients}.
Let $F_j:\mathcal B(\beta_j)\to\mathcal B(\beta_{j+1})$ be their
composition along the subdivision from $\beta_j$ to
$\beta_{j+1}$.
Each factor is an integer-graded quasi-isomorphism by
Lemma~\ref{lem:found-contractible-quotient}, so the same
holds for $F_j$.
Lemma~\ref{lem:found-action} gives
\eqref{eq:reference-transition-action}.
Thus $f_j$ takes values in $\mathcal B(\beta_{j+1})^{<L_j}$.
By \eqref{eq:reference-action-separation},
$\mathcal B(\beta_{j+1})^{<L_j}$ is generated by the spaces
$V_{m,j+1}$ with $mh\le N_j$.
Lemma~\ref{lem:close-reference-triangular} shows that
\[
 f_j:\mathcal B(\beta_j)^{<L_j}\longrightarrow\mathcal B(\beta_{j+1})^{<L_j}
\]
has the form \eqref{eq:layer-triangular-map}, with
each $T_{m,j}:V_{m,j}\to V_{m,j+1}$ invertible.
Lemma~\ref{lem:layer-augmentation-transport} then
implies that this map is a DGA isomorphism.
This proves \textup{(ii)}.
\end{proof}

\begin{corollary}\label{cor:reference-sequence-augmentations}
The algebras $\mathcal B(\beta_j)^{<L_j}$ in
Proposition~\ref{prop:reference-sequence} admit graded
augmentations $\bar{\varepsilon}_{\beta_j}$ satisfying
\begin{equation}\label{eq:augmentation-compatibility}
 \bar{\varepsilon}_{\beta_{j+1}}\circ f_j=\bar{\varepsilon}_{\beta_j}
 \qquad(j\ge1).
\end{equation}
The maps
\[
 Qf_j:Q(\mathcal B(\beta_j)^{<L_j},\bar{\varepsilon}_{\beta_j})
 \longrightarrow Q(\mathcal B(\beta_{j+1})^{<L_{j+1}},\bar{\varepsilon}_{\beta_{j+1}})
\]
preserve the filtrations defined by
\eqref{eq:local-bw-filtration}.
Under the identifications
\eqref{eq:local-bw-filtration-quotient}, $Qf_j$ induces
the chain isomorphism
\[
 T_{m,j}:(V_{m,j},\partial_{m,j})
 \longrightarrow(V_{m,j+1},\partial_{m,j+1})
\]
on the $m$-th quotient, for $1\le m\le\lfloor N_j/h\rfloor$.
\end{corollary}

\begin{proof}
We construct $\bar{\varepsilon}_{\beta_j}$ by induction on $j$.
By \eqref{eq:reference-layer-obstruction-vanishing}
and repeated application of
Lemma~\ref{lem:layer-augmentation-extension},
the augmentation $\mathrm{id}_{\mathbb Q}$ extends to
a graded augmentation $\bar{\varepsilon}_{\beta_1}:\mathcal B(\beta_1)^{<L_1}\to\mathbb Q$.

Suppose that $\bar{\varepsilon}_{\beta_j}$ has been chosen.
By Proposition~\ref{prop:reference-sequence}\textup{(ii)},
$f_j:\mathcal B(\beta_j)^{<L_j}\to\mathcal B(\beta_{j+1})^{<L_j}$
is a DGA isomorphism, so
\[
 \bar{\varepsilon}_{\beta_j}\circ f_j^{-1}:
 \mathcal B(\beta_{j+1})^{<L_j}\longrightarrow\mathbb Q
\]
is a graded augmentation.
Apply Lemma~\ref{lem:layer-augmentation-extension}
successively in increasing order of $m$,
for $N_j<mh\le N_{j+1}$, using
\eqref{eq:reference-layer-obstruction-vanishing}
at each step.
We obtain a graded augmentation
$\bar{\varepsilon}_{\beta_{j+1}}:\mathcal B(\beta_{j+1})^{<L_{j+1}}\to\mathbb Q$ satisfying
\[
 \left.\bar{\varepsilon}_{\beta_{j+1}}\right|_{\mathcal B(\beta_{j+1})^{<L_j}}
 =\bar{\varepsilon}_{\beta_j}\circ f_j^{-1}.
\]
Hence \eqref{eq:augmentation-compatibility} holds.

By Lemma~\ref{lem:layer-augmentation-transport},
the linearization of
$f_j:\mathcal B(\beta_j)^{<L_j}\to\mathcal B(\beta_{j+1})^{<L_j}$
with respect to $\bar{\varepsilon}_{\beta_j}$ and
$\left.\bar{\varepsilon}_{\beta_{j+1}}\right|_{\mathcal B(\beta_{j+1})^{<L_j}}$
preserves the filtrations and induces $T_{m,j}$
on the $m$-th quotient.
The linearization of the inclusion
\[
 \mathcal B(\beta_{j+1})^{<L_j}\subset\mathcal B(\beta_{j+1})^{<L_{j+1}}
\]
also preserves the filtrations and, under
\eqref{eq:local-bw-filtration-quotient}, induces the
identity on $V_{m,j+1}$ for $1\le m\le\lfloor N_j/h\rfloor$.
Thus $Qf_j$ preserves the filtrations and induces
$T_{m,j}$ on the $m$-th quotient in this range.
\end{proof}

\Needspace{6\baselineskip}

\subsection{Homology of the linearized complexes}
\label{subsec:reference-homology}

Consider the algebras and maps of
Proposition~\ref{prop:reference-sequence},
with the augmentations of
Corollary~\ref{cor:reference-sequence-augmentations}.
Define the linearized complexes in Conley--Zehnder grading by
\begin{equation}\label{eq:reference-linearized-complexes}
 (C_j)_k
 =Q(\mathcal B(\beta_j)^{<L_j},\bar{\varepsilon}_{\beta_j})_{k+n-2}.
\end{equation}
The complexes $C_j$ form a direct system with transition maps
$Qf_j:C_j\to C_{j+1}$.
For each integer $m\ge0$, set
\begin{equation}\label{eq:reference-layer-filtration}
 \mathcal F_m C_j
 =\operatorname{span}_{\mathbb Q}
   \{[v-\bar{\varepsilon}_{\beta_j}(v)]:
      v\in V_{r,j},\ 1\le r\le m,\ rh\le N_j\}.
\end{equation}
By Lemma~\ref{lem:local-bw-linearization},
these are subcomplexes, with
$\mathcal F_0C_j=0$ and
$\mathcal F_{\lfloor N_j/h\rfloor}C_j=C_j$.

\begin{proposition}\label{prop:reference-homology}
For every $j\ge1$ and $k\in\mathbb Z$,
\begin{equation}\label{eq:finite-reference-homology}
 H_k(C_j)\cong
 \bigoplus_{m=1}^{\lfloor N_j/h\rfloor}
       H_{k-2\sigma m c_B+n}(B;\mathbb Q).
\end{equation}
The transition maps induce injections
$H_k(C_j)\to H_k(C_{j+1})$, and
\begin{equation}\label{eq:limitBW}
 H_k\!\left(\varinjlim_j C_j\right)
 \cong\bigoplus_{m\ge1}H_{k-2\sigma m c_B+n}(B;\mathbb Q).
\end{equation}
\end{proposition}

\begin{proof}
For each $j$, consider the spectral sequence associated with
the filtration \eqref{eq:reference-layer-filtration},
indexed by the filtration degree $m$ and the total
Conley--Zehnder degree $k$.
By Lemma~\ref{lem:local-bw-linearization} and
\eqref{eq:reference-layer-dga-homology}, its first page is
\begin{equation}\label{eq:pageone}
 E^1_{m,k}=H_{k-2\sigma m c_B+n}(B;\mathbb Q),
 \qquad 1\le m\le\lfloor N_j/h\rfloor.
\end{equation}
All other terms are zero.
Since the filtration is finite, the spectral sequence
converges to $H_*(C_j)$.

For $r\ge1$, the differential maps
$E^r_{m,k}$ to $E^r_{m-r,k-1}$.
If both corresponding $E^1$ terms were nonzero,
\eqref{eq:pageone} would give nonzero homology groups of $B$
in degrees differing by $2\sigma r c_B-1$,
contrary to Lemma~\ref{lem:betti-support}.
Since every later term is a subquotient of the corresponding
$E^1$ term, all differentials vanish.
Thus the spectral sequence collapses at $E^1$, and splitting
the induced filtration on homology gives
\eqref{eq:finite-reference-homology}.

By Corollary~\ref{cor:reference-sequence-augmentations},
$Qf_j$ preserves the filtrations and induces an isomorphism
on every $E^1$ term with
$1\le m\le\lfloor N_j/h\rfloor$.
All other source terms vanish.
Since both spectral sequences collapse at $E^1$,
naturality implies that the induced map on associated
graded homology is injective.
The homology filtrations are finite, so the induced map
$H_k(C_j)\to H_k(C_{j+1})$ is injective.
By \eqref{eq:finite-reference-homology}, its cokernel
has dimension
\[
 \sum_{\lfloor N_j/h\rfloor<m\le\lfloor N_{j+1}/h\rfloor}
 b_{k-2\sigma m c_B+n}(B).
\]

For each $k$, choose an isomorphism in
\eqref{eq:finite-reference-homology} at $j=1$ and,
successively, a complement to the image of
$H_k(C_j)$ in $H_k(C_{j+1})$.
By the preceding dimension formula, each complement can
be identified with the direct sum of
$H_{k-2\sigma m c_B+n}(B;\mathbb Q)$ over
$\lfloor N_j/h\rfloor<m\le\lfloor N_{j+1}/h\rfloor$.
Since $N_j\to\infty$, these choices identify
$\varinjlim_j H_k(C_j)$ with the direct sum
in \eqref{eq:limitBW}.
Exactness of filtered direct limits of vector spaces gives
$H_k(\varinjlim_j C_j)\cong\varinjlim_j H_k(C_j)$,
which proves \eqref{eq:limitBW}.
\end{proof}

\begin{remark}\label{rem:what-local-proof-does-not-claim}
The isomorphisms in \eqref{eq:finite-reference-homology}
and \eqref{eq:limitBW} are noncanonical.
\end{remark}

\subsection{\texorpdfstring{Proof of Theorem~\ref{thm:homology}}%
  {Proof of Theorem \getrefnumber{thm:homology}}}
\label{subsec:reference-comparison}
\label{subsec:homology-proof}

Consider the algebras and maps of
Proposition~\ref{prop:reference-sequence},
with the augmentations of
Corollary~\ref{cor:reference-sequence-augmentations}.
Choose a cobordism map
\[
 F_0:\mathcal B(\alpha)\longrightarrow\mathcal B(\beta_1)
\]
as in Lemmas~\ref{lem:found-contractible-quotient}
and~\ref{lem:found-action}.
The maps $F_j$, $j\ge0$, are quasi-isomorphisms that preserve
the integer grading. For every $L>0$, they satisfy
\[
 \begin{aligned}
  F_0\bigl(\mathcal B(\alpha)^{<L}\bigr)
  &\subset\mathcal B(\beta_1)^{<L},\\
  F_j\bigl(\mathcal B(\beta_j)^{<L}\bigr)
  &\subset\mathcal B(\beta_{j+1})^{<L}
  \qquad(j\ge1).
 \end{aligned}
\]
For $j\ge1$, set
\begin{equation}\label{eq:reference-composites}
 \Phi_j
 =F_{j-1}\circ\cdots\circ F_0:
 \mathcal B(\alpha)\longrightarrow\mathcal B(\beta_j).
\end{equation}
The identity $\Phi_{j+1}=F_j\circ\Phi_j$ gives
the commuting squares
\[
 \begin{array}{ccc}
 \mathcal B(\alpha)^{<L_j}
 &\xrightarrow{\ \Phi_j\ }&
 \mathcal B(\beta_j)^{<L_j}
 \\
 \big\downarrow&&\big\downarrow f_j
 \\
 \mathcal B(\alpha)^{<L_{j+1}}
 &\xrightarrow{\ \Phi_{j+1}\ }&
 \mathcal B(\beta_{j+1})^{<L_{j+1}}.
 \end{array}
\]

Set
\[
 A_\infty
 =\varinjlim_j\mathcal B(\beta_j)^{<L_j},
\]
with its induced integer grading.
By \eqref{eq:augmentation-compatibility}, the augmentations
$\bar{\varepsilon}_{\beta_j}$ induce an augmentation
$\bar{\varepsilon}_\infty:A_\infty\to\mathbb Q$.
The commuting squares induce a map $q$, and we equip
$\mathcal B(\alpha)$ with the pullback augmentation:
\begin{equation}\label{eq:reference-limit-map}
 q:\mathcal B(\alpha)
 =\varinjlim_j\mathcal B(\alpha)^{<L_j}
 \longrightarrow A_\infty,
 \qquad
 \bar{\varepsilon}_\alpha=\bar{\varepsilon}_\infty\circ q.
\end{equation}

Equivalently, for $x\in\mathcal B(\alpha)^{<L_j}$,
\[
 \bar{\varepsilon}_\alpha(x)=\bar{\varepsilon}_{\beta_j}(\Phi_j(x)).
\]
The commuting squares and
\eqref{eq:augmentation-compatibility} make this expression
independent of the choice of $j$.
Let $\varepsilon_\alpha$ be the augmentation of
$\mathcal C(\alpha)$ obtained by pulling back
$\bar{\varepsilon}_\alpha$ along
$\mathcal C(\alpha)\to\mathcal B(\alpha)$.
It vanishes on all generators corresponding to noncontractible Reeb orbits.

\begin{proposition}
\label{prop:found-colimit-qiso}
\label{thm:found-reference-comparison}
The map
$q:(\mathcal B(\alpha),\bar{\varepsilon}_\alpha)
\to(A_\infty,\bar{\varepsilon}_\infty)$
is a quasi-isomorphism and induces a
quasi-isomorphism that preserves the integer grading
\begin{equation}\label{eq:reference-linearization}
 Q(\mathcal B(\alpha),\bar{\varepsilon}_\alpha)
 \xrightarrow{\ Qq\ }
 \varinjlim_j
 Q(\mathcal B(\beta_j)^{<L_j},\bar{\varepsilon}_{\beta_j}).
\end{equation}
The contractible linearized complex at $\varepsilon_\alpha$
is naturally identified with $Q(\mathcal B(\alpha),\bar{\varepsilon}_\alpha)$
after the grading shift in
\eqref{eq:contractible-linearized-complex}.
\end{proposition}

\begin{proof}
Write $A_j=\mathcal B(\beta_j)^{<L_j}$.
To prove injectivity on homology, let
$x\in\mathcal B(\alpha)$ be a cycle whose image under $q$
is a boundary.
Choose $j$ such that $x\in\mathcal B(\alpha)^{<L_j}$.
For some $k\ge j$, there is $b\in A_k$ such that
$\Phi_k(x)=db$.
Since $\Phi_k$ is a quasi-isomorphism,
$[x]=0$ in $H(\mathcal B(\alpha))$.

For surjectivity, represent a class in $H(A_\infty)$
by a cycle $y\in A_j$.
Such a representative exists because filtered direct
limits of vector spaces are exact.
Since $\Phi_j$ is a quasi-isomorphism, there are a cycle
$x\in\mathcal B(\alpha)$ and an element
$z\in\mathcal B(\beta_j)$ such that
\[
 y=\Phi_j(x)+dz.
\]
Each of $x$ and $z$ involves only finitely many generators.
Choose $k>j$ large enough that every generator occurring
in $x$ has $\alpha$-action less than $L_k$ and every
generator occurring in $z$ has $\beta_j$-action less
than $L_k$.
Then $x\in\mathcal B(\alpha)^{<L_k}$, and
Lemma~\ref{lem:found-action} gives
\[
 (F_{k-1}\circ\cdots\circ F_j)(z)\in A_k.
\]
Applying this composition to the equality above shows
that the image of $y$ in $A_k$ is homologous to
$\Phi_k(x)$.
Thus the original class is $Hq([x])$, proving surjectivity.

To show that $A_\infty$ is Sullivan, choose a Sullivan
ordering of the orbit generators of $A_1$ by action.
By Proposition~\ref{prop:reference-sequence}\textup{(ii)},
$f_j$ identifies $A_j$ with the subalgebra
$\mathcal B(\beta_{j+1})^{<L_j}$ of $A_{j+1}$.
Use $f_j$ to transport the chosen generators of $A_j$
and their ordering.
The good contractible orbits of $\beta_{j+1}$ with action
in $[L_j,L_{j+1})$ freely generate $A_{j+1}$ over
$f_j(A_j)$.
Adjoin them in increasing order of action, breaking ties
arbitrarily.
The differential of each new generator lies in the
subalgebra generated by $f_j(A_j)$ and the earlier new
generators, since every orbit occurring in its
differential has strictly smaller action.
This extends the Sullivan ordering at each stage.
The union of these ordered generating sets therefore
gives a Sullivan ordering on $A_\infty$.

By Lemma~\ref{lem:found-contractible-quotient},
$\mathcal B(\alpha)$ is also Sullivan.
Moreover, $q$ preserves the augmentations by
\eqref{eq:reference-limit-map}.
Applying Lemma~\ref{lem:found-forward-sullivan} to $q$
therefore shows that $Qq$ is a quasi-isomorphism of
$\mathbb Z$-graded complexes.

Exactness of filtered direct limits gives
\[
 \varinjlim_j I_{A_j}\cong I_{A_\infty}.
\]
Under this identification,
$\varinjlim_j I_{A_j}^2$ is identified with
$I_{A_\infty}^2$, since every finite sum of products
in $I_{A_\infty}$ can be represented at a common stage.
Passing to the quotients gives the natural isomorphism
\[
 Q(A_\infty)\cong\varinjlim_j Q(A_j),
\]
and hence \eqref{eq:reference-linearization}.
The identification of the contractible linearized complex
for the pulled-back augmentation follows from
Section~\ref{subsec:contact-algebras}.
\end{proof}

\begin{proof}[Proof of Theorem~\ref{thm:homology}]
Let $\alpha$ be the given nondegenerate contact form.
Since $M$ is compact, choose $c>0$ such that
$c\alpha>2\beta$.
Apply the preceding construction to $c\alpha$, obtaining
an augmentation $\bar{\varepsilon}_{c\alpha}$ of
$\mathcal B(c\alpha)$ and its pullback
$\varepsilon_{c\alpha}$ to $\mathcal C(c\alpha)$.
For every $k\in\mathbb Z$, we have
\begin{equation}\label{eq:scaled-contact-homology}
 \begin{aligned}
  \LCH^0_k(M,c\alpha;\varepsilon_{c\alpha})
  &=H_{k+n-2}\bigl(Q(\mathcal B(c\alpha),\bar{\varepsilon}_{c\alpha})\bigr)\\
  &\cong H_k\!\left(\varinjlim_j C_j\right)\\
  &\cong\bigoplus_{m\ge1}
       H_{k-2\sigma mc_B+n}(B;\mathbb Q).
 \end{aligned}
\end{equation}
The equality follows from
\eqref{eq:contractible-linearized-complex}.
The first isomorphism follows from
Proposition~\ref{thm:found-reference-comparison} and the
definition of $C_j$ in \eqref{eq:reference-linearized-complexes},
and the second from
Proposition~\ref{prop:reference-homology}.

For compatible choices of the auxiliary data defining the
contact algebras, the correspondence of Reeb orbits under
constant rescaling gives a graded DGA isomorphism
\begin{equation}\label{eq:contact-rescaling}
 \mathcal B(\alpha)\cong\mathcal B(c\alpha).
\end{equation}
Pulling back $\bar{\varepsilon}_{c\alpha}$ along this isomorphism
gives an augmentation $\bar{\varepsilon}_\alpha$ of
$\mathcal B(\alpha)$.
Let $\varepsilon_\alpha$ be its pullback along
$\mathcal C(\alpha)\to\mathcal B(\alpha)$.
This augmentation vanishes on all generators corresponding to
noncontractible Reeb orbits.
The rescaling isomorphism \eqref{eq:contact-rescaling}
induces an isomorphism of the corresponding linearized complexes.
Thus \eqref{eq:scaled-contact-homology} gives
\eqref{eq:BWpattern} for the original contact form.
\end{proof}

\begin{remark}\label{rem:zero-augmentation}
If $\sigma=1$ and $c_B\ge2$, or if $\sigma=-1$ and
$c_B\ge n$, the augmentations
$\bar{\varepsilon}_{\beta_j}:
\mathcal B(\beta_j)^{<L_j}\to\mathbb Q$
in Corollary~\ref{cor:reference-sequence-augmentations}
may be chosen to vanish on all generators.
Indeed, \eqref{eq:contact-grading} and \eqref{eq:BWindex}
give
\[
 |\gamma_{x,\beta_\epsilon}^{mh}|
 =2\sigma m c_B-2+\operatorname{ind}_f(x),
\]
which is at least $2$ in the first case and at most $-2$
in the second.
These bounds persist under the perturbations of
Corollary~\ref{cor:finite-orbit-persistence}.
Since the differential has degree $-1$ and the cobordism
maps preserve degree, neither has a constant term on
generators. Thus sending every generator to zero defines
augmentations satisfying
\eqref{eq:augmentation-compatibility}.
\end{remark}

\section{Multiplicity of closed Reeb orbits}
\label{sec:complete-mixed-sign-count}

We first prove the multiplicity and non-hyperbolicity
bounds of Theorems~\ref{thm:main} and~\ref{thm:main-nonhyperbolic}
for $c_B>0$. We then treat $c_B=0$ in
Theorem~\ref{thm:zero-chern} and prequantizations of symplectic
orbifolds in Theorem~\ref{thm:main-orbifold}.

\subsection{\texorpdfstring{Proof of Theorem~\ref{thm:main}}%
  {Proof of Theorem \getrefnumber{thm:main}}}
\label{subsec:main-proof}

Let $\alpha$ be a nondegenerate contact form
supporting $\xi$, and fix an integer $p\ge0$. We assume that
$c_1(\xi)|_{\pi_2(M)}=0$, $c_B>0$, and
\[
 b_j(B)b_{j+2\sigma c_B-1}(B)=0
 \qquad\text{for every }j\in\mathbb Z.
\]

The conclusion is immediate if $r_\alpha=\infty$.
We may therefore assume the following condition.
\begin{enumerate}
\renewcommand{\theenumi}{FCI}
\item\label{ass:count-fci}
The contact form $\alpha$ has finitely many geometrically distinct
contractible closed Reeb orbits. Let $\gamma_1,\ldots,\gamma_r$ denote
all geometrically distinct good least contractible closed Reeb
orbits of $\alpha$. For each $i$, we have
$\widehat\mu(\gamma_i)\ne0$ and
$\mu(\gamma_i^k)\notin\{-p,p\}$ whenever $k\ge1$ and
$\gamma_i^k$ is good.
\end{enumerate}

Set $\sigma_i=\operatorname{sign}\widehat\mu(\gamma_i)\in\{-1,1\}$.
Since $r\le r_\alpha$, it suffices to prove the lower bounds of
Theorem~\ref{thm:main} for $r$.

Lemma~\ref{lem:topology}\textup{(i)} gives
$\omega|_{\pi_2(B)}\ne0$. We can therefore choose an augmentation
$\varepsilon_\alpha$ as in Theorem~\ref{thm:homology}.
With the notation of Section~\ref{subsec:morse-mean-euler},
\eqref{eq:BWpattern} gives
\begin{equation}\label{eq:count-bw-pattern}
 b_q=\sum_{m\ge1}b_{q-2\sigma m c_B+n}(B)
 \qquad(q\in\mathbb Z).
\end{equation}
By \eqref{ass:count-fci}, we have $c_p=c_{-p}=0$.
Moreover, \eqref{eq:count-bw-pattern} and
Proposition~\ref{prop:finite-morse} give
$c_{2\sigma c_B-n}\ge b_{2\sigma c_B-n}\ge b_0(B)=1$,
so $r\ge1$.

Choose $0<\eta<1$ such that
\[
 \eta\sum_{i=1}^r
 \left|\frac{\widehat\chi(\gamma_i)}
                 {\widehat\mu(\gamma_i)}\right|<1.
\]
Next, choose $\ell_0\in\mathbb N$ sufficiently large that
\[
 \ell_0\min_{1\le i\le r}|\widehat\mu(\gamma_i)|>n+p+2.
\]
By \eqref{eq:count-mean-estimate}, any iterate $\gamma_i^\ell$
with $|\mu(\gamma_i^\ell)|\le p$ then satisfies $\ell<\ell_0$.

Using trivializations induced by capping disks, apply
Lemma~\ref{lem:index-recurrence} to the symplectic paths
associated with the linearized Reeb flows along the $\gamma_i$,
with parameters $\eta,\ell_0$ and $N=2c_B$.

Choose $j$ sufficiently large that
\[
 \begin{aligned}
 k_{ij}^{\pm}&>2\ell_0
 &&(1\le i\le r),\\
 d_j^\pm&>n+p+2+
 \max_{\substack{1\le i\le r\\1\le\ell\le\ell_0}}
 |\mu(\gamma_i^\ell)|.
 \end{aligned}
\]
Fix this $j$ and write
\[
 (d,k_1,\ldots,k_r)=(d_j^+,k_{1j}^+,\ldots,k_{rj}^+),\qquad
 (d',k'_1,\ldots,k'_r)=(d_j^-,k_{1j}^-,\ldots,k_{rj}^-).
\]
By Lemma~\ref{lem:index-recurrence}, $d=2sc_B$ and $d'=2s'c_B$
for some $s,s'\in\mathbb N$, and all $k_i,k'_i$ are even.

For each $j=0,\ldots,2n$, the copies of $H_j(B;\mathbb Q)$ in
\eqref{eq:BWpattern} occur in degrees
$j+2\sigma mc_B-n$, $m\ge1$.
Fix an integer $L>n$. Such a degree lies in $[-L,L]$ precisely when
\[
 -L+\sigma(n-j)\le 2mc_B\le L+\sigma(n-j).
\]
Since
$-L+\sigma(n-j)\le -L+|n-j|\le -L+n<0$,
the inequality $-L+\sigma(n-j)\le2mc_B$ holds automatically
for every $m\ge1$. Thus the condition reduces to
\[
 1\le m\le\frac{L+\sigma(n-j)}{2c_B}.
\]
The number of integers $m$ in this range is
$\lfloor(L+\sigma(n-j))/(2c_B)\rfloor$.

Therefore,
\[
 \begin{aligned}
 \lim_{L\to\infty}\frac{1}{L}
 \sum_{q=-L}^{L}(-1)^q b_q
 &=\lim_{L\to\infty}\frac{(-1)^n}{L}
   \sum_{j=0}^{2n}(-1)^j b_j(B)
   \left\lfloor\frac{L+\sigma(n-j)}{2c_B}\right\rfloor\\
 &=\frac{(-1)^n}{2c_B}
   \sum_{j=0}^{2n}(-1)^j b_j(B)\\
 &=\frac{(-1)^n\chi(B)}{2c_B}.
 \end{aligned}
\]

Using \eqref{eq:mean-euler-definition} and the resonance relations
\eqref{eq:orbit-resonance}, we obtain
\begin{equation}\label{eq:count-total-resonance}
 \begin{aligned}
 \sum_i\frac{\widehat\chi(\gamma_i)}{|\widehat\mu(\gamma_i)|}
 &=\chi_+(M,\alpha;\varepsilon_\alpha)
   +\chi_-(M,\alpha;\varepsilon_\alpha)\\
 &=\lim_{L\to\infty}\frac{1}{L}
   \sum_{q=-L}^{L}(-1)^q b_q\\
 &=\frac{(-1)^n\chi(B)}{2c_B}.
 \end{aligned}
\end{equation}

For $(d,k_1,\ldots,k_r)$, we have
\begin{equation}\label{eq:count-exact-euler}
 \sum_{i=1}^r\sum_{\ell=1}^{k_i}\chi(\gamma_i^\ell)
 =\sum_{i=1}^r k_i\widehat\chi(\gamma_i)=(-1)^ns\chi(B).
\end{equation}
Indeed, the first equality in \eqref{eq:count-exact-euler}
follows from \eqref{eq:count-period-two}, since each $k_i$ is even.
In particular, $\sum_i k_i\widehat\chi(\gamma_i)$ is an integer.
Multiplying \eqref{eq:count-total-resonance} by $d=2sc_B$
and applying
Lemma~\ref{lem:index-recurrence}\eqref{eq:count-recurrence-mean},
we obtain
\[
 \left|\sum_i k_i\widehat\chi(\gamma_i)-(-1)^ns\chi(B)\right|
 \le \eta\sum_i
 \left|\frac{\widehat\chi(\gamma_i)}
                  {\widehat\mu(\gamma_i)}\right|<1.
\]
Since $(-1)^ns\chi(B)$ is also an integer, the second equality
follows. The same argument gives \eqref{eq:count-exact-euler}
with $k'_i,s'$ in place of $k_i,s$.

We define the following counts for $(d,k_1,\ldots,k_r)$,
using the notation of \cite[(6.18)--(6.19)]{AbM26}.
Let $r_+^{e,p}$ and $r_+^{o,p}$ denote the numbers of good
iterates $\gamma_i^{k_i}$ with even and odd index, respectively,
that satisfy
\[
 \sigma_i\mu(\gamma_i^{k_i})>d+p.
\]
Similarly, let $r_-^{e,p}$ and $r_-^{o,p}$ denote the numbers
of good iterates $\gamma_i^{k_i}$ with even and odd index,
respectively, that satisfy
\[
 \sigma_i\mu(\gamma_i^{k_i})<d-p.
\]
For every $p\ge0$, define $r_\pm^{e,p-1}$ and $r_\pm^{o,p-1}$
as $r_\pm^{e,p}$ and $r_\pm^{o,p}$, with $p-1$
in place of $p$.

\begin{lemma}\label{lem:count-iterates}
For $(d,k_1,\ldots,k_r)$, we have
\begin{equation}\label{eq:count-chain-window}
 \sum_{|q|\le d+p}(-1)^q c_q
 =\sum_{i=1}^r k_i\widehat\chi(\gamma_i)
 +\sum_{|{q}|\le p}(-1)^{q} c_{q}-r_+^{e,p}+r_+^{o,p}.
\end{equation}
We also have
\begin{equation}\label{eq:count-chain-shortened}
 \sum_{|q|\le d+p-1}(-1)^q c_q
 =\sum_{i=1}^r k_i\widehat\chi(\gamma_i)
  +\sum_{|q|\le p}(-1)^q c_q
  -r_+^{e,p-1}+r_+^{o,p-1}.
\end{equation}
Both identities hold with
$d,k_i$ replaced by $d',k'_i$ and subscript $+$ replaced by $-$
in the orbit counts.
\end{lemma}

\begin{proof}
By the definitions of $c_q$
and the local Euler characteristic,
\[
 \sum_{|q|\le d+p}(-1)^q c_q
 =\sum_{i=1}^r
   \sum_{\substack{m\ge1\\|\mu(\gamma_i^m)|\le d+p}}
   \chi(\gamma_i^m).
\]
Fix $i\in\{1,\ldots,r\}$. By \eqref{eq:count-mean-estimate}
and Lemma~\ref{lem:index-recurrence}\eqref{eq:count-recurrence-mean},
the choices of $\eta$ and $\ell_0$ give
\begin{equation}\label{eq:count-outer-indices}
 \begin{aligned}
 |\mu(\gamma_i^m)|&<d-p
 &&\text{if }1\le m<k_i-\ell_0,\\
 |\mu(\gamma_i^m)|&>d+p
 &&\text{if }m>k_i+\ell_0.
 \end{aligned}
\end{equation}
It remains to determine the contributions of $\gamma_i^{k_i}$
and of the pairs $\gamma_i^{k_i-\ell},\gamma_i^{k_i+\ell}$,
where $1\le\ell\le\ell_0$.

For $1\le\ell\le\ell_0$, Lemma~\ref{lem:index-recurrence}
\eqref{eq:count-recurrence-offset} gives
\[
 \mu(\gamma_i^{k_i\pm\ell})
 =\sigma_i d\pm\mu(\gamma_i^\ell).
\]
The choice of $d$ ensures that these indices have sign $\sigma_i$.
Their absolute values, in some order, are
\[
 d-|\mu(\gamma_i^\ell)|
 \qquad\text{and}\qquad
 d+|\mu(\gamma_i^\ell)|.
\]
Since $k_i$ is even, these two iterates are good precisely when
$\gamma_i^\ell$ is good, and \eqref{eq:count-period-two} gives
\[
 \chi(\gamma_i^{k_i\pm\ell})=\chi(\gamma_i^\ell).
\]
At least one of $\gamma_i^{k_i-\ell}$ and $\gamma_i^{k_i+\ell}$
has index in $[-d-p,d+p]$, and both do exactly when
$|\mu(\gamma_i^\ell)|\le p$. Consequently,
\[
 \sum_{\substack{m\in\{k_i-\ell,k_i+\ell\}\\
                 |\mu(\gamma_i^m)|\le d+p}}
 \chi(\gamma_i^m)
 =\begin{cases}
   \chi(\gamma_i^\ell), & |\mu(\gamma_i^\ell)|>p,\\
   2\chi(\gamma_i^\ell), & |\mu(\gamma_i^\ell)|\le p.
  \end{cases}
\]

Each pair therefore contributes $\chi(\gamma_i^{k_i-\ell})$,
with an additional $\chi(\gamma_i^\ell)$ when
$|\mu(\gamma_i^\ell)|\le p$.
Together with the contributions from $1\le m<k_i-\ell_0$,
the terms $\chi(\gamma_i^{k_i-\ell})$ thus give
$\sum_{m=1}^{k_i-1}\chi(\gamma_i^m)$.
The remaining iterate $\gamma_i^{k_i}$ contributes
$\chi(\gamma_i^{k_i})$ if its index lies in $[-d-p,d+p]$,
and zero otherwise.
We may therefore include $\chi(\gamma_i^{k_i})$ in the sum
and subtract it when $|\mu(\gamma_i^{k_i})|>d+p$.
Summing over $i$, we obtain
\begin{equation}\label{eq:count-euler-decomposition}
 \begin{aligned}
 \sum_{|q|\le d+p}(-1)^q c_q
 &=\sum_{i=1}^r\sum_{m=1}^{k_i}\chi(\gamma_i^m)\\
 &\quad+\sum_{i=1}^r
   \sum_{\substack{1\le\ell\le\ell_0\\
                   |\mu(\gamma_i^\ell)|\le p}}
   \chi(\gamma_i^\ell)\\
 &\quad-\sum_{\substack{1\le i\le r\\
                   |\mu(\gamma_i^{k_i})|>d+p}}
   \chi(\gamma_i^{k_i}).
 \end{aligned}
\end{equation}

To evaluate the second term on the right-hand side of
\eqref{eq:count-euler-decomposition}, recall that every iterate
with $|\mu(\gamma_i^\ell)|\le p$ satisfies $\ell<\ell_0$.
Hence, by the definitions of $c_{q}$ and the local Euler characteristic,
\[
 \sum_{i=1}^r
 \sum_{\substack{1\le\ell\le\ell_0\\
                 |\mu(\gamma_i^\ell)|\le p}}
 \chi(\gamma_i^\ell)
 =\sum_{|{q}|\le p}(-1)^{q} c_{q}.
\]

To evaluate the last sum in \eqref{eq:count-euler-decomposition},
we first use \eqref{eq:count-mean-estimate} and
Lemma~\ref{lem:index-recurrence}\eqref{eq:count-recurrence-mean}
to obtain
\[
 \sigma_i\mu(\gamma_i^{k_i})>d-\eta-n>0.
\]
Thus $|\mu(\gamma_i^{k_i})|>d+p$ is equivalent to
$\sigma_i\mu(\gamma_i^{k_i})>d+p$.
Among the iterates satisfying this inequality, good iterates
of even index contribute $1$, good iterates of odd index
contribute $-1$, and bad iterates contribute zero.
By the definitions of $r_+^{e,p}$ and $r_+^{o,p}$, we therefore have
\[
 \sum_{\substack{1\le i\le r\\
                 |\mu(\gamma_i^{k_i})|>d+p}}
 \chi(\gamma_i^{k_i})
 =r_+^{e,p}-r_+^{o,p}.
\]
Substituting the formulas for the two correction terms into
\eqref{eq:count-euler-decomposition} and using
\eqref{eq:count-exact-euler} proves
\eqref{eq:count-chain-window}.

We next prove \eqref{eq:count-chain-shortened}.
Suppose that $\gamma_i^m$ is good, $m\ne k_i$, and
$|\mu(\gamma_i^m)|=d+p$.
By \eqref{eq:count-outer-indices}, we have
$m=k_i\pm\ell$ for some $1\le\ell\le\ell_0$.
Since $\mu(\gamma_i^{k_i\pm\ell})$ has sign $\sigma_i$,
Lemma~\ref{lem:index-recurrence}
\eqref{eq:count-recurrence-offset} gives
$|\mu(\gamma_i^\ell)|=p$.
Moreover, $\gamma_i^\ell$ is good because $k_i$ is even,
contrary to \eqref{ass:count-fci}.
Thus every good iterate in degrees $\pm(d+p)$ is among
$\gamma_1^{k_1},\ldots,\gamma_r^{k_r}$.

Since the indices are integers, replacing $p$ by $p-1$
in the counts $r_+^{e,p},r_+^{o,p}$ adds precisely the good
iterates $\gamma_i^{k_i}$ with
$|\mu(\gamma_i^{k_i})|=d+p$.
Since $d$ is even,
$(-1)^{d+p}=(-1)^{-d-p}=(-1)^p$, so
\[
 (-1)^p(c_{d+p}+c_{-d-p})
 =
 \bigl(r_+^{e,p-1}-r_+^{o,p-1}\bigr)
 -
 \bigl(r_+^{e,p}-r_+^{o,p}\bigr).
\]
Subtracting $(-1)^p(c_{d+p}+c_{-d-p})$ from
\eqref{eq:count-chain-window} therefore proves
\eqref{eq:count-chain-shortened}.

The same arguments apply to $(d',k'_1,\ldots,k'_r)$.
To identify the corresponding counts, we use
Lemma~\ref{lem:index-recurrence}
\eqref{eq:count-recurrence-companion}, which gives
\[
 \sigma_i\mu(\gamma_i^{k'_i})-d'
 =-\bigl(\sigma_i\mu(\gamma_i^{k_i})-d\bigr).
\]
Consequently,
\[
 \begin{aligned}
 \sigma_i\mu(\gamma_i^{k'_i})>d'+p
 &\quad\Longleftrightarrow\quad
 \sigma_i\mu(\gamma_i^{k_i})<d-p,\\
 \sigma_i\mu(\gamma_i^{k'_i})>d'+p-1
 &\quad\Longleftrightarrow\quad
 \sigma_i\mu(\gamma_i^{k_i})<d-p+1.
 \end{aligned}
\]
Since $k_i$ and $k'_i$ are even,
$\gamma_i^{k_i}$ and $\gamma_i^{k'_i}$ are either both good
or both bad, and their indices have the same parity.
The definitions of
$r_-^{e,p},r_-^{o,p},r_-^{e,p-1},r_-^{o,p-1}$
therefore show that the correction terms for
$(d',k'_1,\ldots,k'_r)$ are
\[
 -r_-^{e,p}+r_-^{o,p}
 \qquad\text{and}\qquad
 -r_-^{e,p-1}+r_-^{o,p-1}
\]
in \eqref{eq:count-chain-window} and
\eqref{eq:count-chain-shortened}, respectively.
\end{proof}

Define
\[
 A_p=-2(-1)^{n+p}\sum_{j=0}^{n-p-1}(-1)^j b_j(B)
 \qquad(p\ge0),
\]
where empty sums are understood to be zero.
The same quantity $A_p$ appears in
Theorem~\ref{thm:main-nonhyperbolic}.

For $q\in\mathbb Z$, put
\begin{equation}\label{eq:count-periodic-ranks}
 E_q=\sum_{m\in\mathbb Z}b_{n+q-2mc_B}(B).
\end{equation}

\begin{lemma}\label{lem:count-ranks}
We have
\begin{equation}\label{eq:count-homology-window}
 \sum_{|q|\le d+p}(-1)^q b_q
 =s(-1)^n\chi(B)+\sum_{|q|\le p}(-1)^q b_q
  +\frac{(-1)^pA_p}{2}.
\end{equation}
We also have
\begin{equation}\label{eq:count-shortened-windows}
 \sum_{|q|\le d+p-1}(-1)^q b_q
 =(-1)^ns\chi(B)+\sum_{|q|\le p}(-1)^q b_q
  +(-1)^p\left(\frac{A_p}{2}-E_p\right).
\end{equation}
Both identities hold with $d',s'$ in place of $d,s$.
\end{lemma}

\begin{proof}
By \eqref{eq:count-bw-pattern} and Poincar\'e duality,
\[
 b_{-q}
 =\sum_{m\ge1}b_{n-q-2\sigma mc_B}(B)
 =\sum_{m\ge1}b_{n+q+2\sigma mc_B}(B).
\]
Since $\sigma\in\{-1,1\}$, separating the term with $m=0$
in \eqref{eq:count-periodic-ranks} gives
\begin{equation}\label{eq:nonhyperbolic-rank-splitting}
 \begin{aligned}
 E_q
 &=b_{n+q}(B)
   +\sum_{m\ge1}b_{n+q-2\sigma mc_B}(B)
   +\sum_{m\ge1}b_{n+q+2\sigma mc_B}(B)\\
 &=b_q+b_{-q}+b_{n-q}(B).
 \end{aligned}
\end{equation}
Here we used $b_{n+q}(B)=b_{n-q}(B)$.
The same identity and Poincar\'e duality give $E_{-q}=E_q$.

By \eqref{eq:count-periodic-ranks}, we have
$E_{q+2c_B}=E_q$. Since $2c_B$ is even, the sequence
$(-1)^qE_q$ also has period $2c_B$. Moreover,
\[
 \begin{aligned}
 \sum_{q=0}^{2c_B-1}(-1)^qE_q
 &=\sum_{q=0}^{2c_B-1}\sum_{m\in\mathbb Z}
   (-1)^q b_{n+q-2mc_B}(B)\\
 &=(-1)^n\sum_{j=0}^{2n}(-1)^j b_j(B)\\
 &=(-1)^n\chi(B).
 \end{aligned}
\]

Each of the intervals $[-d-p,-p-1]$ and $[p+1,d+p]$
contains $d=2sc_B$ consecutive integers. The sum of
$(-1)^qE_q$ over each interval is therefore
$s(-1)^n\chi(B)$.

By the change of variable $q\mapsto-q$, we have
\[
 \sum_{p<|q|\le d+p}(-1)^q b_{-q}
 =\sum_{p<|q|\le d+p}(-1)^q b_q.
\]
Also, $b_{n-q}(B)=0$ for $|q|>n$, so $d+p>n$ gives
\[
 \sum_{|q|\le d+p}(-1)^q b_{n-q}(B)
 =(-1)^n\sum_{j=0}^{2n}(-1)^j b_j(B)
 =(-1)^n\chi(B).
\]
Thus \eqref{eq:nonhyperbolic-rank-splitting} yields
\begin{equation}\label{eq:count-period-sum}
 \begin{aligned}
 2s(-1)^n\chi(B)
 &=\sum_{p<|q|\le d+p}(-1)^qE_q\\
 &=2\sum_{p<|q|\le d+p}(-1)^q b_q
   +(-1)^n\chi(B)
   -\sum_{|q|\le p}(-1)^q b_{n-q}(B).
 \end{aligned}
\end{equation}

By Poincar\'e duality and the change of variable $j=n-q$,
\[
 \begin{aligned}
 \sum_{|q|>p}(-1)^q b_{n-q}(B)
 &=2\sum_{q>p}(-1)^q b_{n-q}(B)\\
 &=2(-1)^n\sum_{j=0}^{n-p-1}(-1)^j b_j(B)\\
 &=-(-1)^pA_p.
 \end{aligned}
\]
Consequently,
\begin{equation}\label{eq:count-base-window}
 \begin{aligned}
 (-1)^p\sum_{|q|\le p}(-1)^q b_{n-q}(B)
 &=(-1)^p\left(
   (-1)^n\chi(B)-\sum_{|q|>p}(-1)^q b_{n-q}(B)
   \right)\\
 &=A_p+(-1)^{n+p}\chi(B).
 \end{aligned}
\end{equation}
Substituting \eqref{eq:count-base-window} into
\eqref{eq:count-period-sum}, we obtain
\[
 \sum_{p<|q|\le d+p}(-1)^q b_q
 =s(-1)^n\chi(B)+\frac{(-1)^pA_p}{2}.
\]
Adding $\sum_{|q|\le p}(-1)^q b_q$ to both sides proves
\eqref{eq:count-homology-window}.

Moreover, \eqref{eq:nonhyperbolic-rank-splitting},
$d+p>n$, and $d=2sc_B$ give
\[
 b_{d+p}+b_{-d-p}=E_{d+p}=E_p.
\]
Subtracting $(-1)^pE_p$ from
\eqref{eq:count-homology-window} proves
\eqref{eq:count-shortened-windows}.
The same argument applies with $d',s'$ in place of $d,s$.
\end{proof}

By \eqref{ass:count-fci}, $c_p=c_{-p}=0$.
The Morse inequalities \eqref{eq:count-finite-morse},
together with \eqref{eq:nonhyperbolic-rank-splitting}
at $q=p$, give
\begin{equation}\label{eq:count-hole-rank}
 E_p=b_{n-p}(B),\qquad
 \sum_{|q|\le p}(-1)^q c_q
 =\sum_{|q|\le p}(-1)^q b_q.
\end{equation}

\Needspace{9\baselineskip}
We will use the following lemma for $n=2$ in the
proofs of Theorems~\ref{thm:main} and~\ref{thm:main-nonhyperbolic}.

\begin{lemma}\label{lem:five-dimensional-index-zero}
Let $\alpha$ be a nondegenerate contact form on $(M^5,\xi)$,
where $c_1(\xi)|_{\pi_2(M)}=0$, and let $\gamma$ be a good
contractible closed Reeb orbit.
\begin{enumerate}
\item If $\mu(\gamma)=0$ and $\widehat\mu(\gamma)\ne0$, then
for every positive even integer $p$, there is a good
iterate $\gamma^k$ such that
\[
 \mu(\gamma^k)=\operatorname{sign}\widehat\mu(\gamma)\,p.
\]
\item If $\mu(\gamma)\in\{-1,1\}$, then for every positive
odd integer $p$, there is a good iterate $\gamma^k$ such that
\[
 \mu(\gamma^k)=p\mu(\gamma).
\]
\end{enumerate}
In particular, suppose that, for some positive integer $p$,
every good contractible closed Reeb orbit has nonzero mean
index and index different from $\pm p$.
If $p$ is even, no such orbit has index zero.
If $p$ is odd, no such orbit has index $\pm1$.
\end{lemma}

\begin{proof}
The nondegenerate index iteration formula
\cite[Theorem~2.1]{GK16} gives integers $a$ and
$e\in\{0,1,2\}$ and irrational numbers
$\theta_1,\ldots,\theta_e\in(0,1)$ such that
\begin{equation}\label{eq:five-dimensional-iteration}
 \mu(\gamma^k)
 =ka+2\sum_{j=1}^{e}\lfloor k\theta_j\rfloor+e
 \qquad(k\ge1).
\end{equation}
Moreover, $a$ is even when $e=2$.
Since $\gamma$ is good, its index has the same parity
as that of its underlying simple orbit.
The definition of good orbits and
\eqref{eq:index-trivialization-change} therefore give
\[
 \gamma^k\text{ is good}
 \quad\Longleftrightarrow\quad
 \mu(\gamma^k)\equiv\mu(\gamma)\pmod2.
\]

For \textup{(i)}, $\mu(\gamma)=0$ gives $a=-e$.
If $e=0$, all iterates have index zero, contradicting
$\widehat\mu(\gamma)\ne0$.
Thus $e=1$ or $e=2$.
If $e=2$, every iterate has even index, and
\eqref{eq:five-dimensional-iteration} gives
\[
 \mu(\gamma^{k+1})-\mu(\gamma^k)
 =-2+2\sum_{j=1}^{2}
   \bigl(\lfloor(k+1)\theta_j\rfloor
         -\lfloor k\theta_j\rfloor\bigr)
 \in\{-2,0,2\}.
\]
If $e=1$, every odd iterate has even index, and
\[
 \mu(\gamma^{2m+3})-\mu(\gamma^{2m+1})
 =-2+2\bigl(\lfloor(2m+3)\theta_1\rfloor
           -\lfloor(2m+1)\theta_1\rfloor\bigr)
 \in\{-2,0,2\}
 \qquad(m\ge0).
\]
In both cases, we obtain a sequence of indices of good
iterates that starts at zero and has successive differences
in $\{-2,0,2\}$.
By \eqref{eq:count-mean-estimate}, this sequence tends to
$+\infty$ if $\widehat\mu(\gamma)>0$ and to $-\infty$
if $\widehat\mu(\gamma)<0$.
It therefore attains
$\operatorname{sign}\widehat\mu(\gamma)\,p$ for every
positive even integer $p$.

For \textup{(ii)}, $a+e=\mu(\gamma)\in\{-1,1\}$.
Since $a$ is even when $e=2$, we must have $e\le1$.
If $e=0$, then $\mu(\gamma^k)=k\mu(\gamma)$.
Taking $k=p$ gives the required index, and $\gamma^p$
is good because $p$ is odd.
If $e=1$, \eqref{eq:five-dimensional-iteration} becomes
\[
 \mu(\gamma^k)=
 \begin{cases}
  2\lfloor k\theta_1\rfloor+1,&\mu(\gamma)=1,\\
  -2k+2\lfloor k\theta_1\rfloor+1,&\mu(\gamma)=-1.
 \end{cases}
\]
Every iterate has odd index and is therefore good.
When $\mu(\gamma)=1$, these indices start at $1$,
tend to $+\infty$, and increase by $0$ or $2$ at each step.
When $\mu(\gamma)=-1$, they start at $-1$,
tend to $-\infty$, and decrease by $0$ or $2$ at each step.
Thus they attain $p\mu(\gamma)$ for every positive
odd integer $p$.
\end{proof}

\begin{proof}[Proof of Theorem~\ref{thm:main}]
\textit{Step 1.} For each integer $q$ with $|q|\le p$,
the number of good iterates $\gamma_i^{k_i}$ satisfying
$\mu(\gamma_i^{k_i})\in\{d+q,-d+q\}$ is
\begin{equation}\label{eq:count-central-ranks}
 c_{d+q}+c_{-d+q}-c_q-c_{-q}.
\end{equation}

By \eqref{eq:count-outer-indices}, an iterate $\gamma_i^m$
with $m\ne k_i$ and
$\mu(\gamma_i^m)\in\{d+q,-d+q\}$ must satisfy
$m=k_i\pm\ell$ for some $1\le\ell\le\ell_0$.
For $1\le\ell\le\ell_0$, the choice of $d$ ensures that
$\mu(\gamma_i^{k_i\pm\ell})$ has sign $\sigma_i$.
Since $d>p$, Lemma~\ref{lem:index-recurrence}
\eqref{eq:count-recurrence-offset} gives
\begin{equation}\label{eq:count-central-correspondence}
 \begin{aligned}
 \mu(\gamma_i^{k_i+\ell})\in\{d+q,-d+q\}
 &\quad\Longleftrightarrow\quad
 \mu(\gamma_i^\ell)=q,\\
 \mu(\gamma_i^{k_i-\ell})\in\{d+q,-d+q\}
 &\quad\Longleftrightarrow\quad
 \mu(\gamma_i^\ell)=-q.
 \end{aligned}
\end{equation}
Since $k_i$ is even, \eqref{eq:count-period-two} implies that
$\gamma_i^{k_i\pm\ell}$ is good if and only if
$\gamma_i^\ell$ is good.
Moreover, \eqref{eq:count-mean-estimate} and the choice of
$\ell_0$ imply $\ell<\ell_0$ whenever
$\mu(\gamma_i^\ell)\in\{q,-q\}$.
Thus \eqref{eq:count-central-correspondence} gives
$c_q$ good iterates of the form $\gamma_i^{k_i+\ell}$
and $c_{-q}$ of the form $\gamma_i^{k_i-\ell}$
with index in $\{d+q,-d+q\}$.
Since $k_i-\ell<k_i<k_i+\ell'$ for all
$1\le\ell,\ell'\le\ell_0$, the two lines of
\eqref{eq:count-central-correspondence} count distinct iterates,
even when $q=0$.

Since $c_{d+q}+c_{-d+q}$ counts all good iterates with index
in $\{d+q,-d+q\}$, subtracting $c_q+c_{-q}$ proves
\eqref{eq:count-central-ranks}.
The same argument proves \eqref{eq:count-central-ranks}
with $d,k_i$ replaced by $d',k'_i$.

We next compute $b_{d+q}+b_{-d+q}$.
For $|q|\le p$, \eqref{eq:count-bw-pattern} and $d>n+p$ give
\[
 b_{-\sigma d+q}=b_{-\sigma d-q}=0,
 \qquad
 b_{n-\sigma d-q}(B)=0.
\]
Since $d=2sc_B$, \eqref{eq:count-periodic-ranks} gives
$E_{\sigma d+q}=E_q$.
Applying \eqref{eq:nonhyperbolic-rank-splitting}
at $\sigma d+q$, we obtain
\[
 \begin{aligned}
 b_{\sigma d+q}
 &=E_{\sigma d+q}
   -b_{-\sigma d-q}-b_{n-\sigma d-q}(B)\\
 &=E_q.
 \end{aligned}
\]
Therefore,
\begin{equation}\label{eq:count-large-band}
 b_{d+q}+b_{-d+q}=E_q\qquad(|q|\le p).
\end{equation}
The same calculation proves \eqref{eq:count-large-band}
with $d$ replaced by $d'$.

\textit{Step 2.} We prove the bound
\begin{equation}\label{eq:count-final-bound}
 r\ge
 \begin{cases}
  -(-1)^n\chi(B)+2b_n(B),&p=0,\\
  -(-1)^{n+p}\chi(B)+4b_{n-p}(B),&p>0.
 \end{cases}
\end{equation}
We follow the argument in \cite[proof of Theorem~5]{AbM26}.

We first estimate the number of good iterates
$\gamma_i^{k_i}$ satisfying
$|\mu(\gamma_i^{k_i})|\notin[d-p,d+p]$.

By \eqref{eq:count-mean-estimate} and
Lemma~\ref{lem:index-recurrence}\eqref{eq:count-recurrence-mean},
we have
$|\mu(\gamma_i^{k_i})-\sigma_i d|<n+\eta<d$,
so $|\mu(\gamma_i^{k_i})|=\sigma_i\mu(\gamma_i^{k_i})$.
Thus $r_+^{e,p}+r_+^{o,p}$ and $r_-^{e,p}+r_-^{o,p}$
count the good iterates satisfying
$|\mu(\gamma_i^{k_i})|>d+p$ and
$|\mu(\gamma_i^{k_i})|<d-p$, respectively.

Since $d$ is even, applying \eqref{eq:count-finite-morse}
on $[-d-p,d+p]$ gives
\[
 0\le (-1)^p\sum_{|q|\le d+p}(-1)^q(c_q-b_q).
\]
Substituting \eqref{eq:count-chain-window},
\eqref{eq:count-exact-euler}, and \eqref{eq:count-homology-window},
and using \eqref{eq:count-hole-rank}, we obtain
\[
 (-1)^p(r_+^{o,p}-r_+^{e,p})\ge\frac{A_p}{2}.
\]

For $(d',k'_1,\ldots,k'_r)$, Lemma~\ref{lem:count-iterates}
gives \eqref{eq:count-chain-window} with
$r_-^{e,p},r_-^{o,p}$ in place of $r_+^{e,p},r_+^{o,p}$.
The same calculation gives
\[
 (-1)^p(r_-^{o,p}-r_-^{e,p})\ge\frac{A_p}{2}.
\]
Since $r_\pm^{e,p},r_\pm^{o,p}\ge0$, it follows that
\begin{equation}\label{eq:count-outside-bound}
 \begin{aligned}
 r_+^{e,p}+r_+^{o,p}+r_-^{e,p}+r_-^{o,p}
 &\ge (-1)^p
 \bigl(r_+^{o,p}-r_+^{e,p}+r_-^{o,p}-r_-^{e,p}\bigr)\\
 &\ge A_p.
 \end{aligned}
\end{equation}

It remains to estimate the number of good iterates
$\gamma_i^{k_i}$ satisfying
$|\mu(\gamma_i^{k_i})|\in[d-p,d+p]$.
Summing \eqref{eq:count-central-ranks}
over $|q|\le p$ and using \eqref{eq:count-outside-bound},
we obtain
\begin{equation}\label{eq:count-band-bound}
 r\ge A_p+\sum_{|q|\le p}
       (c_{d+q}+c_{-d+q}-c_q-c_{-q}).
\end{equation}

If $p=0$, then $c_0=0$. By the Morse inequalities
\eqref{eq:count-finite-morse},
\eqref{eq:count-large-band}, and \eqref{eq:count-hole-rank},
\[
 c_d+c_{-d}\ge b_d+b_{-d}=E_0=b_n(B).
\]
Equation~\eqref{eq:count-base-window} gives
$A_0=b_n(B)-(-1)^n\chi(B)$.
Hence \eqref{eq:count-band-bound} implies
\[
 r\ge A_0+c_d+c_{-d}
 \ge 2b_n(B)-(-1)^n\chi(B).
\]

Suppose now that $p>0$.
Let $Q$ and $I$ denote the numbers of good iterates
$\gamma_i^{k_i}$ satisfying
\[
 |\mu(\gamma_i^{k_i})|\in\{d-p,d+p\},
 \qquad
 d-p<|\mu(\gamma_i^{k_i})|<d+p,
\]
respectively.
By \eqref{eq:count-central-ranks} and $c_p=c_{-p}=0$,
\[
 \begin{aligned}
 Q&=c_{d-p}+c_{d+p}+c_{-d-p}+c_{-d+p},\\
 I&=\sum_{|q|<p}(c_{d+q}+c_{-d+q}-c_q-c_{-q}).
 \end{aligned}
\]
Thus \eqref{eq:count-band-bound} gives $r\ge A_p+Q+I$.
Using the Morse inequalities \eqref{eq:count-finite-morse},
\eqref{eq:count-large-band}, and \eqref{eq:count-hole-rank},
together with $E_{-p}=E_p$, we obtain
\begin{equation}\label{eq:count-endpoint-bound}
 Q\ge E_{-p}+E_p=2b_{n-p}(B).
\end{equation}

Applying \eqref{eq:count-finite-morse} on
$[d-p+1,d+p-1]$ and $[-d-p+1,-d+p-1]$ gives
\[
 (-1)^{p+1}\sum_{|q|<p}(-1)^q(c_{d+q}+c_{-d+q})
 \ge
 (-1)^{p+1}\sum_{|q|<p}(-1)^q(b_{d+q}+b_{-d+q}).
\]
Since $c_{\pm p}=b_{\pm p}=0$,
\eqref{eq:count-hole-rank} also gives
\[
 \sum_{|q|<p}(-1)^q(c_q+c_{-q})
 =\sum_{|q|<p}(-1)^q(b_q+b_{-q}).
\]
Each summand defining $I$ is nonnegative by
\eqref{eq:count-central-ranks}. Consequently,
\begin{equation}\label{eq:count-interior-comparison}
 \begin{aligned}
 I
 &\ge (-1)^{p+1}\sum_{|q|<p}(-1)^q
       (c_{d+q}+c_{-d+q}-c_q-c_{-q})\\
 &\ge (-1)^{p+1}\sum_{|q|<p}(-1)^q
       (b_{d+q}+b_{-d+q}-b_q-b_{-q}).
 \end{aligned}
\end{equation}
Using \eqref{eq:count-interior-comparison},
\eqref{eq:count-large-band},
\eqref{eq:nonhyperbolic-rank-splitting}, and
\eqref{eq:count-base-window}, together with
Poincar\'e duality, we obtain
\begin{equation}\label{eq:count-interior-bound}
 \begin{aligned}
 I
 &\ge (-1)^{p+1}\sum_{|q|<p}(-1)^q(E_q-b_q-b_{-q})\\
 &=(-1)^{p+1}\sum_{|q|<p}(-1)^q b_{n-q}(B)\\
 &=2b_{n-p}(B)
   -(-1)^p\sum_{|q|\le p}(-1)^q b_{n-q}(B)\\
 &=2b_{n-p}(B)-A_p-(-1)^{n+p}\chi(B).
 \end{aligned}
\end{equation}
Finally, \eqref{eq:count-band-bound},
\eqref{eq:count-endpoint-bound}, and
\eqref{eq:count-interior-bound} give
\[
 \begin{aligned}
 r&\ge A_p+Q+I\\
  &\ge 4b_{n-p}(B)-(-1)^{n+p}\chi(B).
 \end{aligned}
\]
This proves \eqref{eq:count-final-bound}.

\textit{Step 3.} We prove that at least $L_p$
of the iterates $\gamma_1^{k_1},\ldots,\gamma_r^{k_r}$ are good
and satisfy $\mu(\gamma_i^{k_i})\equiv p\pmod2$, where
\begin{equation}\label{eq:count-parity-bound}
 L_p=
 \begin{cases}
  \displaystyle\frac{(-1)^n\chi(B)+b_n(B)}{2},
  &p=0,\\[6pt]
  \displaystyle\left\lceil
    \frac{(-1)^{n+p}\chi(B)}{2}
  \right\rceil+b_{n-p}(B),
  &p>0.
 \end{cases}
\end{equation}
In particular, $r\ge L_p$.

Suppose first that $p=0$.
By \eqref{ass:count-fci}, \eqref{eq:count-hole-rank}, and
\eqref{eq:count-base-window},
we have $c_0=b_0=0$, $E_0=b_n(B)$, and
$A_0=b_n(B)-(-1)^n\chi(B)$.
Equations~\eqref{eq:count-chain-shortened},
\eqref{eq:count-exact-euler}, and \eqref{eq:count-shortened-windows}
therefore give
\[
 \sum_{|q|\le d-1}(-1)^q(c_q-b_q)
 =-r_+^{e,-1}+r_+^{o,-1}+\frac{(-1)^n\chi(B)+b_n(B)}{2}.
\]
Applying \eqref{eq:count-finite-morse} on $[-d+1,d-1]$,
whose endpoints are odd, yields
\[
 r\ge r_+^{e,-1}\ge r_+^{e,-1}-r_+^{o,-1}
 \ge L_0.
\]
Poincar\'e duality implies that $L_0$ is an integer.

Suppose now that $p>0$.
The endpoints of $[-d-p+1,d+p-1]$ have parity $p-1$.
Applying \eqref{eq:count-finite-morse} on this interval
and using \eqref{eq:count-chain-shortened},
\eqref{eq:count-exact-euler}, \eqref{eq:count-shortened-windows},
and \eqref{eq:count-hole-rank}, we obtain
\[
 (-1)^p\bigl(r_+^{e,p-1}-r_+^{o,p-1}\bigr)
 \ge E_p-\frac{A_p}{2}.
\]
Applying the same calculation with $d,k_i,s$ replaced by
$d',k'_i,s'$, and using
Lemma~\ref{lem:index-recurrence}
\eqref{eq:count-recurrence-companion}, gives
\[
 (-1)^p\bigl(r_-^{e,p-1}-r_-^{o,p-1}\bigr)
 \ge E_p-\frac{A_p}{2}.
\]
Adding these inequalities and using
$E_p=b_{n-p}(B)$ from \eqref{eq:count-hole-rank}, we obtain
\begin{equation}\label{eq:count-same-outer}
 (-1)^p
 \bigl(r_+^{e,p-1}-r_+^{o,p-1}
       +r_-^{e,p-1}-r_-^{o,p-1}\bigr)
 \ge 2b_{n-p}(B)-A_p.
\end{equation}

Applying \eqref{eq:count-finite-morse} on
$[d-p,d+p]$ and $[-d-p,-d+p]$, and using
\eqref{eq:count-hole-rank}, gives
\begin{equation}\label{eq:count-same-central}
 \begin{aligned}
 &(-1)^p\sum_{|q|\le p}(-1)^q
       (c_{d+q}+c_{-d+q}-c_q-c_{-q})\\
 &\qquad\ge
 (-1)^p\sum_{|q|\le p}(-1)^q
       (b_{d+q}+b_{-d+q}-b_q-b_{-q})\\
 &\qquad=
 (-1)^p\sum_{|q|\le p}(-1)^q b_{n-q}(B)\\
 &\qquad=A_p+(-1)^{n+p}\chi(B),
 \end{aligned}
\end{equation}
where the equalities follow from
\eqref{eq:count-large-band},
\eqref{eq:nonhyperbolic-rank-splitting}, and
\eqref{eq:count-base-window}.

Let $r^e$ and $r^o$ denote the numbers of good iterates
among $\gamma_1^{k_1},\ldots,\gamma_r^{k_r}$
with even and odd index, respectively.
Recall that $Q$ counts the good iterates
$\gamma_i^{k_i}$ satisfying
$|\mu(\gamma_i^{k_i})|\in\{d-p,d+p\}$.
The counts in \eqref{eq:count-same-outer} include exactly
the good iterates $\gamma_i^{k_i}$ satisfying
\[
 |\mu(\gamma_i^{k_i})|\le d-p
 \qquad\text{or}\qquad
 |\mu(\gamma_i^{k_i})|\ge d+p.
\]
By \eqref{eq:count-central-ranks}, adding
\eqref{eq:count-same-outer} and
\eqref{eq:count-same-central} therefore counts every good
$\gamma_i^{k_i}$ once, except for the $Q$ iterates with
$|\mu(\gamma_i^{k_i})|\in\{d-p,d+p\}$, which are counted twice.
These $Q$ iterates satisfy $\mu(\gamma_i^{k_i})\equiv p\pmod2$, so
\begin{equation}\label{eq:count-same-signed}
 (-1)^{n+p}\chi(B)+2b_{n-p}(B)
 \le (-1)^p(r^e-r^o)+Q
 \le
 \begin{cases}
  2r^e-r^o,&p\text{ even},\\
  2r^o-r^e,&p\text{ odd}.
 \end{cases}
\end{equation}
Since $r^e,r^o\ge0$, \eqref{eq:count-same-signed}
and the definition \eqref{eq:count-parity-bound} imply that
at least $L_p$ of $\gamma_1^{k_1},\ldots,\gamma_r^{k_r}$ are good
and satisfy $\mu(\gamma_i^{k_i})\equiv p\pmod2$.

Thus $r\ge L_p$ for every $p\ge0$.

\textit{Step 4.} We now complete the proof
of \textup{(i)} and \textup{(ii)}.

For \textup{(i)}, suppose $p\not\equiv n\pmod2$ and $r=L_p$.
Then every $\gamma_i^{k_i}$ is good and satisfies $\mu(\gamma_i^{k_i})\equiv p\pmod2$.
All good iterates of $\gamma_i$ have the same index parity,
so $c_q=0$ for $q\not\equiv p\pmod2$.
Since $2\sigma c_B-n$ has parity $n$,
\eqref{eq:count-rank-nullity} and \eqref{eq:count-bw-pattern} give
\[
 0=c_{2\sigma c_B-n}
 \ge b_{2\sigma c_B-n}\ge b_0(B)=1,
\]
a contradiction. Thus $r\ge L_p+1$.
Combining this estimate with \eqref{eq:count-final-bound}
and $r_\alpha\ge r$ proves \eqref{eq:main-bound}.

For \textup{(ii)}, suppose $p\equiv n\pmod2$.
Assume now that $p=0$ or $n\ge3$, and suppose $r=L_p$.
Then every $\gamma_i^{k_i}$ is good and satisfies $\mu(\gamma_i^{k_i})\equiv n\pmod2$.
Since all good iterates of each $\gamma_i$ have the same
index parity, \eqref{eq:count-rank-nullity} gives
\[
 c_q=b_q=0\qquad(q\not\equiv n\pmod2).
\]
For odd $j$, \eqref{eq:count-bw-pattern} therefore gives
\[
 b_j(B)\le b_{j+2\sigma c_B-n}=0.
\]
Thus $H_{\mathrm{odd}}(B;\mathbb Q)=0$.

If $n$ is odd, then $c_0=0$; if $n$ is even, then
$c_1=c_{-1}=0$.
Applying \eqref{eq:count-final-bound} with $p=0$ when $n$
is odd and with $p=1$ when $n$ is even gives
$r\ge\chi(B)$, since $H_{\mathrm{odd}}(B;\mathbb Q)=0$.

By Poincar\'e duality,
\[
 \chi(B)=2\sum_{j=0}^{n-1}b_j(B)+b_n(B).
\]
If $p=0$, then $n$ is even, and
\eqref{eq:count-parity-bound} gives
\[
 \chi(B)-L_0
 =\sum_{j=0}^{n-1}b_j(B)
 \ge b_0(B)=1.
\]
If $p>0$ and $n\ge3$, then $b_0(B),b_2(B)\ge1$ and
$0,2<n$. Since $n-p<n$, at least one of the degrees
$0$ and $2$ differs from $n-p$. Hence
\[
 \sum_{j=0}^{n-1}b_j(B)\ge b_{n-p}(B)+1,
\]
also when $n-p<0$, since then $b_{n-p}(B)=0$.
Using \eqref{eq:count-parity-bound}, we obtain
\[
 \begin{aligned}
 \chi(B)-L_p
 &=\left\lfloor\frac{\chi(B)}{2}\right\rfloor-b_{n-p}(B)\\
 &=\sum_{j=0}^{n-1}b_j(B)
   +\left\lfloor\frac{b_n(B)}{2}\right\rfloor-b_{n-p}(B)\\
 &\ge1.
 \end{aligned}
\]
Thus $r\ge\chi(B)\ge L_p+1$, contradicting $r=L_p$.
We conclude that $r\ge L_p+1$ when $p=0$ or $n\ge3$.
Combining this inequality with \eqref{eq:count-final-bound}
and $r_\alpha\ge r$ proves
\eqref{eq:main-same-zero} for $p=0$,
\eqref{eq:main-same-positive} for $n\ge3$ and $p\ge2$,
and the last two bounds in \eqref{eq:main-odd-one} for $p=1$.

For $n=2$ and positive even $p$,
Lemma~\ref{lem:five-dimensional-index-zero}\textup{(i)} gives $c_0=0$.
The $p=0$ case of \eqref{eq:main-same-zero} therefore applies.
When $n=p=2$, \eqref{eq:count-final-bound} also gives
\eqref{eq:main-same-positive-lowdim}.

It remains to prove $r_\alpha\ge\chi(B)+2b_n(B)$
when $p=1$. Since $p\equiv n\pmod2$, $n$ is odd.
By \eqref{ass:count-fci}, $c_1=c_{-1}=0$.
Equation~\eqref{eq:count-bw-pattern} gives
\[
 b_0=\sum_{m\ge1}b_{n-2\sigma mc_B}(B).
\]
If $2c_B\le n$, every term vanishes by
Lemma~\ref{lem:betti-support}\textup{(i)}, because
$n-2\sigma mc_B$ is odd.
If $2c_B>n$, then
$n-2\sigma mc_B\notin[0,2n]$ for every $m\ge1$,
so every term again vanishes.
Thus $b_0=0$.
Since $c_1=c_{-1}=0$, the Morse inequality
\eqref{eq:count-finite-morse} on $[-1,1]$ gives
\[
 0\le c_0\le b_0=0.
\]
Applying \eqref{eq:count-final-bound} with $p=0$,
and using that $n$ is odd, we obtain
\[
 r_\alpha\ge r\ge\chi(B)+2b_n(B).
\]
This proves the first bound in \eqref{eq:main-odd-one}
and completes the proof of
Theorem~\ref{thm:main}\textup{(ii)}.
\end{proof}

\subsection{\texorpdfstring{Proof of Theorem~\ref{thm:main-nonhyperbolic}}%
  {Proof of Theorem \getrefnumber{thm:main-nonhyperbolic}}}
\label{subsec:nonhyperbolic}

\begin{proof}[Proof of Theorem~\ref{thm:main-nonhyperbolic}]
Let $\gamma_1,\ldots,\gamma_r$ and $d,k_1,\ldots,k_r$
be as in Section~\ref{subsec:main-proof}. By
Lemma~\ref{lem:index-recurrence}\eqref{eq:count-recurrence-mean},
$\mu(\gamma_i^{k_i})=\sigma_i d$ whenever the simple orbit
underlying $\gamma_i$ is hyperbolic. Thus a good iterate
$\gamma_i^{k_i}$ with $|\mu(\gamma_i^{k_i})|\ne d$ lies on a
non-hyperbolic simple orbit. In particular, every good
$\gamma_i^{k_i}$ of odd index lies on a non-hyperbolic simple orbit,
since $d$ is even.
The simple orbits underlying $\gamma_1,\ldots,\gamma_r$ are distinct.

Recall that $r^o$ and $r^e$ count the good iterates
among $\gamma_1^{k_1},\ldots,\gamma_r^{k_r}$ of odd and even
index, respectively.
We will use the implication
\begin{equation}\label{eq:nonhyperbolic-equality-parity}
 \left.
 \begin{aligned}
  r_\alpha^{\mathrm{nh}}&=r^o,\\
  r^e&=0
 \end{aligned}
 \right\}
 \quad\Longrightarrow\quad
 n\text{ is odd},\quad H_{\mathrm{odd}}(B;\mathbb Q)=0,
 \quad c_0=0.
\end{equation}
To prove \eqref{eq:nonhyperbolic-equality-parity}, assume that
$r_\alpha^{\mathrm{nh}}=r^o$ and $r^e=0$.
\par
Since $r^o=r_\alpha^{\mathrm{nh}}$ and the underlying simple
orbits are distinct, $\gamma_i^{k_i}$ is good and has odd index
whenever the simple orbit underlying $\gamma_i$ is non-hyperbolic.

Since every $\gamma_i$ is good, the definition of good and bad
iterates gives
\[
 \gamma_i^k\text{ is good}
 \quad\Longleftrightarrow\quad
 \mu(\gamma_i^k)\equiv\mu(\gamma_i)\pmod2
 \qquad(k\ge1).
\]
If $\gamma_i$ is non-hyperbolic, the good iterate
$\gamma_i^{k_i}$ has odd index, so $\mu(\gamma_i)$ is odd.
If $\gamma_i$ is hyperbolic, \eqref{eq:index-hyperbolic} and
the fact that $k_i$ is even give
\[
 \mu(\gamma_i^{k_i})=k_i\mu(\gamma_i)\in2\mathbb Z.
\]
Since $r^e=0$, $\gamma_i^{k_i}$ is bad, so $\mu(\gamma_i)$
is again odd. Thus every good iterate of every $\gamma_i$
has odd index, and $c_q=0$ for every even $q$.
The Morse inequality \eqref{eq:count-finite-morse} with $a=b=q$
then gives $b_q=0$ for every even $q$.

If $n$ is even, \eqref{eq:count-bw-pattern} gives
$b_{2\sigma c_B-n}\ge b_0(B)=1$ in an even degree, a contradiction.

If $n$ is odd, \eqref{eq:count-bw-pattern} gives
\[
 b_{j+2\sigma c_B-n}\ge b_j(B)
 \qquad(0\le j\le2n).
\]
For odd $j$, the degree $j+2\sigma c_B-n$ is even, so
$b_j(B)=0$. Hence $H_{\mathrm{odd}}(B;\mathbb Q)=0$ and
$b_n(B)=0$.
This proves \eqref{eq:nonhyperbolic-equality-parity}.

\textit{Step 1.} Let $O$ count the good iterates $\gamma_i^{k_i}$ with
$|\mu(\gamma_i^{k_i})|\notin[d-p,d+p]$. Then
\[
 O=r_+^{e,p}+r_-^{e,p}+r_+^{o,p}+r_-^{o,p}.
\]
Equations~\eqref{eq:count-central-ranks} and
\eqref{eq:count-outside-bound} give
\begin{equation}\label{eq:nonhyperbolic-band}
 \begin{split}
 r_\alpha^{\mathrm{nh}}&\ge O+
       \sum_{0<|q|\le p}(c_{d+q}+c_{-d+q}-c_q-c_{-q}),\\
 O&\ge\max\{0,A_p\}.
 \end{split}
\end{equation}

If $p=0$, \eqref{eq:count-base-window} gives
$A_0=b_n(B)-(-1)^n\chi(B)$,
so \eqref{eq:nonhyperbolic-band} proves \eqref{eq:main-nonhyperbolic}.

Assume now that $p>0$.
As in Section~\ref{subsec:main-proof}, let $Q$ count the good
iterates $\gamma_i^{k_i}$ satisfying
$|\mu(\gamma_i^{k_i})|\in\{d-p,d+p\}$.
Equations~\eqref{eq:count-central-ranks} and
\eqref{eq:count-endpoint-bound} give
\begin{equation}\label{eq:nonhyperbolic-boundary}
 Q=c_{d-p}+c_{d+p}+c_{-d-p}+c_{-d+p}
   \ge2E_p=2b_{n-p}(B).
\end{equation}
Together with \eqref{eq:nonhyperbolic-band}, this yields
$O+Q\ge2b_{n-p}(B)+A_p$.
Since the indices are integers, $r_-^{e,p-1}$ and $r_-^{o,p-1}$
count the good $\gamma_i^{k_i}$ with
$|\mu(\gamma_i^{k_i})|\le d-p$, while
$r_+^{e,p-1}$ and $r_+^{o,p-1}$ count those with
$|\mu(\gamma_i^{k_i})|\ge d+p$.
Since $d-p<d+p$, no iterate is counted in both ranges. Hence
\[
 O+Q=r_+^{e,p-1}+r_-^{e,p-1}+r_+^{o,p-1}+r_-^{o,p-1}.
\]
Equation~\eqref{eq:count-same-outer} therefore gives
\begin{equation}\label{eq:nonhyperbolic-inclusive-signed}
 \begin{aligned}
 O+Q
 &\ge(-1)^p\bigl(r_+^{e,p-1}+r_-^{e,p-1}
                    -r_+^{o,p-1}-r_-^{o,p-1}\bigr)\\
 &\ge2b_{n-p}(B)-A_p.
 \end{aligned}
\end{equation}
Combining \eqref{eq:nonhyperbolic-inclusive-signed}
with $O+Q\ge2b_{n-p}(B)+A_p$, we obtain
\begin{equation}\label{eq:nonhyperbolic-outer-boundary}
 O+Q\ge2b_{n-p}(B)+|A_p|.
\end{equation}
Every iterate counted by $O+Q$ satisfies
$|\mu(\gamma_i^{k_i})|\ne d$. Consequently,
$r_\alpha^{\mathrm{nh}}\ge O+Q\ge2b_{n-p}(B)+|A_p|$.

\textit{Step 2.} Suppose that $p>0$ is even.
By \eqref{eq:nonhyperbolic-rank-splitting},
\eqref{eq:count-base-window}, and \eqref{eq:count-hole-rank},
\begin{equation}\label{eq:nonhyperbolic-central-euler}
 \begin{split}
 \sum_{|q|\le p}(-1)^q(E_q-c_q-c_{-q})
 &=\sum_{|q|\le p}(-1)^q b_{n-q}(B)\\
 &=(-1)^n\chi(B)+A_p.
 \end{split}
\end{equation}
Since $c_p=c_{-p}=0$ and $E_{\pm p}=b_{n-p}(B)$,
removing the terms at $q=\pm p$ from
\eqref{eq:nonhyperbolic-central-euler} gives
\begin{equation}\label{eq:nonhyperbolic-outer-euler}
 \sum_{|q|<p}(-1)^q(E_q-c_q-c_{-q})
 =(-1)^n\chi(B)+A_p-2b_{n-p}(B).
\end{equation}

Apply
\eqref{eq:count-finite-morse} to the two intervals
$[d-p+1,d+p-1]$ and $[-d-p+1,-d+p-1]$, whose endpoints are odd.
By \eqref{eq:count-large-band} and \eqref{eq:nonhyperbolic-outer-euler},
\begin{equation}\label{eq:nonhyperbolic-even-middle}
 \begin{aligned}
 \sum_{|q|<p}(-1)^q(c_{d+q}+c_{-d+q}-c_q-c_{-q})
 &\le\sum_{|q|<p}(-1)^q(E_q-c_q-c_{-q})\\
 &=A_p+(-1)^n\chi(B)-2b_{n-p}(B).
 \end{aligned}
\end{equation}
By \eqref{eq:count-central-ranks}, the number of good
$\gamma_i^{k_i}$ with odd index and
$d-p<|\mu(\gamma_i^{k_i})|<d+p$ is at least
\[
 \max\{0,\,2b_{n-p}(B)-(-1)^n\chi(B)-A_p\}.
\]
These good iterates are not counted by $O+Q$, and their
underlying simple orbits are non-hyperbolic.
Adding this lower bound to \eqref{eq:nonhyperbolic-outer-boundary}
gives
\begin{equation}\label{eq:nonhyperbolic-even-general}
 r_\alpha^{\mathrm{nh}}\ge
 \max\left\{
 \begin{aligned}
  &2b_{n-p}(B)+|A_p|,\\
  &4b_{n-p}(B)+|A_p|-A_p-(-1)^n\chi(B)
 \end{aligned}
 \right\}.
\end{equation}

For even $n$, we also have
\begin{equation}\label{eq:nonhyperbolic-even-improved}
 r_\alpha^{\mathrm{nh}}\ge1-\chi(B).
\end{equation}
If $p\le n$, then $n-p$ is even, so $b_{n-p}(B)\ge1$,
as observed in the proof of Lemma~\ref{lem:betti-support}.
Since $|A_p|-A_p\ge0$, the second entry in
\eqref{eq:nonhyperbolic-even-general} is at least $4-\chi(B)$.
This proves \eqref{eq:nonhyperbolic-even-improved} when $p\le n$.

Suppose now that $p>n$. Then $A_p=b_{n-p}(B)=0$.
Equations~\eqref{eq:nonhyperbolic-even-middle} and
\eqref{eq:count-central-ranks} give at least $-\chi(B)$ good
iterates $\gamma_i^{k_i}$ with odd index and
$d-p<|\mu(\gamma_i^{k_i})|<d+p$. Hence
\[
 r_\alpha^{\mathrm{nh}}\ge O+Q-\chi(B).
\]
Suppose that $r_\alpha^{\mathrm{nh}}=-\chi(B)$.
Then $O=Q=0$, so every good $\gamma_i^{k_i}$ satisfies
$d-p<|\mu(\gamma_i^{k_i})|<d+p$.
By \eqref{eq:count-central-ranks}, the sum on the left-hand side
of \eqref{eq:nonhyperbolic-even-middle} is now $r^e-r^o$.
Thus
\[
 r_\alpha^{\mathrm{nh}}\ge r^o\ge r^e-\chi(B)
 =r^e+r_\alpha^{\mathrm{nh}}.
\]
It follows that $r^e=0$ and $r_\alpha^{\mathrm{nh}}=r^o$,
contradicting \eqref{eq:nonhyperbolic-equality-parity}, since $n$ is even.
Therefore $r_\alpha^{\mathrm{nh}}\ge1-\chi(B)$, proving
\eqref{eq:nonhyperbolic-even-improved}.
\par

\textit{Step 3.} Suppose that $p$ is odd.
Equation~\eqref{eq:count-same-central} gives
\begin{equation}\label{eq:nonhyperbolic-odd-middle}
 \sum_{|q|\le p}(-1)^q(c_{d+q}+c_{-d+q}-c_q-c_{-q})
 \le -A_p+(-1)^n\chi(B).
\end{equation}
Equations~\eqref{eq:nonhyperbolic-odd-middle} and
\eqref{eq:count-central-ranks} show that at least $A_p-(-1)^n\chi(B)$ good iterates
$\gamma_i^{k_i}$ have odd index and satisfy
$|\mu(\gamma_i^{k_i})|\in[d-p,d+p]$.
These good iterates are not counted by $O$, and their
underlying simple orbits are non-hyperbolic.
Using $O\ge\max\{0,A_p\}$ from \eqref{eq:nonhyperbolic-band},
we obtain
\begin{equation}\label{eq:nonhyperbolic-odd-general}
 r_\alpha^{\mathrm{nh}}
 \ge A_p+\max\{0,A_p\}-(-1)^n\chi(B).
\end{equation}

Suppose that either $n$ is even or $p<n$.
Equations~\eqref{eq:nonhyperbolic-odd-middle} and
\eqref{eq:count-central-ranks} give
\[
 r_\alpha^{\mathrm{nh}}\ge O+A_p-(-1)^n\chi(B).
\]
Suppose that $r_\alpha^{\mathrm{nh}}=A_p-(-1)^n\chi(B)$.
Then $O=0$, so every good $\gamma_i^{k_i}$ satisfies
$|\mu(\gamma_i^{k_i})|\in[d-p,d+p]$.
By \eqref{eq:count-central-ranks}, the sum on the left-hand side
of \eqref{eq:nonhyperbolic-odd-middle} is now $r^e-r^o$.
Hence
\[
 r_\alpha^{\mathrm{nh}}\ge r^o
 \ge r^e+A_p-(-1)^n\chi(B)
 =r^e+r_\alpha^{\mathrm{nh}}.
\]
It follows that $r^e=0$ and $r_\alpha^{\mathrm{nh}}=r^o$.
By \eqref{eq:nonhyperbolic-equality-parity}, $n$ is odd,
$H_{\mathrm{odd}}(B;\mathbb Q)=0$, and $c_0=0$.
In particular, $b_n(B)=0$, so the $p=0$ bound
\eqref{eq:main-nonhyperbolic} gives
$r_\alpha^{\mathrm{nh}}\ge\chi(B)$.
However, $p<n$ and the vanishing of odd homology imply
\[
 A_p=-2\sum_{j=0}^{n-p-1}b_j(B)\le-2,
\]
contradicting $r_\alpha^{\mathrm{nh}}=A_p+\chi(B)$.
Therefore, whenever $p$ is odd and either $n$ is even or $p<n$,
\begin{equation}\label{eq:nonhyperbolic-odd-improved}
 r_\alpha^{\mathrm{nh}}\ge A_p-(-1)^n\chi(B)+1.
\end{equation}
\par

We also prove
\begin{equation}\label{eq:main-nonhyperbolic-strict}
 r_\alpha^{\mathrm{nh}}\ge b_{n-p}(B)+1
       -\left\lceil\frac{(-1)^n\chi(B)}2\right\rceil.
\end{equation}
Since $p$ is odd,
\eqref{eq:count-same-signed} gives
\begin{equation}\label{eq:nonhyperbolic-strict-signed}
 2r^o-r^e\ge-(-1)^n\chi(B)+2b_{n-p}(B).
\end{equation}
Every iterate counted by $r^o$ lies on a non-hyperbolic
simple orbit. Since $r^e\ge0$, it follows that
\[
 r_\alpha^{\mathrm{nh}}\ge r^o\ge b_{n-p}(B)
       +\left\lceil-\frac{(-1)^n\chi(B)}2\right\rceil.
\]
If $\chi(B)$ is odd, then
\[
 \left\lceil-\frac{(-1)^n\chi(B)}2\right\rceil
 =1-\left\lceil\frac{(-1)^n\chi(B)}2\right\rceil,
\]
so \eqref{eq:main-nonhyperbolic-strict} follows.
Assume now that $\chi(B)$ is even. Since $p$ is odd,
\eqref{eq:count-parity-bound} gives
\[
 L_p=-\frac{(-1)^n\chi(B)}2+b_{n-p}(B)\in\mathbb Z.
\]
We have $r_\alpha^{\mathrm{nh}}\ge r^o\ge L_p$.
Suppose that $r_\alpha^{\mathrm{nh}}=L_p$.
Then $r^o=L_p$, and \eqref{eq:nonhyperbolic-strict-signed} gives
$r^e=0$.
By \eqref{eq:nonhyperbolic-equality-parity}, $n$ is odd,
$H_{\mathrm{odd}}(B;\mathbb Q)=0$, and $c_0=0$.
In particular, $b_n(B)=0$.
Since $c_0=0$, the $p=0$ bound
\eqref{eq:main-nonhyperbolic}, proved in Step~1, gives
$r_\alpha^{\mathrm{nh}}\ge\chi(B)$.
Poincar\'e duality now gives
\[
 \frac{\chi(B)}2=\sum_{0\le j<n}b_j(B).
\]
Using the Betti-number estimate from Step~4 of the proof of
Theorem~\ref{thm:main}, we obtain
\[
 \frac{\chi(B)}2\ge b_{n-p}(B)+1,
 \qquad \chi(B)>L_p.
\]
Hence $r_\alpha^{\mathrm{nh}}\ge\chi(B)>L_p$, a contradiction.
Thus $r_\alpha^{\mathrm{nh}}\ge L_p+1$, proving
\eqref{eq:main-nonhyperbolic-strict}.

\textit{Step 4.} We now obtain the bounds in
Theorem~\ref{thm:main-nonhyperbolic}\textup{(i)--(iii)}.

For \textup{(i)}, suppose $n=2$.
If $p>0$ is even, Lemma~\ref{lem:five-dimensional-index-zero}\textup{(i)}
gives $c_0=0$. Applying \eqref{eq:main-nonhyperbolic} with $p=0$
yields $r_\alpha^{\mathrm{nh}}\ge b_2(B)-\chi(B)$.
Since $A_p=0$, this estimate and
\eqref{eq:nonhyperbolic-even-general} give
\eqref{eq:main-nonhyperbolic-five-even}.

If $p$ is odd, Lemma~\ref{lem:five-dimensional-index-zero}\textup{(ii)}
gives $c_1=c_{-1}=0$.
Apply the estimate \eqref{eq:nonhyperbolic-outer-boundary}
with $1$ in place of $p$.
Since $A_1=2b_0(B)=2$, we obtain
\[
 r_\alpha^{\mathrm{nh}}\ge2b_1(B)+|A_1|=2b_1(B)+2,
\]
which proves \eqref{eq:main-nonhyperbolic-five-odd}.

For \textup{(ii)}, suppose $n\ge4$ is even.
When $p$ is even, \eqref{eq:nonhyperbolic-even-general} and
\eqref{eq:nonhyperbolic-even-improved} give
\eqref{eq:main-nonhyperbolic-even}.
When $p$ is odd, combine
\eqref{eq:nonhyperbolic-outer-boundary},
\eqref{eq:nonhyperbolic-odd-general},
\eqref{eq:nonhyperbolic-odd-improved}, and
\eqref{eq:main-nonhyperbolic-strict} to obtain
\eqref{eq:main-nonhyperbolic-same}.

For \textup{(iii)}, suppose $n\ge3$ is odd.
If $p$ is even, \eqref{eq:nonhyperbolic-even-general} gives
\eqref{eq:main-nonhyperbolic-odd-even}.
If $p<n$ is odd, \eqref{eq:nonhyperbolic-outer-boundary},
\eqref{eq:nonhyperbolic-odd-general},
\eqref{eq:nonhyperbolic-odd-improved}, and
\eqref{eq:main-nonhyperbolic-strict} give
\eqref{eq:main-nonhyperbolic-odd-low}, using
$-\lceil-\chi(B)/2\rceil=\lfloor\chi(B)/2\rfloor$.
If $p\ge n$ is odd, then $A_p=0$, and
\eqref{eq:nonhyperbolic-outer-boundary},
\eqref{eq:nonhyperbolic-odd-general}, and
\eqref{eq:main-nonhyperbolic-strict} give
\eqref{eq:main-nonhyperbolic-odd-high}.
\par
\end{proof}

\subsection{\texorpdfstring{Proof of Theorem~\ref{thm:zero-chern}}%
  {Proof of Theorem \getrefnumber{thm:zero-chern}}}
\label{subsec:zero-chern-proof}
\label{subsec:nonzero-area-proof}
\label{subsec:aspherical-proof}

Let $\alpha$ be a nondegenerate contact form supporting $\xi$.
Since $c_B=0$, we have $c_1(TB)|_{\pi_2(B)}=0$.
The identification $\xi\simeq\pi^*TB$ then gives
$c_1(\xi)|_{\pi_2(M)}=0$.

For \textup{(i)}, suppose $\omega|_{\pi_2(B)}=0$.
We use linearized contact homology in the fiber classes, following
the augmentation construction in \cite[Section~6.1]{GGM18}.
The homotopy exact sequence of the circle bundle contains
\begin{equation}\label{eq:aspherical-homotopy}
 0\longrightarrow\pi_2(M)\xrightarrow{\pi_*}\pi_2(B)
 \xrightarrow{-\langle[\omega],\cdot\rangle}\mathbb Z
 \longrightarrow\pi_1(M)\longrightarrow\pi_1(B)
 \longrightarrow1.
\end{equation}
Let $\mathfrak f\in\pi_1(M)$ denote the positively oriented fiber class.
Since $-\langle[\omega],\cdot\rangle=0$,
\eqref{eq:aspherical-homotopy} gives the central extension
\begin{equation}\label{eq:aspherical-central-extension}
 1\longrightarrow\langle \mathfrak{f}\rangle\cong\mathbb Z
 \longrightarrow\pi_1(M)\xrightarrow{\pi_*}\pi_1(B)
 \longrightarrow1.
\end{equation}
Centrality follows from the principal circle-bundle structure.
In particular, the conjugacy class of $\mathfrak{f}^k$ is the singleton
$\{\mathfrak{f}^k\}$, and $\mathfrak{f}$ has infinite order.

For a closed orbit $z$ with $\pi\circ z$ contractible in $B$,
choose a capping disk for $\pi\circ z$ and a symplectic
trivialization of $TB$ over that disk. The identification
$\xi\simeq\pi^*TB$ induces a trivialization along $z$ and hence
indices $\mu(z)$ and $\widehat\mu(z)$.
Since $c_1(TB)|_{\pi_2(B)}=0$, the indices are independent
of the chosen disk. When $z$ is contractible in $M$, the indices
$\mu(z)$ and $\widehat\mu(z)$ agree with those defined in
Section~\ref{subsec:index-iteration}.
For $z^\ell$, use the iterated trivialization; then
\eqref{eq:index-homogeneity} gives
$\widehat\mu(z^\ell)=\ell\widehat\mu(z)$.

\textit{Step 1.} We construct reference forms $\beta_j$ and compute
their linearized contact homology in the classes $\mathfrak f^k$.

Let $\beta$ and $f$ be as in
Section~\ref{subsec:finite-action-perturbations}.
After multiplying $\alpha$ by a positive constant, assume
$\alpha>2\beta$.
The construction in the proof of Proposition~\ref{prop:reference-sequence},
applied to all iterates of the simple fiber orbits in the
full contact algebras, gives nondegenerate forms
\[
 \alpha>\beta_1>\beta_2>\cdots>\beta
\]
and increasing positive integers $N_j\to\infty$.
Set $L_j=N_j+\tfrac12$.
By Lemma~\ref{lem:local-bw-finite} and
Corollary~\ref{cor:finite-orbit-persistence}, the closed orbits
of $\beta_j$ with action less than $L_j$ are precisely
$\gamma_{x,\beta_j}^{\ell}$, where $x\in\operatorname{Crit}(f)$
and $1\le\ell\le N_j$, and all are good.
Set
\[
 A_j=\mathcal C(\beta_j)^{<L_j},\qquad
 V_{\ell,j}=\operatorname{span}_{\mathbb Q}
 \{\gamma_{x,\beta_j}^{\ell}:x\in\operatorname{Crit}(f)\}.
\]
Each orbit $\gamma_{x,\beta_j}^{\ell}$ represents $\mathfrak f^\ell$.
For a genus-zero curve with positive end in
$\mathfrak f^\ell$ and negative ends in
$\mathfrak f^{\ell_1},\ldots,\mathfrak f^{\ell_s}$,
the relation between the boundary loops gives
\[
 \mathfrak f^\ell=\mathfrak f^{\ell_1+\cdots+\ell_s},
 \qquad \ell=\ell_1+\cdots+\ell_s,
\]
because $\mathfrak f$ is central and has infinite order.
Since $\ell\ge1$, no constant term occurs in the differentials
on $A_j$ or in the cobordism maps between the reference algebras.
Thus setting every generator equal to zero defines augmentations
$\varepsilon_j:A_j\to\mathbb Q$ compatible with these maps.

The index calculation in the proof of Lemma~\ref{lem:local-bw-finite}
gives
\[
 \mu(\gamma_{x,\beta_j}^{\ell})=-n+\operatorname{ind}_f(x).
\]
The cylinder calculations in the proofs of
Proposition~\ref{prop:local-bw-counts} and
Lemma~\ref{lem:local-bw-perturbation-persistence} apply to
each $1\le\ell\le N_j$.
Indeed, a limiting cylinder projects to a holomorphic sphere in $B$.
Since $\omega|_{\pi_2(B)}=0$, this sphere has zero area and is constant.
No plane can occur, since the asymptotic orbits with action less
than $L_j$ are noncontractible in $M$.
Only cylinders contribute to the linearized differential
$\partial_{\varepsilon_j}$, because $\varepsilon_j$ vanishes
on every generator. Dividing the cylinder counts in
\eqref{eq:local-cylinder-count} by $\ell$ gives
\[
 \partial_{\varepsilon_j}\gamma_{x,\beta_j}^{\ell}
 =\sum_{\operatorname{ind}_f(y)=\operatorname{ind}_f(x)-1}
   n_f(x,y)\,\gamma_{y,\beta_j}^{\ell}.
\]
Thus the linearized complex in $\mathfrak f^\ell$, with
Conley--Zehnder grading, is the Morse complex of $f$ shifted by $-n$.

\textit{Step 2.} We construct $\varepsilon_\alpha$ and compute
$\LCH_q^{\mathfrak f^k}(M,\alpha;\varepsilon_\alpha)$.

The construction in the proof of
Lemma~\ref{lem:close-reference-triangular} gives transition
maps $f_j:A_j\to A_{j+1}$ whose linear terms define isomorphisms
$V_{\ell,j}\to V_{\ell,j+1}$ for $\ell\le N_j$.
For a generator in $V_{\ell,j}$, every nonlinear term
of its image under $f_j$ involves only orbits of covering
multiplicity less than $\ell$. Hence $f_j$ identifies
$A_j$ with $\mathcal C(\beta_{j+1})^{<L_j}$.
The map $Qf_j$ therefore induces an isomorphism on homology
in each class $\mathfrak f^\ell$ with $\ell\le N_j$.
Let
\[
 A_\infty=\varinjlim_j A_j,\qquad
 \varepsilon_\infty:A_\infty\longrightarrow\mathbb Q
\]
be the direct limit algebra and the augmentation induced by
the $\varepsilon_j$.
The cobordism maps from $\alpha$ to $\beta_j$ induce
\[
 \Psi:\mathcal C(\alpha)\longrightarrow A_\infty,
 \qquad \varepsilon_\alpha=\varepsilon_\infty\circ\Psi.
\]
As in the proof of Proposition~\ref{prop:found-colimit-qiso},
the action bounds and the quasi-isomorphisms between the full
contact algebras imply that $\Psi$ is a quasi-isomorphism.
Since $f_j(A_j)=\mathcal C(\beta_{j+1})^{<L_j}$,
we extend the Sullivan orderings by adjoining the new generators
in increasing order of action.
Thus $A_\infty$ is Sullivan, and
Lemma~\ref{lem:found-forward-sullivan} gives a quasi-isomorphism
\[
 Q\Psi:Q(\mathcal C(\alpha),\varepsilon_\alpha)
 \longrightarrow\varinjlim_jQ(A_j,\varepsilon_j).
\]

A constant term in a cobordism map can occur only when the
positive orbit is contractible in $M$.
Consequently, $\varepsilon_\alpha$ vanishes on every noncontractible
generator and factors through $\mathcal B(\alpha)$.
Since $\mathcal B(\beta_j)^{<L_j}=\mathbb Q$,
Lemma~\ref{lem:found-contractible-quotient} implies that the
augmentation of $\mathcal B(\alpha)$ preserves the integer grading.
By Lemma~\ref{lem:free-homotopy-splitting}, $Q\Psi$ preserves
each free homotopy class $\mathfrak f^k$.
On each class $\mathfrak f^k$, the map $Q\Psi$ also preserves
the grading defined by $\mu$.
Indeed, projecting a curve contributing to a linearized differential
or a linearized cobordism map to $B$ and capping its projected ends gives a sphere
of zero first Chern number.
A negative asymptotic orbit with nonzero augmentation has contact
degree zero. The index formula therefore gives degree $-1$ for
the linearized differentials and degree zero for the linearized
cobordism maps.

For $k\ge1$ and $N_j\ge k$, Step~1 gives
$H_{q+n}(B;\mathbb Q)$ in Conley--Zehnder degree $q$ and class
$\mathfrak f^k$. Since each $Qf_j$ induces an isomorphism on
these homology groups, $Q\Psi$ gives
\begin{equation}\label{eq:aspherical-fiber-homology}
 \LCH_q^{\mathfrak f^k}(M,\alpha;\varepsilon_\alpha)
 \cong H_{q+n}(B;\mathbb Q)
 \qquad(k\ge1,\ q\in\mathbb Z).
\end{equation}
The correspondence of Reeb orbits under constant rescaling
identifies the linearized complexes in each $\mathfrak f^k$, so
\eqref{eq:aspherical-fiber-homology} also holds for the original $\alpha$.
In particular,
\begin{equation}\label{eq:aspherical-top-class}
 \LCH_n^{\mathfrak f^k}(M,\alpha;\varepsilon_\alpha)
 \cong H_{2n}(B;\mathbb Q)\cong\mathbb Q
 \qquad(k\ge1).
\end{equation}

\textit{Step 3.} Suppose, for a contradiction, that $\alpha$ has
only finitely many geometrically distinct simple closed Reeb
orbits $x_1,\ldots,x_r$.
Choose $g_i\in\pi_1(M)$ representing the free homotopy class of $x_i$,
and set
\[
 I=\{i\in\{1,\ldots,r\}:\pi_*(g_i)\text{ has finite order}\}.
\]
For $i\in I$, let $d_i$ be the order of $\pi_*(g_i)$.
By \eqref{eq:aspherical-central-extension}, there is a unique
$a_i\in\mathbb Z$ such that

\begin{equation}\label{eq:aspherical-root-data}
 g_i^{d_i}=\mathfrak{f}^{a_i}.
\end{equation}
If $x_i^m$ represents $\mathfrak f^k$, then
$g_i^m=\mathfrak f^k$.
Applying $\pi_*$ gives $i\in I$ and $m=d_i\ell$ for some
$\ell\ge1$. Since $\mathfrak f$ has infinite order,
\eqref{eq:aspherical-root-data} gives

\begin{equation}\label{eq:aspherical-root-exponents}
 a_i\ell=k.
\end{equation}
For any prime $k$ satisfying
\begin{equation}\label{eq:aspherical-prime-bound}
 k>\max\bigl(\{1\}\cup\{|a_i|:i\in I\}\bigr),
\end{equation}
equation~\eqref{eq:aspherical-root-exponents} forces $a_i=1$
and $\ell=k$.
Set
\[
 J=\{i\in I:a_i=1\},\qquad y_i=x_i^{d_i}\quad(i\in J).
\]
For every prime $k$ satisfying
\eqref{eq:aspherical-prime-bound}, the closed orbits in
$\mathfrak f^k$ are precisely $y_i^k$, $i\in J$.
Each $y_i$ represents $\mathfrak f$; it need not be simple.

By \eqref{eq:count-mean-estimate},
\[
 \mu(y_i^k)=n
 \quad\Longrightarrow\quad
 |n-k\widehat\mu(y_i)|<n
 \quad\Longrightarrow\quad
 0<k\widehat\mu(y_i)<2n.
\]
For each $i\in J$, the equality $\mu(y_i^k)=n$ holds
for only finitely many $k\ge1$.
Choose a prime $k$ satisfying \eqref{eq:aspherical-prime-bound}
such that $\mu(y_i^k)\ne n$ for every $i\in J$.
There is therefore no closed orbit in $\mathfrak f^k$ of index $n$, so
\[
 \LCC_n^{\mathfrak f^k}(M,\alpha;\varepsilon_\alpha)=0,
 \qquad
 \LCH_n^{\mathfrak f^k}(M,\alpha;\varepsilon_\alpha)=0,
\]
contradicting \eqref{eq:aspherical-top-class}.

For \textup{(ii)}, suppose $\omega|_{\pi_2(B)}\ne0$.
By hypothesis, $H_{\mathrm{odd}}(B;\mathbb Q)=0$, and hence
$b_j(B)b_{j-1}(B)=0$ for every $j\in\mathbb Z$.
Since $\sigma=0$, Theorem~\ref{thm:homology} gives an augmentation
$\varepsilon_\alpha$ for which \eqref{eq:BWpattern} becomes
\begin{equation}\label{eq:zero-chern-homology}
 \LCH^0_k(M,\alpha;\varepsilon_\alpha)
 \cong\bigoplus_{m\ge1}H_{k+n}(B;\mathbb Q).
\end{equation}
Since $H_0(B;\mathbb Q)\cong\mathbb Q$,
\eqref{eq:zero-chern-homology} gives
\begin{equation}\label{eq:zero-chern-infinite-rank}
 \LCH^0_{-n}(M,\alpha;\varepsilon_\alpha)
 \cong\bigoplus_{m\ge1}H_0(B;\mathbb Q)
 \cong\bigoplus_{m\ge1}\mathbb Q.
\end{equation}
Suppose, for a contradiction, that only finitely many simple closed
Reeb orbits have a contractible iterate. Let $\gamma_1,\ldots,\gamma_r$
be their least contractible iterates. Every contractible closed orbit
is then of the form $\gamma_i^k$ for some $1\le i\le r$ and $k\ge1$.
By \eqref{eq:count-mean-estimate},
\[
 \mu(\gamma_i^k)=-n
 \quad\Longrightarrow\quad
 -2n<k\widehat\mu(\gamma_i)<0.
\]
For each $i$, the equality $\mu(\gamma_i^k)=-n$ holds
for only finitely many $k\ge1$.
Thus only finitely many good contractible orbits have index $-n$, and
\[
 \dim\LCH^0_{-n}(M,\alpha;\varepsilon_\alpha)
 \le\dim\LCC^0_{-n}(M,\alpha;\varepsilon_\alpha)<\infty,
\]
contradicting \eqref{eq:zero-chern-infinite-rank}.

\subsection{\texorpdfstring{Proof of Theorem~\ref{thm:main-orbifold}}%
  {Proof of Theorem \getrefnumber{thm:main-orbifold}}}
\label{sec:homotopy-orbifold}

\begin{lemma}\label{lem:orbifold-topology}
Let $\pi:M^{2n+1}\to(B^{2n},\omega)$, $n\ge1$, be a
prequantization $S^1$-orbibundle over a closed connected symplectic
orbifold, with $M$ smooth.
Suppose that $B$ admits a Hamiltonian circle action with isolated
fixed points, and let $H:B\to\mathbb R$ be its Hamiltonian.
Then
\[
 H_{\mathrm{odd}}(B;\mathbb Q)=0,
 \qquad
 \#\operatorname{Crit}(H)=b(B),
\]
and $\pi_1(M)$ is finite cyclic.
In particular, every loop in $M$ has a contractible iterate.
\end{lemma}

\begin{proof}
By orbifold Morse theory \cite[Lemma~5.3 and Remark~5.4]{LeT97}, the Hamiltonian $H$
is a perfect Morse function over $\mathbb Q$, and all its
critical points have even index. Thus
\[
 b_j(B)
 =\#\{x\in\operatorname{Crit}(H):\operatorname{ind}_xH=j\};
\]
see also \cite[Section~8, case~\textup{(PO)}]{AbM26}.
It follows that
\[
 H_{\mathrm{odd}}(B;\mathbb Q)=0,
 \qquad
 \#\operatorname{Crit}(H)=b(B).
\]
Since $b_0(B)=1$, exactly one critical point of $H$ has index zero.

We next show that $\pi_1(M)$ is cyclic.
The function $H\circ\pi$ is Morse--Bott, with critical circles
$\pi^{-1}(x)$, $x\in\operatorname{Crit}(H)$, of index
$\operatorname{ind}_xH$.
Indeed, the bundle has the local form
$(S^1\times U)/G\to U/G$, where $G$ is finite and acts freely
on $S^1\times U$.
The lift of $H\circ\pi$ to $S^1\times U$ depends only on the
$U$-coordinate and equals the local lift of $H$.

Perturb $H\circ\pi$ near each critical circle using a Morse
function on the circle with one minimum and one maximum.
For each $x\in\operatorname{Crit}(H)$, the perturbed function
has two critical points of indices
$\operatorname{ind}_xH$ and $\operatorname{ind}_xH+1$.
Since the indices of $H$ are even and exactly one is zero,
the resulting handle decomposition of $M$ has one $0$-handle
and one $1$-handle. Attaching the remaining handles only
imposes relations on the fundamental group, so $\pi_1(M)$
is a quotient of $\mathbb Z$.

To prove that $\pi_1(M)$ is finite, it remains to show that
$H^1(M;\mathbb R)=0$.
Fix a connection form $\beta$ with $d\beta=\pi^*\omega$,
and let $\eta$ be a closed one-form on $M$.
Averaging over the circle action generated by $R_\beta$
preserves $[\eta]$, so we may assume that $\eta$ is invariant.
Cartan's formula gives
\[
 d\bigl(\eta(R_\beta)\bigr)=0.
\]
Since $M$ is connected, $\eta(R_\beta)=c$ for some constant $c$.
The form $\eta-c\beta$ is invariant and vanishes on $R_\beta$,
so there is an orbifold one-form $\theta$ on $B$ such that
\[
 \eta-c\beta=\pi^*\theta,
 \qquad
 d\theta=-c\omega.
\]
By Stokes' theorem,
\[
 0=\int_B d(\theta\wedge\omega^{n-1})
   =-c\int_B\omega^n.
\]
Since $\int_B\omega^n>0$, we have $c=0$. Hence
$d\theta=0$ and $\eta=\pi^*\theta$.
The orbifold de Rham theorem and $H^1(B;\mathbb R)=0$
imply that $\theta$ is exact, so $\eta$ is exact as well.
Therefore $H^1(M;\mathbb R)=0$, and the cyclic group
$\pi_1(M)$ is finite.
\end{proof}

\begin{proof}[Proof of Theorem~\ref{thm:main-orbifold}]
We may assume that $r_\alpha<\infty$.

\textit{Step 1.} We construct $\widetilde\beta$ and compute
$\LCH_*^0(M,\alpha;\varepsilon_\alpha)$.
Let $H:B\to\mathbb R$ generate the Hamiltonian circle action,
so $\operatorname{Crit}(H)=B^{S^1}$.
By Lemma~\ref{lem:orbifold-topology}, $\pi_1(M)$ is finite cyclic
and $\#B^{S^1}=b(B)$.
The construction in the proof of \cite[Proposition~5]{AbM26}
gives arbitrarily small nondegenerate perturbations
$\widetilde\beta$ of $\beta$ whose simple closed Reeb orbits
are precisely the fibers over $B^{S^1}$ and all of whose closed
orbits are elliptic. For each $x\in B^{S^1}$, let
$\gamma_{x,\widetilde\beta}$ denote the least contractible
closed Reeb orbit of $\widetilde\beta$ along the fiber over $x$.
The orbits $\gamma_{x,\beta}$ and $\gamma_{x,\widetilde\beta}$
differ only by reparametrization, so we use the same capping
disk to compute their mean indices.
Since $B^{S^1}$ is finite, continuity of the mean index
(Section~\ref{subsec:index-iteration}) and the hypothesis
$\widehat\mu(\gamma_{x,\beta})\ne0$ allow us to choose
$\widetilde\beta$ sufficiently close to $\beta$ that
\[
 \widehat\mu(\gamma_{x,\widetilde\beta})\ne0
 \qquad(x\in B^{S^1}).
\]
Every contractible closed orbit of $\widetilde\beta$ is an
iterate of some $\gamma_{x,\widetilde\beta}$ and hence has
nonzero mean index by \eqref{eq:index-homogeneity}.

Rescale $\widetilde\beta$ by a positive constant so that
$\alpha>\widetilde\beta$, without changing its orbit images or indices.
Lemmas~\ref{lem:found-action} and~\ref{lem:found-contractible-quotient}
give an integer-graded quasi-isomorphism
\[
 \bar F:\mathcal B(\alpha)\longrightarrow\mathcal B(\widetilde\beta).
\]
Every contractible closed orbit of $\widetilde\beta$ is good and
has index congruent to $n$ modulo $2$ by
\eqref{eq:index-elliptic-parity}.
Thus every generator of $\mathcal B(\widetilde\beta)$ has even
degree $\mu+n-2$, and its differential is zero.
Sending every generator to zero defines a graded augmentation
$\bar\varepsilon_{\widetilde\beta}$. Set
$\bar\varepsilon_\alpha=\bar\varepsilon_{\widetilde\beta}\bar F$,
and let $\varepsilon_\alpha$ and $\varepsilon_{\widetilde\beta}$
denote the pullbacks to the full contact algebras, as in
Section~\ref{subsec:contact-algebras}.

The contractible quotients are Sullivan by
Lemma~\ref{lem:found-contractible-quotient}, so
Lemma~\ref{lem:found-forward-sullivan} implies that
$Q\bar F$ is a quasi-isomorphism. Consequently,
\begin{equation}\label{eq:elliptic-reference-homology}
 \begin{aligned}
 \LCH_q^0(M,\alpha;\varepsilon_\alpha)
 &\cong\LCH_q^0(M,\widetilde\beta;\varepsilon_{\widetilde\beta})\\
 &\cong\bigoplus_{\substack{x\in B^{S^1},\ k\ge1\\
                  \mu(\gamma_{x,\widetilde\beta}^k)=q}}\mathbb Q.
 \end{aligned}
\end{equation}
We use the notation of Section~\ref{subsec:morse-mean-euler}:
\[
 b_q=\dim_{\mathbb Q}\LCH_q^0(M,\alpha;\varepsilon_\alpha),
 \qquad
 c_q=\dim_{\mathbb Q}\LCC_q^0(M,\alpha;\varepsilon_\alpha).
\]
Equation~\eqref{eq:elliptic-reference-homology} gives
$b_q=0$ for $q\not\equiv n\pmod2$, and the index hypothesis
on $\alpha$ gives $c_p=c_{-p}=0$.

\textit{Step 2.} Let $\gamma_1,\ldots,\gamma_r$ be all geometrically
distinct good least contractible closed Reeb orbits of $\alpha$.
The homology isomorphism \eqref{eq:elliptic-reference-homology}
and the resonance relations \eqref{eq:orbit-resonance} give
\begin{equation}\label{eq:orbifold-resonance}
 \sum_{i=1}^r\frac{\widehat\chi(\gamma_i)}{|\widehat\mu(\gamma_i)|}
 =\sum_{x\in B^{S^1}}
   \frac{\widehat\chi(\gamma_{x,\widetilde\beta})}
        {|\widehat\mu(\gamma_{x,\widetilde\beta})|}.
\end{equation}
Choose $0<\eta<1$ such that
\[
 \eta\left(
 \sum_{i=1}^r
   \left|\frac{\widehat\chi(\gamma_i)}{\widehat\mu(\gamma_i)}\right|
 +\sum_{x\in B^{S^1}}
   \left|\frac{\widehat\chi(\gamma_{x,\widetilde\beta})}
                  {\widehat\mu(\gamma_{x,\widetilde\beta})}\right|
 \right)<1.
\]
Choose $\ell_0\in\mathbb N$ such that
\[
 \begin{aligned}
 \ell_0|\widehat\mu(\gamma_i)|&>n+p+2
 &&(1\le i\le r),\\
 \ell_0|\widehat\mu(\gamma_{x,\widetilde\beta})|&>n+p+2
 &&(x\in B^{S^1}).
 \end{aligned}
\]
Apply Lemma~\ref{lem:index-recurrence} simultaneously to the
symplectic paths associated with $\gamma_1,\ldots,\gamma_r$
and $\gamma_{x,\widetilde\beta}$, $x\in B^{S^1}$, with
parameters $\eta,\ell_0$, and $N=2$.
Denote the common even integers by $d,d'$, the even exponents
for $\gamma_i$ by $k_i,k'_i$, and those for
$\gamma_{x,\widetilde\beta}$ by $\widetilde k_x,\widetilde k'_x$.
By Lemma~\ref{lem:index-recurrence}, we may choose these integers
so that
\[
 \begin{aligned}
 k_i,k'_i&>2\ell_0
 &&(1\le i\le r),\\
 \widetilde k_x,\widetilde k'_x&>2\ell_0
 &&(x\in B^{S^1}),\\
 d,d'&>n+p+2+|\mu(\gamma_i^\ell)|
 &&(1\le i\le r,\ 1\le\ell\le\ell_0),\\
 d,d'&>n+p+2+|\mu(\gamma_{x,\widetilde\beta}^\ell)|
 &&(x\in B^{S^1},\ 1\le\ell\le\ell_0).
 \end{aligned}
\]

By \eqref{eq:orbifold-resonance}, the choice of $\eta$, and
Lemma~\ref{lem:index-recurrence}\eqref{eq:count-recurrence-mean},
\[
 \left|
 \sum_{i=1}^r k_i\widehat\chi(\gamma_i)
 -\sum_{x\in B^{S^1}}\widetilde k_x
       \widehat\chi(\gamma_{x,\widetilde\beta})
 \right|<1.
\]
Since $k_i$ and $\widetilde k_x$ are even, both sums are
integers by \eqref{eq:count-period-two}. Hence
\[
 \sum_{i=1}^r k_i\widehat\chi(\gamma_i)
 =\sum_{x\in B^{S^1}}\widetilde k_x
       \widehat\chi(\gamma_{x,\widetilde\beta}).
\]
The same argument gives
\[
 \sum_{i=1}^r k'_i\widehat\chi(\gamma_i)
 =\sum_{x\in B^{S^1}}\widetilde k'_x
       \widehat\chi(\gamma_{x,\widetilde\beta}).
\]
Define
\[
 \begin{aligned}
 K(d)
 &=\sum_{|q|\le d+p}(-1)^q b_q
   -\sum_{i=1}^r k_i\widehat\chi(\gamma_i),\\
 K(d')
 &=\sum_{|q|\le d'+p}(-1)^q b_q
   -\sum_{i=1}^r k'_i\widehat\chi(\gamma_i).
 \end{aligned}
\]

\textit{Step 3.} We prove $r_\alpha\ge b(B)$.
Set $\sigma_i=\operatorname{sign}\widehat\mu(\gamma_i)$ and use
$r_\pm^{e,p},r_\pm^{o,p}$ for $\gamma_i^{k_i}$ as in
Section~\ref{subsec:main-proof}.
The proof of Lemma~\ref{lem:count-iterates} applies to these
contractible orbits. Combining \eqref{eq:count-chain-window}
with \eqref{eq:count-finite-morse} gives
\[
 \begin{aligned}
 (-1)^p(r_+^{o,p}-r_+^{e,p})
 &\ge(-1)^p\left(K(d)-\sum_{|q|\le p}(-1)^q c_q\right),\\
 (-1)^p(r_-^{o,p}-r_-^{e,p})
 &\ge(-1)^p\left(K(d')-\sum_{|q|\le p}(-1)^q c_q\right).
 \end{aligned}
\]
Here the second inequality uses the exponents $k'_i$ and the
companion identity in
Lemma~\ref{lem:index-recurrence}\eqref{eq:count-recurrence-companion}.
The counts with subscripts $+$ and $-$ concern the good
$\gamma_i^{k_i}$ with absolute index greater than $d+p$ and
less than $d-p$, respectively.
The simple orbits underlying the good iterates $\gamma_i^{k_i}$
are pairwise distinct.
Adding the lower bounds for these two counts and the count in
\eqref{eq:count-central-ranks} for $|q|\le p$ gives
\begin{equation}\label{eq:reference-target-window}
 \begin{aligned}
 r_\alpha
 &\ge(-1)^p\left(K(d)+K(d')
                   -2\sum_{|q|\le p}(-1)^q c_q\right)\\
 &\quad+\sum_{|q|\le p}
                 (c_{d+q}+c_{-d+q}-c_q-c_{-q}).
 \end{aligned}
\end{equation}

Apply the same counting identities to
$\gamma_{x,\widetilde\beta}^{\widetilde k_x}$, $x\in B^{S^1}$.
These iterates are all good and satisfy
$\mu(\gamma_{x,\widetilde\beta}^{\widetilde k_x})\equiv n\pmod2$,
while $n\not\equiv p\pmod2$ by hypothesis.
Their underlying simple orbits are precisely the $b(B)$
distinct simple closed orbits of $\widetilde\beta$.
Since the differential of
$\LCC_*^0(M,\widetilde\beta;\varepsilon_{\widetilde\beta})$
is zero, the Morse inequalities are equalities. Thus
\begin{equation}\label{eq:reference-orbit-count}
 \begin{aligned}
 b(B)
 &=(-1)^p\left(
     K(d)+K(d')-2\sum_{|q|\le p}(-1)^q b_q
   \right)\\
 &\quad+(-1)^{p+1}\sum_{|q|\le p}(-1)^q
       (b_{d+q}+b_{-d+q}-b_q-b_{-q})\\
 &=(-1)^p\bigl(K(d)+K(d')\bigr)
   +(-1)^{p+1}\sum_{|q|\le p}(-1)^q
       (b_{d+q}+b_{-d+q}).
 \end{aligned}
\end{equation}

If $p=0$, then $n$ is odd. Hence $b_d=b_{-d}=0$ and $c_0=0$.
Equations~\eqref{eq:reference-target-window} and
\eqref{eq:reference-orbit-count} give
\[
 r_\alpha\ge K(d)+K(d')=b(B).
\]

Suppose $p>0$. Every summand in the last sum of
\eqref{eq:reference-target-window} is nonnegative by
\eqref{eq:count-central-ranks}. The unsigned sum equals
$(-1)^{p+1}$ times its alternating sum plus twice the sum over
$q\equiv p\pmod2$. Since $c_p=c_{-p}=0$, the two summands
with $q=\pm p$ have sum
\[
 Q=c_{d-p}+c_{d+p}+c_{-d-p}+c_{-d+p}.
\]
Keeping these two summands and cancelling the
alternating sums of $c_q$ over $[-p,p]$ gives
\[
 r_\alpha\ge(-1)^p\bigl(K(d)+K(d')\bigr)
 +(-1)^{p+1}\sum_{|q|\le p}(-1)^q(c_{d+q}+c_{-d+q})+2Q.
\]
The contribution of $q=\pm p$ to
$(-1)^{p+1}\sum_{|q|\le p}(-1)^q(c_{d+q}+c_{-d+q})$
is $-Q$. Apply \eqref{eq:count-finite-morse} on
$[d-p+1,d+p-1]$ and $[-d-p+1,-d+p-1]$, whose endpoints have
parity $p-1$. We obtain
\begin{equation}\label{eq:orbifold-total-bound}
 \begin{aligned}
 r_\alpha
 &\ge(-1)^p\bigl(K(d)+K(d')\bigr)
 +(-1)^{p+1}\sum_{|q|<p}(-1)^q(b_{d+q}+b_{-d+q})+Q\\
 &=b(B)+Q.
 \end{aligned}
\end{equation}
The equality follows from \eqref{eq:reference-orbit-count},
since $b_{d-p}=b_{d+p}=b_{-d-p}=b_{-d+p}=0$ by parity.
This proves $r_\alpha\ge b(B)$.

\textit{Step 4.} Suppose $n$ is odd, so $p$ is even.
If the simple orbit underlying $\gamma_i$ is hyperbolic, then
Lemma~\ref{lem:index-recurrence}\eqref{eq:count-recurrence-mean}
gives $|\mu(\gamma_i^{k_i})|=d$.
Thus every good $\gamma_i^{k_i}$ satisfying $|\mu(\gamma_i^{k_i})|\notin[d-p,d+p]$
or counted in \eqref{eq:count-central-ranks} with $q\ne0$ is non-hyperbolic.
Consequently, \eqref{eq:reference-target-window} holds with
$r_\alpha^{\mathrm{nh}}$ on the left if the $q=0$ summand is
omitted from the last sum on the right.
For $p=0$, this gives
\[
 r_\alpha^{\mathrm{nh}}\ge K(d)+K(d')=b(B).
\]

For $p>0$, the nonnegativity of the counts in
\eqref{eq:count-central-ranks} and the evenness of $p$ give
\[
 \begin{aligned}
 &\sum_{0<|q|\le p}
       (c_{d+q}+c_{-d+q}-c_q-c_{-q})\\
 &\qquad\ge
 -\sum_{|q|\le p}(-1)^q
       (c_{d+q}+c_{-d+q}-c_q-c_{-q})+2Q.
 \end{aligned}
\]
Indeed, the difference between the left- and right-hand sides is
\[
 c_d+c_{-d}-2c_0
 +2\sum_{\substack{0<|q|<p\\q\ \mathrm{even}}}
       (c_{d+q}+c_{-d+q}-c_q-c_{-q}),
\]
which is nonnegative.
Consequently, \eqref{eq:reference-target-window}, with
$r_\alpha^{\mathrm{nh}}$ on the left and the $q=0$ count
omitted on the right, gives
\[
 r_\alpha^{\mathrm{nh}}
 \ge K(d)+K(d')
 -\sum_{|q|\le p}(-1)^q(c_{d+q}+c_{-d+q})+2Q.
\]
Applying the Morse inequalities on
$[d-p+1,d+p-1]$ and $[-d-p+1,-d+p-1]$, as in
\eqref{eq:orbifold-total-bound}, yields
\[
 r_\alpha^{\mathrm{nh}}\ge b(B)+Q\ge b(B).
\]
\end{proof}

\section*{Acknowledgements}
The author is partially supported by the National
Natural Science Foundation of China (Grant Nos.~12631007, 12271268,
and 12361141812) and the Natural Science Foundation of Tianjin
Municipality (Grant No.~25JCZDJC01030).

% Display labels use surname abbreviations and the final two year digits.
% For two authors, use the first two letters of the first surname and the initial of the second.
% Particles de/van are omitted; equal author/year labels receive a/b suffixes.

\end{document}